\documentclass[sn-jnl, 11pt]{sn-jnl}

\usepackage{graphicx}
\usepackage{amsfonts,amssymb,amsmath}
\usepackage{url}
\usepackage{caption}
\usepackage{subcaption}
\usepackage{graphicx}
\usepackage{placeins}
\usepackage{xcolor}
\usepackage{ulem}
\usepackage[utf8]{inputenc}

\graphicspath{{Figures/}}

\usepackage{xcolor}
\definecolor{revisionColor}{rgb}{0, 0.5, 0}
\usepackage{cite}

\jyear{2021}%

\theoremstyle{thmstyleone}%

\theoremstyle{thmstyletwo}%

\theoremstyle{thmstylethree}%

\begin{document}

\title[Locally refined quad meshing based on convolutional neural networks]{Locally refined quad meshing for linear elasticity problems based on convolutional neural networks}

\author*[1]{\fnm{Chiu Ling} \sur{Chan}}\email{chiu\_ling.chan@jku.at}

\author[2]{\fnm{Felix} \sur{Scholz}}\email{felix.scholz@tafsm.org}
\equalcont{These authors contributed equally to this work.}

\author[1,3]{\fnm{Thomas} \sur{Takacs}}\email{thomas.takacs@ricam.oeaw.ac.at}
\equalcont{These authors contributed equally to this work.}

\affil[1]{\orgdiv{Institute of Applied Geometry}, \orgname{Johannes Kepler University Linz}, \orgaddress{\street{Altenberger Str. 69}, \postcode{4040} \city{Linz}, \country{Austria}}}

\affil[2]{\orgdiv{Waseda Research Institute for Science and Engineering}, \orgname{Waseda University}, \orgaddress{\street{3-4-1 Okubo}, \city{Shinjuku}, \postcode{169-8555}, \state{Tokyo}, \country{Japan}}}

\affil[3]{\orgdiv{Johann Radon Institute for Computational and Applied Mathematics}, \orgname{Austrian Academy of Sciences}, \orgaddress{\street{Altenberger Str. 69}, \postcode{4040} \city{Linz}, \country{Austria}}}

\abstract{In this paper we propose a method to generate suitably refined finite element meshes using neural networks. 
	As a model problem we consider a linear elasticity problem on a planar domain (possibly with holes) having a polygonal boundary.
	We impose boundary conditions by fixing the position of a part of the boundary and applying a force on another part of the boundary. 
	The resulting displacement and distribution of stresses depend on the geometry of the domain and on the boundary conditions.
	When applying a standard Galerkin discretization using quadrilateral finite elements, one usually has to perform adaptive refinement to properly resolve maxima of the stress distribution. Such an adaptive scheme requires a local error estimator and a corresponding local refinement strategy. The overall costs of such a strategy are high. We propose to reduce the costs of obtaining a suitable discretization by training a neural network whose evaluation replaces this adaptive refinement procedure. {We set up a single network for} a large class of possible domains and boundary conditions and not on a single domain of interest. The computational domain and boundary conditions are interpreted as images, which are suitable inputs for  convolution neural networks. In our approach we use the U-net architecture and we {devise} training strategies by dividing the possible inputs into different categories based on their overall {geometric} complexity. {Thus, we compare different training strategies based on varying geometric complexity.} One of the advantages of the proposed approach is the interpretation of input and output as images, which do not depend on the underlying discretization scheme. {Another is the generalizability and geometric flexibility. The network can be applied to previously unseen geometries, even with different topology and level of detail. Thus, training can easily be extended to other classes of geometries.}}


\maketitle

\section{Introduction}\label{sec1}
   
In this work, we consider the problem of finding an optimal discretization for a linear elasticity problem for a given planar domain and boundary conditions. 
Previously, this problem has often been solved by adaptive refinement, which means that the partial differential equation (PDE) is first solved on an initial mesh, then the error is estimated using a local error estimator, and the mesh is subsequently refined in the regions of the domain where the error is large. This process is repeated until the obtained uniform error is below a given tolerance. Adaptive refinement is computationally costly, since at each step, a new mesh is generated and progressively larger linear systems need to be assembled and solved. The objective of this study is therefore to find an optimal mesh {a priori} without needing to first solve the PDE, thereby avoiding the computational costs involved in iterative mesh refinements.

For this purpose, we train a neural network that takes as input the geometry and boundary conditions, and predicts a relative mesh size field for linear elasticity problems. It is well known that for a complex geometry, certain areas need to be refined more in order to obtain acceptable accuracy, e.g. around re-entrant corners or near to the areas where fixed boundary conditions are prescribed. Using a machine learning approach for the prediction of mesh refinement has several advantages. First, one can avoid an iterative refinement scheme by obtaining a suitable initial guess directly from the initial model. Furthermore, neural networks (especially convolutional neural networks) allow training on several classes of simple configurations, while they can be evaluated on more complex domains. {Based on the output from the neural network, we construct a quadrilateral mesh with local refinement properties. Quadrilateral meshes are adopted because of their superior performance in engineering applications. Moreover, we intend to extend the method to isogeometric discretizations based on (curvilinear) quadrilateral patches. There is a rich literature on this topic, for example automatic generation of high quality quadrilateral meshes from given planar curves \cite{liang2012matching}. The method emphasizes matching interior and exterior boundaries and avoids distorted quads by introducing angle bounds.}

Machine learning and especially artificial neural networks (ANN) have been employed recently to geometric problems as well as to problems that arise when solving PDEs. For instance, machine learning techniques are employed in model order reduction~\cite{tang2019map123, li2019clustering}, and dynamic model decomposition \cite{Barros2021}. Specifically, convolutional neural networks (CNN) are often applied in computational problems~\cite{han2019deep,finol2019deep,Li2020}. The combination of techniques from computational mathematics and ANN has resulted in interesting contributions for both fields, since the theoretical results from classical approximation theory can be used to derive results on the approximation properties of ANN~\cite{1908.03833,he2019mgnet,he2019relu,Panghal2020,Chi2021}.

A straight-forward application of deep learning to the solution of PDEs consists of training a neural network to represent the solution of a single boundary value problem~\cite{SIRIGNANO20181339}.
Especially for very high-dimensional problems, this approach results in efficient methods, compared to classical Galerkin methods.  But ANN cannot only be used to describe the discretization spaces for solutions of a PDE, machine learning techniques have also been used in several ways to facilitate the efficient and accurate solution of PDEs using classical methods, such as FEM.
Deep learning has been used for the system matrix assembly in Galerkin methods~\cite{compmechanics} as well as for accelerating the solution process of the resulting system~\cite{Tompson:2017:AEF:3305890.3306035}.

Another very promising application of ANN in this context is mesh generation. 
In order to ensure an efficient approximation, a locally refined discretization space is often needed. So far, the efforts in this direction have focused on generating polygonal/polyhedral meshes for finite element discretizations~\cite{717791,MANEVITZ2005447,doi:10.1142/S021821300800431X, pfluger2010spatially}.
In \cite{chen2020developing} an ANN is trained to predict the mesh quality of a given finite element mesh. In our approach, we train an ANN to predict the optimal local mesh density for a given polygonal geometry that can have holes, and given boundary conditions in the form of one traction boundary and one fixed boundary. Moreover, we study in this paper the dependence of the output of our neural network on the complexity of both the computational domain and different training strategies. We represent all input data, i.e. the domain as well as the boundary conditions, as images, and we employ a convolutional neural network (CNN) that maps to a pixel-wise estimate for the optimal mesh density.

The use of images and convolutional neural networks has several advantages:
CNNs are able to detect local features in the data, which is especially well-suited to the problem of local refinement.
Moreover, it is possible to rescale the models in order to apply them to different input sizes. 
Finally, the use of images as our input data makes the method completely independent of the employed scheme for numerically approximating the solution to the PDE.
This means that the training data for our model can be produced by using finite elements based on different polygons as well as using isogeometric analysis~\cite{hughes2005isogeometric}. {Moreover, image-based engineering and science is a growing research area. Generating meshes from scanned images has wide applications in computation biology, medicine and materials science \cite{zhang2018geometric}. As an example, the construction of quadrilateral and hexahedral meshes from volumetric image data is discussed in \cite{zhang2006adaptive}. }

{In recent works, the concept of CNN has been generalized from planar image data to discrete manifolds and graphs~\cite{bronstein2017geometric, 10.1145/3326362}, resulting in geometric deep learning.
This makes it possible to train and evaluate neural networks on data sets that consist of discrete surfaces. 
These techniques might allow a future generalization of our method to the problem of local refinement on surfaces.
However for the processing of planar domains, as we do in the present paper, as well as of volumetric domains, image data and CNNs are more suitable.
In particular, the representation of the geometry as an image is more flexible and does not depend on a specific parameterization.}

The approach we present in this paper is similar to~\cite{zhang2020meshingnet, zhang2021meshingnet3d}, where a neural network is used to facilitate the generation of locally refined finite element meshes. The network predicts the distribution of the a posteriori error for a given set of geometry, material properties and boundary conditions, {which demonstrates the feasibility of neural networks for generating locally refined meshes.} We want to highlight here that the network architecture and the structure of input and output are important factors when generalizing the approach. {The approach developed in~\cite{zhang2020meshingnet, zhang2021meshingnet3d} is based on a fully connected neural network with a fixed structure, which encodes the geometry. Thus, the possible geometries are taken from a parameterized family of possible geometries. In our approach, we consider the domain and input parameters to be images and employ CNNs. Using this setup, we can extend the network to other geometries by incorporating them in the training process, while keeping the network structure fixed. We believe that this facilitates the extension to more complex domains.}

The remainder of this paper is organized as follows: 
In Section~\ref{sec:problem-formulation}, we describe the model problem and objective which we consider throughout the paper.
In Section~\ref{complexity}, we analyze and classify different measures for geometric complexity that occur in the problems that we consider.
We describe our methods for generating the training data in Section~\ref{datagen} and present the architecture of our neural network in Section~\ref{networkarchitecture}.
Finally, we present numerical experiments in Section~\ref{numericalexamples}, {discuss directions for further research on the use of artificial neural networks for generating locally refined meshes in Section~\ref{futurework} and conclude the paper in Section~\ref{conclusions}.}

\section{Problem formulation}
\label{sec:problem-formulation}
Given a computational domain and corresponding boundary conditions, we want to obtain a quadrilateral mesh over which we can define a finite element space to represent the solution of a linear elasticity problem. 
The goal is to find a mesh that yields a numerical solution of high quality without comprising the computational complexity, i.e., the uniform error should be below a given tolerance and the number of elements should be as small as possible. In order to do this, we train a neural network that takes as an input the geometry, the Dirichlet and Neumann boundary conditions and whose output predicts the relative local mesh resolution.

\subsection{The PDE model problem}\label{pdeProblem}

As a model problem, we consider the problem of linear elasticity on a planar domain $\Omega\subset\mathbb{R}^2$ with polygonal boundary $\partial\Omega$. The segments $\Gamma_D \subset \partial\Omega$ and $\Gamma_N\subset \partial\Omega$ are the Dirichlet and Neumann boundaries, respectively, where $\Gamma_D\cap\Gamma_N =\emptyset$. The governing equations are the equilibrium equation~\eqref{eq:equilibrium}, strain-displacement relation~\eqref{eq:strain-displacement}, and the constitutive law~\eqref{eq:constitutive}:
\begin{subequations}\label{eq:gov_elast}
	\begin{align}
		&\boldsymbol{-\nabla\cdot\sigma}=\boldsymbol{f}{,} \label{eq:equilibrium} \\
		&\boldsymbol{\epsilon}=\frac{1}{2}(\boldsymbol{\nabla u }+ \boldsymbol{\nabla} \boldsymbol{u}^T ){,}\label{eq:strain-displacement} \\
		&\boldsymbol{\sigma}=2\mu\boldsymbol{\epsilon}+\lambda(\boldsymbol{\nabla} \cdot \boldsymbol{u})\boldsymbol{I}.\label{eq:constitutive} 
	\end{align}
\end{subequations}
In the governing equations, $\boldsymbol{\sigma}$ and $\boldsymbol{\epsilon}$ denote the stress and the strain tensors respectively, while  $\boldsymbol{u}:\Omega \rightarrow \mathbb{R}^2$ is the displacement field and $\boldsymbol{f}:\Omega \rightarrow \mathbb{R}^2 $ is the body force. The constants $\lambda$ and $\mu$ are the Lam\'e parameters that satisfy
\begin{align*} 
	\lambda&=\frac{\nu E}{(1+\nu)(1-2\nu)},\\
	\mu&=\frac{E}{2(1+\nu)},
\end{align*}
where $E$ is Young's modulus and $\nu$ is Poisson's ratio. We moreover consider boundary conditions of the form
\begin{align*}
	\boldsymbol{u}&=\boldsymbol{u}_0 \quad \mbox{on} \quad  \Gamma_D  \subset \partial\Omega,\\
	\boldsymbol{\sigma}\boldsymbol{n} &= \boldsymbol{t} \quad \mbox{on}  \quad  \Gamma_N \subset \partial\Omega,
\end{align*}
where $\boldsymbol{u}_0$ are the prescribed displacements, $\boldsymbol{t}$ is the prescribed traction, and $\boldsymbol{n}$ is the outer unit normal vector to $\Gamma_N$. 
Then, after combining the equations in \eqref{eq:gov_elast}, multiplying by a test function $\boldsymbol{v}$ and applying Green's theorem, we can derive the weak form of  \eqref{eq:gov_elast} as: Find $\boldsymbol{u}\in \mathcal{U}$ such that 
\begin{equation*}
	\int_\Omega \frac{\mu}{2}(\boldsymbol{\nabla u}+\boldsymbol{\nabla u}^T)\cdot(\boldsymbol{\nabla v}+\boldsymbol{\nabla v}^T)\,d\Omega + \int_\Omega \lambda(\boldsymbol{\nabla} \cdot \boldsymbol{u})(\boldsymbol{\nabla} \cdot \boldsymbol{v})\, d\Omega =  \int_\Omega \boldsymbol{v}\cdot\boldsymbol{f}\,d\Omega + \int_{\Gamma_N}\boldsymbol{v}\cdot \boldsymbol{t}\,d\Gamma
\end{equation*}
for $\boldsymbol{v} \in \mathcal{V}$. Here, 
\begin{align*}
	\mathcal{U} &= \lbrace \boldsymbol{u} \in (H^1)^2 : \boldsymbol{u}=\boldsymbol{u_0} \mbox{ on } \Gamma_D \rbrace \mbox{ and } \\
	\mathcal{V} &= \lbrace \boldsymbol{v} \in (H^1)^2 : \boldsymbol{v}=\boldsymbol{0} \mbox{ on } \Gamma_D \rbrace,
\end{align*}
where $H^1$ is the space of functions with square integrable derivatives.

In this study, we consider two-dimensional problems with fixed isotropic materials under plane stress conditions, with $E=10^{3}$ and $\nu=0.3$.
The body force is set to $\boldsymbol{f}=0$.
Moreover, we consider fixed boundary conditions $\boldsymbol{u_0} = \boldsymbol{0}$ and pressure loading where $\boldsymbol{t}=10\boldsymbol{n}$.
Due to the linearity of the problem, different material properties or load magnitudes will still result in similar refinement patterns. 
We believe that the approach presented in this paper has the potential to be extended to a larger class of PDEs.

\subsection{The optimization problem}\label{optimization}
The aim of our method is to find the optimal mesh density for a given geometry, Dirichlet and Neumann boundaries.
This means that the mesh density produced by our method should {ideally} result in a discretization that achieves an optimal accuracy with a minimum number of degrees of freedom.
To formulate this objective in a mathematically precise way is a difficult undertaking, as there is no well-defined notion of optimality for the mesh density. Our proposed method can be regarded as an instance of operator learning. We therefore opt to describe the operator that we want to approximate by referring to a standard adaptive algorithm.

As described in the previous section, the input of our method consists of a geometry $\Omega$ as well as information on the fixed boundary $\Gamma_D\subset\partial\Omega$ and the traction boundary $\Gamma_N\subset\partial\Omega$.
We denote the set of possible input data as $X$, where each element $$x = (\Omega, \Gamma_D, \Gamma_N)\in X$$ consists of a geometry  and the Dirichlet and Neumann boundaries.

The operator that we want to approximate is therefore of the form
\begin{equation} \mathcal{F} : X \rightarrow D, \label{operator}\end{equation}
where $$D = \bigcup_{x\in X} D_x$$ and an element $d_x\in D_x$ is a function
$d_x: \Omega\rightarrow [0,1]$
that represents a relative local mesh density. 

The precise definition of the operator $\mathcal F$ depends on the chosen method for generating the locally refined mesh.
For a data point $x\in X$, we denote by $W_x$ the set of possible discretizations for solving the PDE using our chosen computational method.
In our examples, $W_x$ is the set of quadrilateral meshes that exactly represent the boundary $\partial\Omega$.

We write $$W = \bigcup_{x\in X}W_x$$
and we denote by
\begin{equation}\label{exactA}
	A: X \rightarrow W
\end{equation}
an algorithm for generating a discretization for the given input data.
For example, $A$ can be a classical adaptive method based on solving the PDE, marking, and refining the elements iteratively until a stopping criterion is reached.
The operator $\mathcal F$ is therefore defined as $$\mathcal F (x) = \mathcal D\circ A(x),$$
where
\begin{equation*}
	\mathcal D: W\rightarrow D
\end{equation*}
maps a discretization to its local relative mesh density.

Since it is fully deterministic, the operator $\mathcal F$ can be easily evaluated by applying the refinement algorithm $A$.
However, since the evaluation of $A$ typically consists of solving the PDE multiple times on meshes with increasing mesh density, this can be prohibitively computationally complex.
Therefore, we want to find a way to evaluate (or approximate) $\mathcal F$ in an  efficient way.
Since it is a highly non-linear operator that is impossible to evaluate directly, we opt to approximate $\mathcal F$ with an artificial neural network $\mathcal F_{\text{NN}}$.

To be able to use neural networks based on convolution, we consider as a discretized input a discrete approximation of $x = (\Omega, \Gamma_D, \Gamma_N)$ at a fixed grid of points $(\mathbf{p}_1, \ldots, \mathbf{p}_N) \in \Omega$, i.e., an image.
Likewise,  we regard the mesh density as a discrete function at the same set of points.
Therefore, for an input $x\in X$ and  target mesh density $d_x$, the loss function can be defined as:
\begin{equation} \mathcal{L} = \sum_{i=1}^N \left(\mathcal{F}_{NN}(x)(\mathbf{p}_i) - d_x(\mathbf{p}_i)\right)^2.\end{equation}

During the training, we try to find the optimal approximation $\mathcal{F}_{NN}$ by optimizing the network parameters.
Since we select the data points $\mathbf{x}_i$ as uniformly spaced, we investigate the use of image-based convolutional networks such as U-net, which are described in detail in Section~\ref{networkarchitecture}. 

\section{Training strategy and data representation}\label{complexity}
As shown in Section~\ref{optimization}, the problem of finding an optimal mesh density is complicated and largely depends on the domain $\Omega$ as well as the boundary conditions. The basic idea of our training strategy is to train the network with simple geometries, for which one can more easily obtain target discretizations. The trained network can then be applied to more complex domains. 
In this section, we give an overview of several measures of geometric complexity, the training strategy and the data representation.

\subsection{Measures of complexity}\label{complexityGeometry}

Since in our model problem we restrict ourselves to constant coefficients and zero body forces, the complexity of the mesh generation problem depends largely on the complexity of the domain.
While there is no way to objectively quantify the complexity of a geometry in this context, we aim to categorize possible input geometries by a number of measures corresponding to different aspects of geometric complexity.

One important advantage of using neural networks for this problem is their potential to generalize knowledge obtained from training on data from a number of simple data classes to complex unseen data.
We consider computational domains whose geometric complexity is measured using the following criteria:
\begin{itemize}
	\item \emph{Convexity}: Domains with only convex boundaries and domains with non-convex boundaries.
	\item \emph{Genus}: Simply connected domains and domains with a number of $k$ holes.
	\item \emph{Smoothness of Boundary:} Domains with polygonal boundaries and domains whose boundaries are piece-wise polynomial curves of degree $d\geq 2$.
\end{itemize}
For each of the geometric complexity criteria we consider a few simple cases, which are listed as follows:

\begin{tabular}[t]{lrl}
	``\emph{Convexity}'' & $=$ & $\left\{ 
	\begin{array}{ll}
		0 & \mbox{ if the outer boundary is convex{,}} \\
		1 & \mbox{ if the outer boundary is non-convex{,}}
	\end{array}
	\right.$ \\
	``\emph{Genus}'' & $=$ & $k \in \{0,1,2\}$, where $k=$ number of holes{,} \\
	``\emph{Smoothness of boundary}'' & $=$ & $\left\{ 
	\begin{array}{ll}
		0 & \mbox{ if all boundaries are polygonal{,}} \\
		1 & \mbox{ if at least one boundary curve is a spline{.}}
	\end{array}
	\right.$
\end{tabular}
\\

We then obtain geometric complexity classes, which we can denote with triples
\[
(\mbox{``\emph{Convexity}''},\mbox{``\emph{Genus}''},\mbox{``\emph{Smoothness of boundary}''}), 
\]
where e.g., $(0,0,0)$ represents convex, simply-connected, polygonal domains, $(1,2,0)$ represents non-convex polygonal domains with holes or $(0,0,1)$ represents convex, simply-connected, spline domains. {We randomly generated 10,000 sets of data each for the geometric complexity classes $(0,0,0)$ and $(1,0,0)$. While for the classes $(0,k,0)$ and $(1,k,0)$, which have higher geometric complexity, we generated 20,000 sets of data each. We then train and validate the neural networks by a random 90\% - 10\% split between the training and testing data.}

While we restrict ourselves to these criteria in this paper, the method can be extended arbitrarily by adding more criteria of complexity. Possible generalizations are discussed in more detail in Section~\ref{futurework}. Furthermore, when applying the method in practice, using real-world training data, the choice of geometric complexity classes does not need to be a conscious decision, as the network can learn from any data in the available data set.

\subsection{Training strategy}\label{sec:training_strategy}

In this study, we aim to train the neural networks such that they recognize the following refinement rules:
\begin{itemize}
	\item The refinement should be concentrated at points where the stress values are expected to be high, such as the reentrant corners and at the endpoints of the Dirichlet boundary.
	\item The refinement level depends on the position of the corner relative to the Dirichlet and Neumann boundaries.	
	\item The strength of the refinement also depends on the length of the support and traction boundaries and on their relative position.  
\end{itemize}

Instead of trying to implement these rules directly, we consider a more general data driven approach. A data set that contains domains of all possible geometric complexity classes, i.e., combinations of cases of complexity criteria listed in Section~\ref{complexityGeometry}, would need to be unfeasibly large and therefore difficult to generate and process. In this study, we exploit the neural network's potential for generalization to unseen data in order to handle complicated domains that are complex in more than one of these classes. 

We train the neural network on a data set that consists of domains that are simple in all but one of the different classes of geometric complexity. When evaluating the network on a domain that is complex in more than one of the listed classes, the network is able to apply the knowledge learned from this training data also to this case. An example is given in Section \ref{sec:exp_120}, where the training data set contains convex polygonal domains with one hole as well as non-convex polygonal domains without holes (i.e., the classes $(0,1,0)$ and $(1,0,0)$, respectively). When evaluating the network on a domain that is non-convex and has several holes, e.g. from the class $(1,2,0)$, the network can give a good prediction, although it has never seen such a domain.

\subsection{Data representation}
In our method, we represent the input as well as the output of the neural network as pixel-based images.
This choice leads to two significant advantages: On the one hand, we can employ powerful deep learning techniques that have been developed for image processing in the recent years. In particular, we can use convolutional neural networks and related network structures. On the other hand, using only image data makes our method completely independent from the representation of the geometry and from the discretization spaces that are used for solving the partial differential equation numerically.

For example, training data for our method can be obtained using a finite element method, isogeometric analysis, or any other adaptive method for solving PDEs. Likewise, the trained neural network has the potential to guide the generation of locally refined meshes for any of these numerical methods. What has to be taken into account is the interpretation of the output data and translation into a locally refined mesh.  Since we train our method with bilinear finite elements over quadrilaterals, the local mesh grading needs to be adapted for higher order or isogeometric elements.

\begin{figure}[th]
	\centering
	\begin{subfigure}[b]{0.2\textwidth}
		\includegraphics[width=1\textwidth]{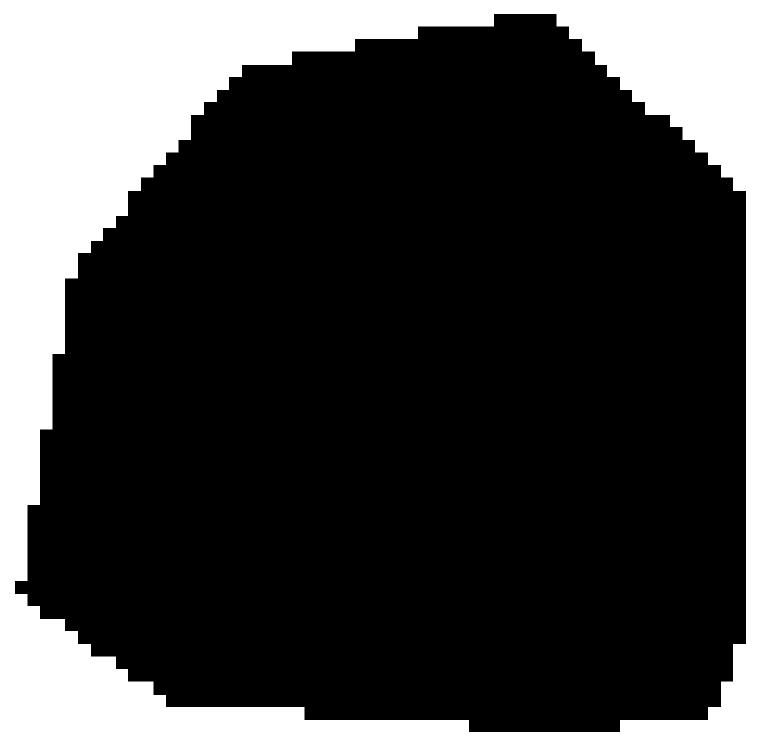}
		\caption{}
		\label{fig:geo_convexExp}
	\end{subfigure} 
	\begin{subfigure}[b]{0.2\textwidth}
		\includegraphics[width=1\textwidth]{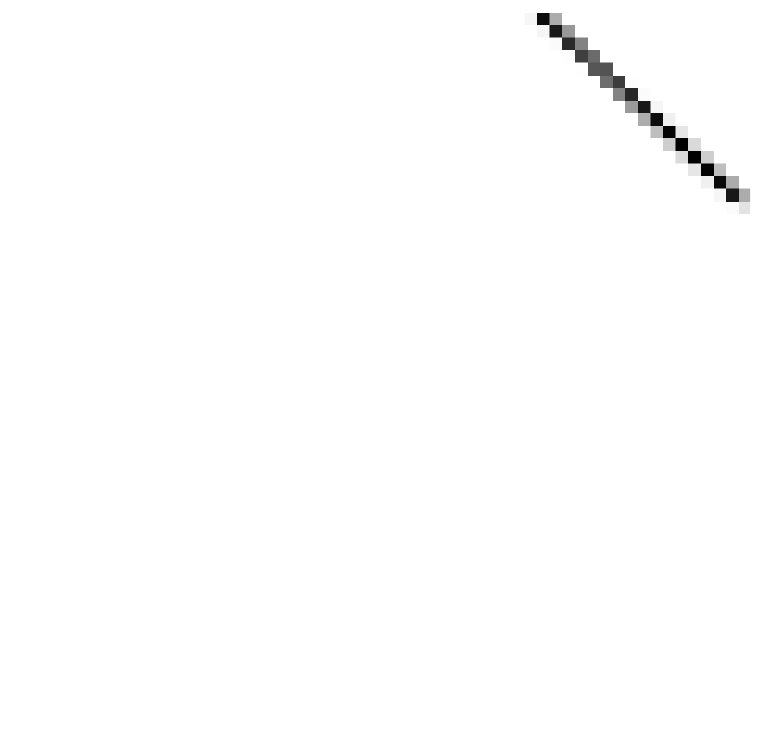}
		\caption{}
		\label{fig:fix_convexExp}
	\end{subfigure} 
	\begin{subfigure}[b]{0.2\textwidth}
		\includegraphics[width=1\textwidth]{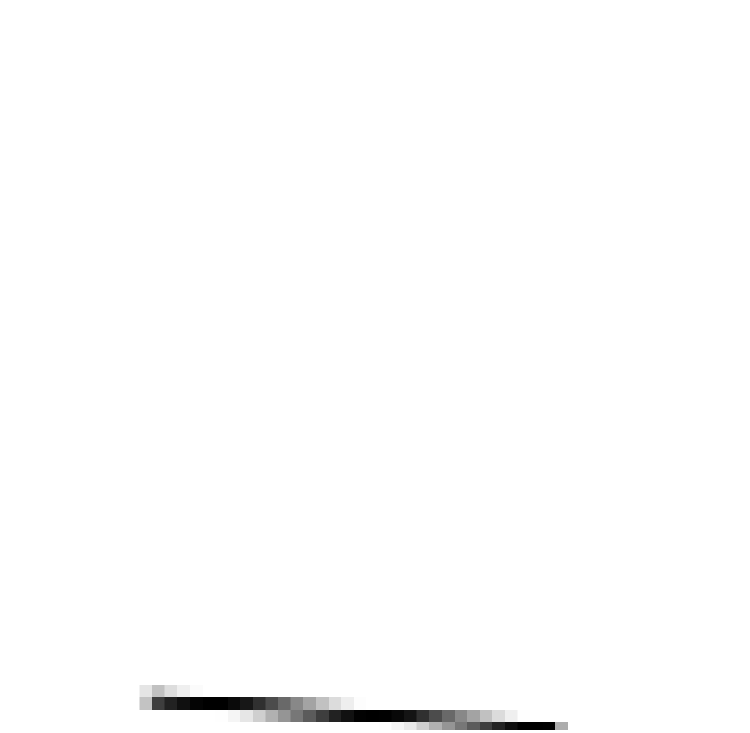}
		\caption{}
		\label{fig:trac_convexExp}
	\end{subfigure} 
	\begin{subfigure}[b]{0.25\textwidth}
		\includegraphics[width=1\textwidth]{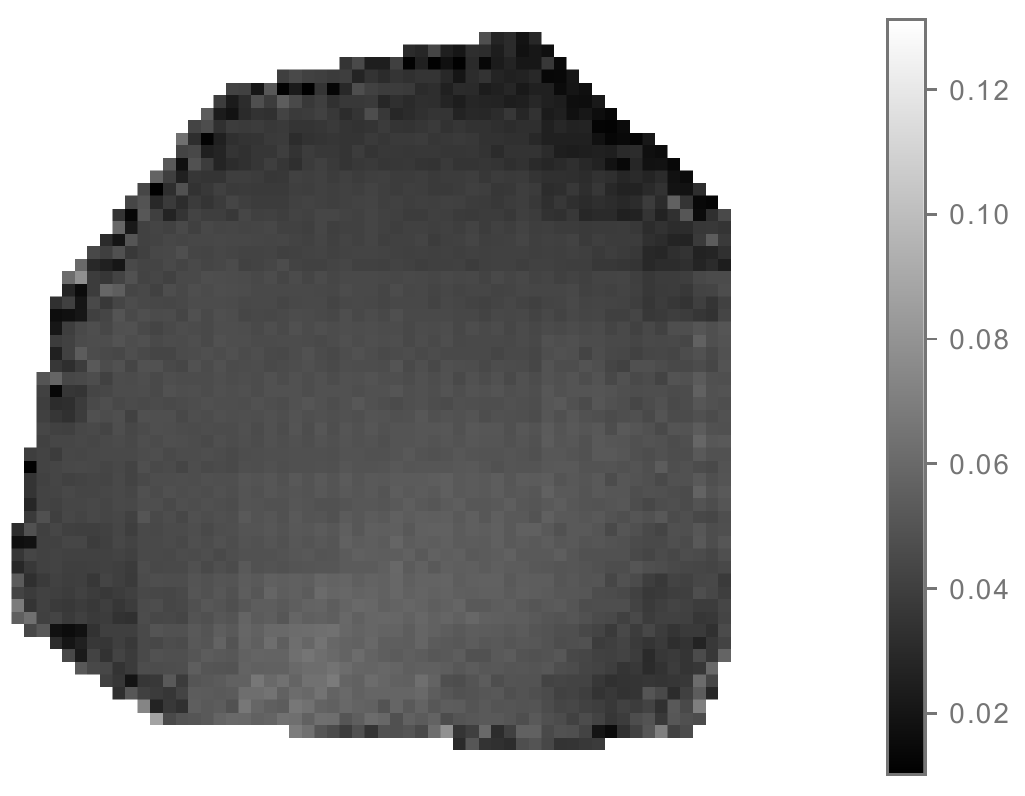}
		\caption{}
		\label{fig:outputExample}
	\end{subfigure} 
	\caption{An example for the input and output of our neural network: (a) Geometry, (b) Dirichlet boundary, (c) Neumann boundary and (d) predicted mesh density.}
	\label{fig:inputExample}
\end{figure}

The input of our neural network consists of three grayscale images of the same resolution:
One image contains the geometry of the computational domain $\Omega$.
The other two images contain only the Dirichlet boundary $\Gamma_D$ and the Neumann boundary $\Gamma_N$, respectively.
See Figures~\ref{fig:geo_convexExp} - \ref{fig:trac_convexExp} for an example input data.
The output of the network consists of a single grayscale image of the same size as the input data. The scalar value assigned to each pixel represents the network's prediction of the local mesh-size at that point of the input image representing the computational domain, see Figure~\ref{fig:outputExample}. 

\section{Data generation and adaptive refinement}\label{datagen}
In this section we explain how the training data sets, that is, the input geometry and the target mesh density, are generated.

\subsection{Generating the geometry} \label{sec:mesh_generation}

We have two generators which generate the data for 
geometry classes with varying \emph{Convexity} and \emph{Genus} (as discussed in Section \ref{complexityGeometry}), i.e., for classes $\{(i_1,i_2,0),i_1\in\{0,1\},i_2\in\{0,1,2\}\}$. In the following, we describe the procedures that we use to generate geometries for different geometric complexity classes. More precisely, we produce random convex polygonal domains, non-convex domains as well as domains with voids.

\begin{figure}[th]
	\centering
	\begin{subfigure}[b]{0.325\textwidth}
		\includegraphics[width=1\textwidth]{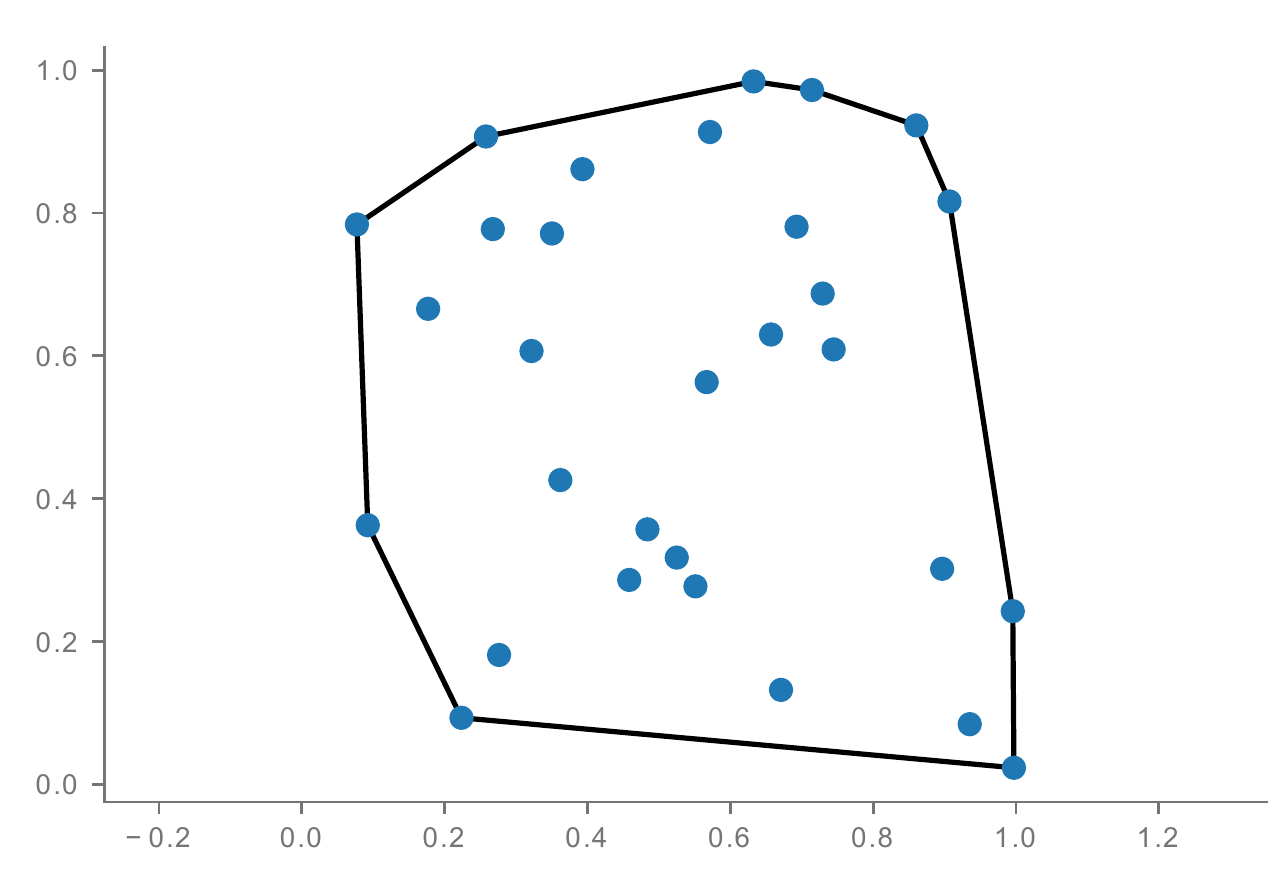}
		\caption{}
		\label{fig:generate_convex_polygon}
	\end{subfigure} 
	\begin{subfigure}[b]{0.2\textwidth}
		\includegraphics[width=1\textwidth]{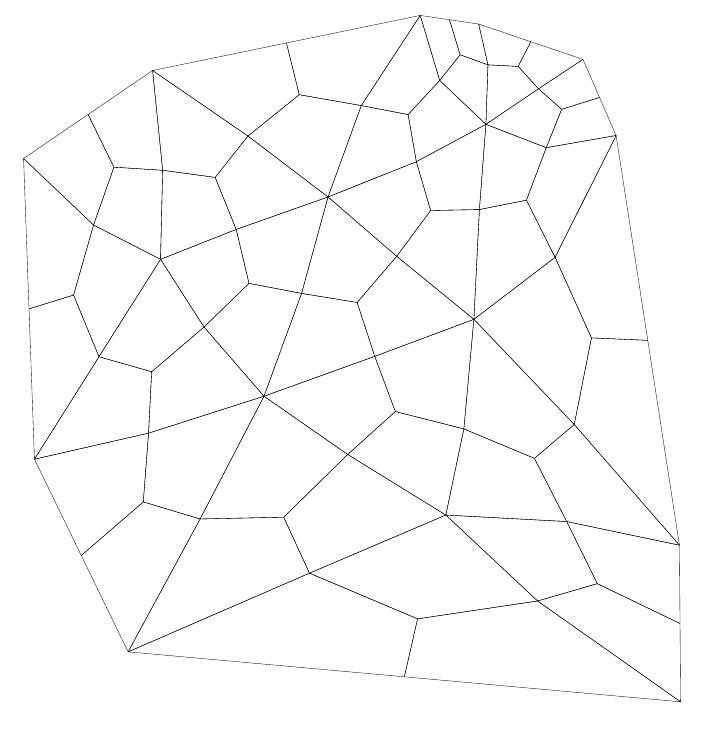}
		\caption{}
		\label{fig:convex_initialMesh}
	\end{subfigure} 
	\caption{Generating a convex domain (of class $(0,0,0)$) and initial mesh: (a) sample points and vertices of the convex hull and (b) resulting initial quad mesh obtained from Gmsh.}
	\label{fig:mesh generation}
\end{figure}

To create a domain with a convex boundary, we first generate $n$ sample points distributed randomly within the unit square. In our implementation we use $n=30$ for all training data sets. The domain is then given by the convex hull of the point set, see Figure \ref{fig:generate_convex_polygon}, where the blue dots are the sampled points and the black line represents the boundary of the convex hull.  The blue dots that lie on the black line are the vertices of the convex hull. They are used to construct a quadrilateral mesh in Gmsh, as shown in Figure \ref{fig:convex_initialMesh}.

\begin{figure}[th]
	\centering
	\begin{subfigure}[b]{0.325\textwidth}
		\includegraphics[width=1\textwidth]{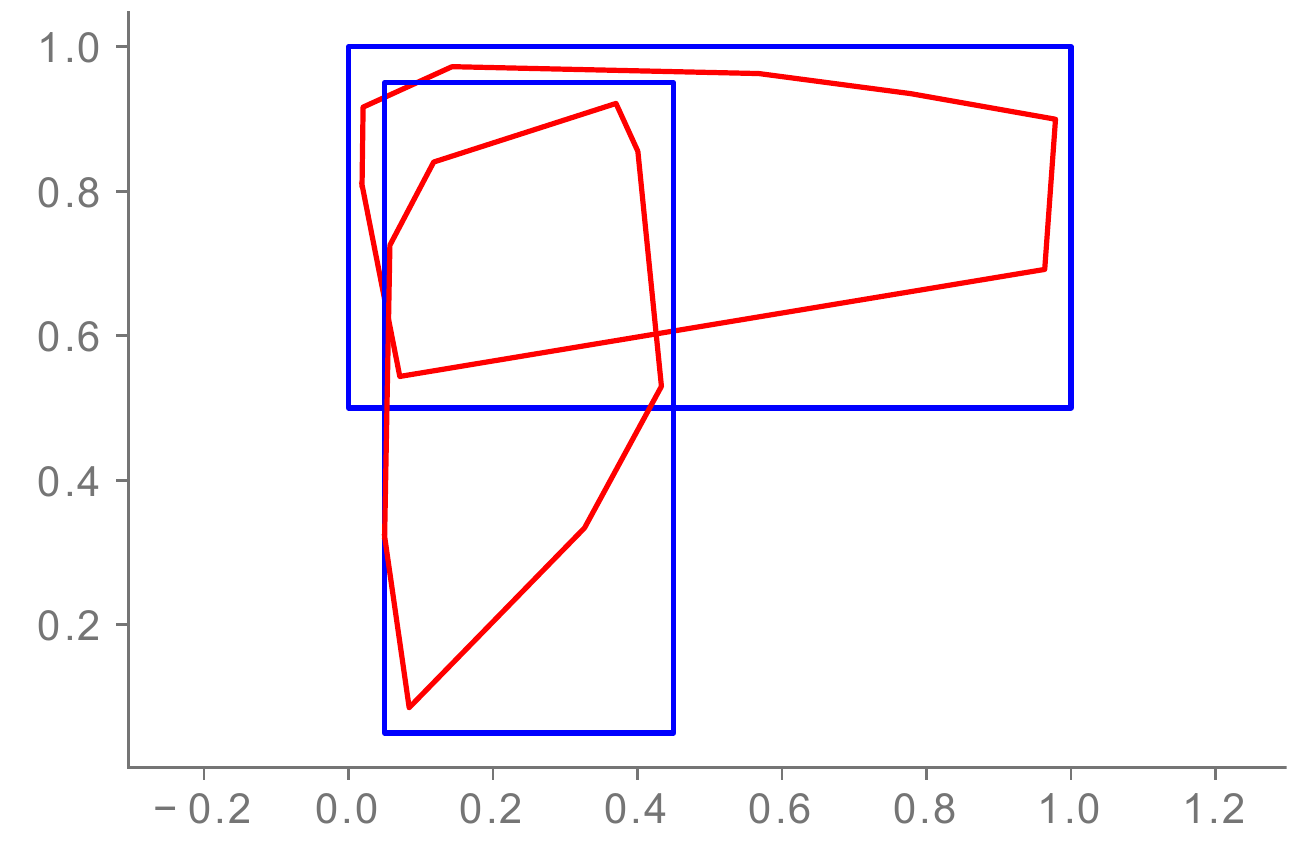}
		\caption{}
		\label{fig:concave_geometry}
	\end{subfigure} \quad\quad
	\begin{subfigure}[b]{0.2\textwidth}
		\includegraphics[width=1\textwidth]{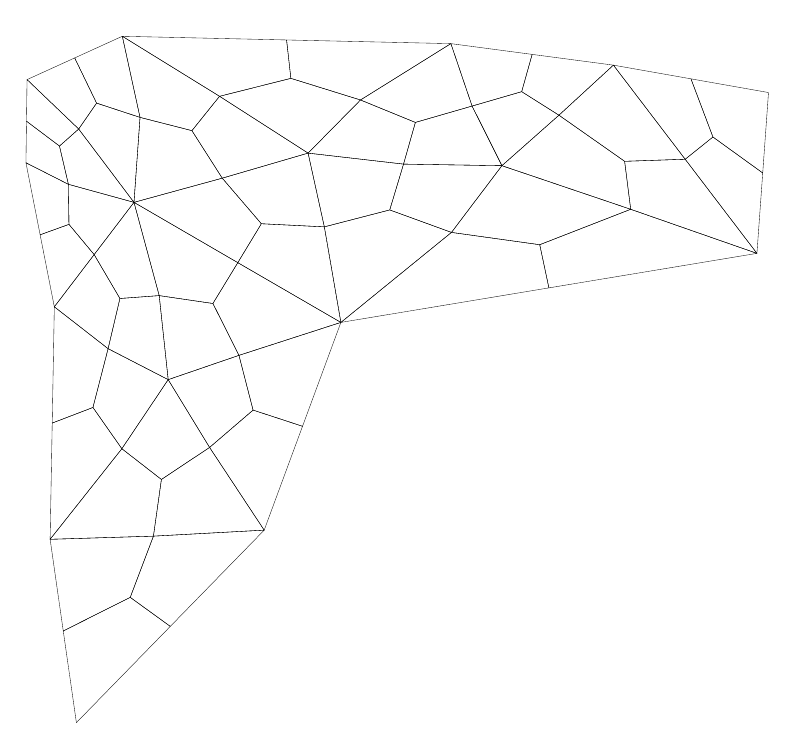}
		\caption{}
		\label{fig:concave_geometry_mesh}
	\end{subfigure} 
	\caption{Example for generating a non-convex domain (of class $(1,0,0)$).}
	\label{fig:concave_geometry_mesh_exp}
\end{figure}

In Figure \ref{fig:concave_geometry}, we illustrate the creation of a non-convex domain as a union of two convex domains. Two rectangles that intersect and lie within the unit square are first created. Then, as described above, a convex domain is created for each rectangle as the convex hull (represented by the red polygon) of random sample points. The non-convex domain is then given as the union of the two convex sub-domains; the corresponding mesh is shown in Figure \ref{fig:concave_geometry_mesh}. If the resulting domain is not simply connected and non-convex, it is discarded.

\begin{figure}[th]
	\centering
	\begin{subfigure}[b]{0.325\textwidth}
		\includegraphics[width=1\textwidth]{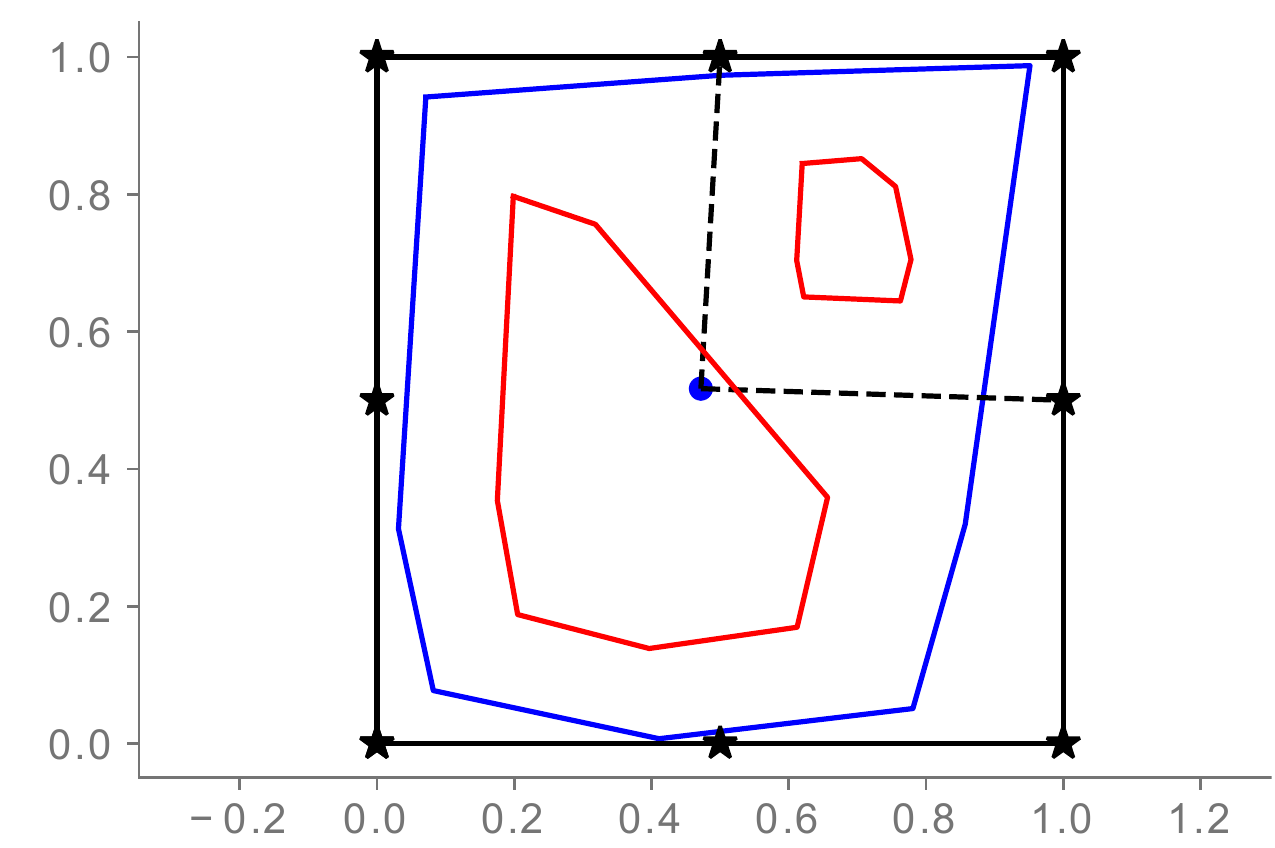}
		\caption{}
		\label{fig:create_void}
	\end{subfigure} \quad\quad
	\begin{subfigure}[b]{0.2\textwidth}
		\includegraphics[width=1\textwidth]{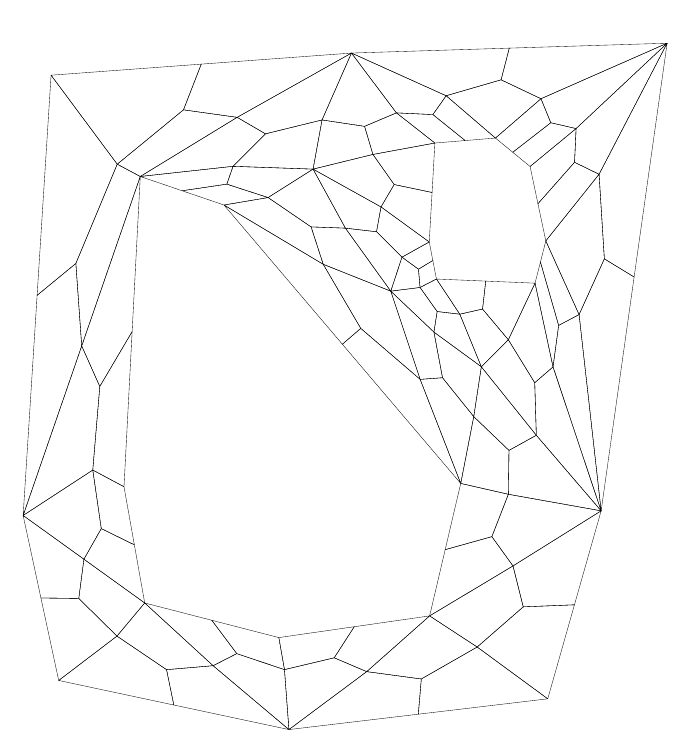}
		\caption{}
		\label{fig:convex_void_exp}
	\end{subfigure} 
	\caption{Example of generating a domain with voids (of class $(0,2,0)$).}
	\label{fig:mesh generation_void}
\end{figure}

Creating a domain of non-zero genus is illustrated in Figure \ref{fig:mesh generation_void}. First, a random convex or non-convex domain is generated, which is then divided into two or more sub-domains. In Figure \ref{fig:create_void}, a convex domain (blue polygon) is divided into two sub-domains by the dashed lines. The division lines are created by randomly selecting points that lie on the boundary of the unit square (marked by stars) and connecting them to the center of the polygon (marked by the blue dot). Then, for each sub-domain, a small convex domain is created (represented by the red polygon), again as the convex hull of random sample points (with $n=10$). These small domains are then cut out from the main polygon, forming the voids. Figure \ref{fig:convex_void_exp} shows the resulting mesh. 

The boundary conditions are generated by randomly selecting two edges from the domain. Here we only consider the outer boundary of the domain. The restriction that we make is that the two boundaries should not be adjacent to each other. One of the selected boundaries is set as the Dirichlet boundary, while the other is set as the Neumann boundary. The domain as well as the edges selected to impose boundary conditions are converted to images which are used as input for the considered data set, cf. Figure~\ref{fig:inputExample}. 

Once the domain is constructed and the boundary conditions are fixed, we create an initial mesh over the domain and perform a numerical simulation for the linear elasticity problem (see Section \ref{pdeProblem} for more details) using a finite element solver. This is explained in more detail in the following subsection.

\subsection{Generating the target mesh density}\label{generating-target-density}

For each domain, the sampled points that lie on the boundary are selected and used to construct an initial mesh using Gmsh\footnote{\url{https://gitlab.onelab.info/gmsh/gmsh/-/blob/master/api/gmsh.py}}. Gmsh is an open-source program written in C++ with an extensive Python interface. It produces unstructured triangular and quad (or mixed) meshes based on a given boundary representation. In our experiments, an all-quad mesh is constructed from a triangle mesh by subdivision. 

For the numerical analysis we use SolidsPy\footnote{\url{https://github.com/AppliedMechanics-EAFIT/SolidsPy}} which is written in Python/NumPy. It is an educational software which is easy to use and modify. SolidsPy supports triangle and quad meshes. In addition, it performs stress averaging for plotting continuous stress fields. For the given mesh and boundary conditions, information such as displacement, stress and strain can be obtained. The von Mises stress can be computed from the $\sigma_{xx}$, $\sigma_{yy}$ and $\sigma_{xy}$ stress fields.

\begin{figure}[th]
	\centering
	\begin{subfigure}[b]{0.24\textwidth}
		\includegraphics[width=1\textwidth]{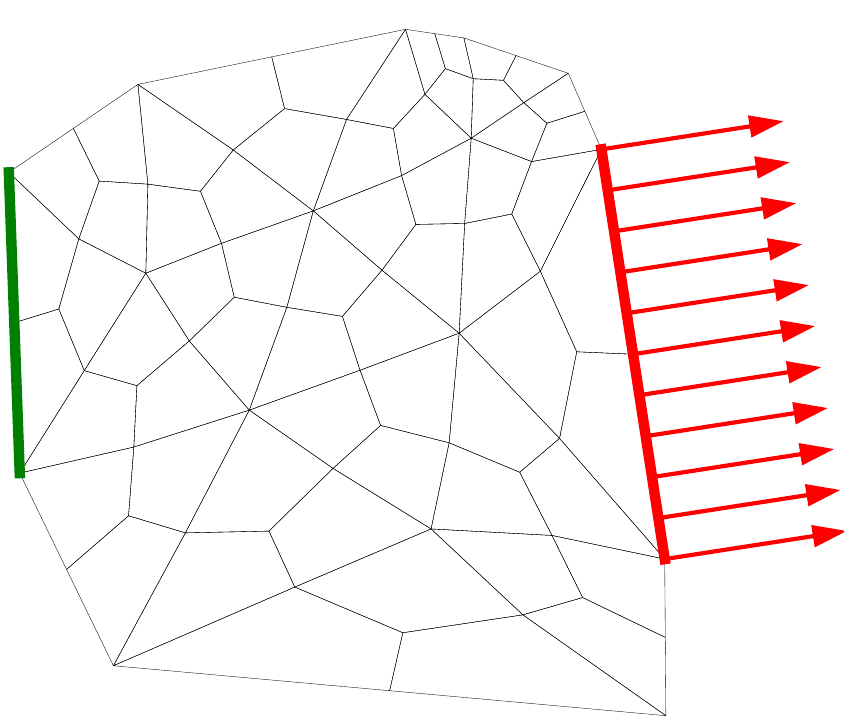}
		\caption{}
		\label{fig:convex_adaptive_0}
	\end{subfigure} 
	\begin{subfigure}[b]{0.24\textwidth}
		\includegraphics[width=1\textwidth]{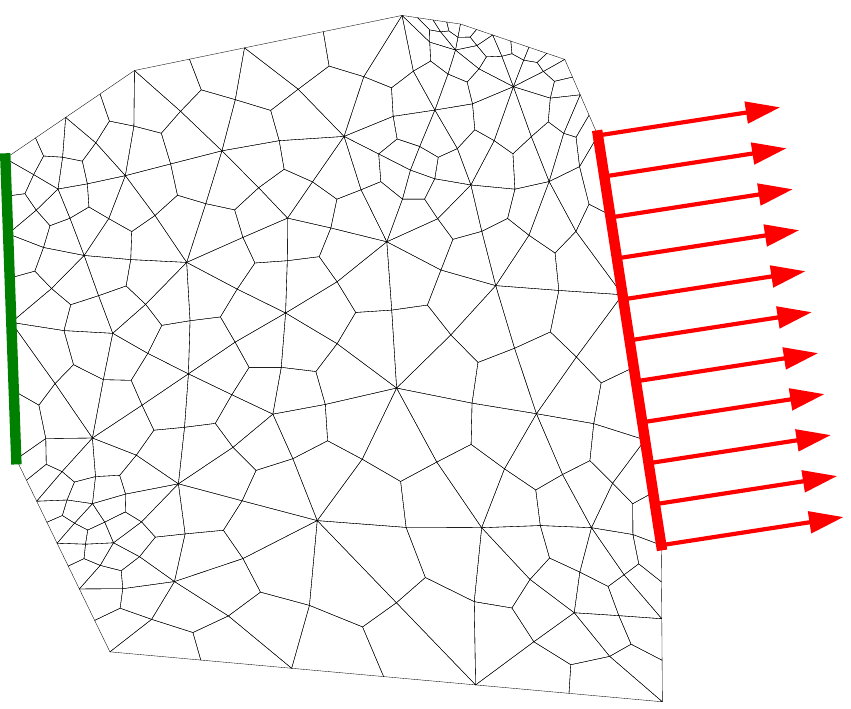}
		\caption{}
		\label{fig:convex_adaptive_1}
	\end{subfigure} 
	\begin{subfigure}[b]{0.24\textwidth}
		\includegraphics[width=1\textwidth]{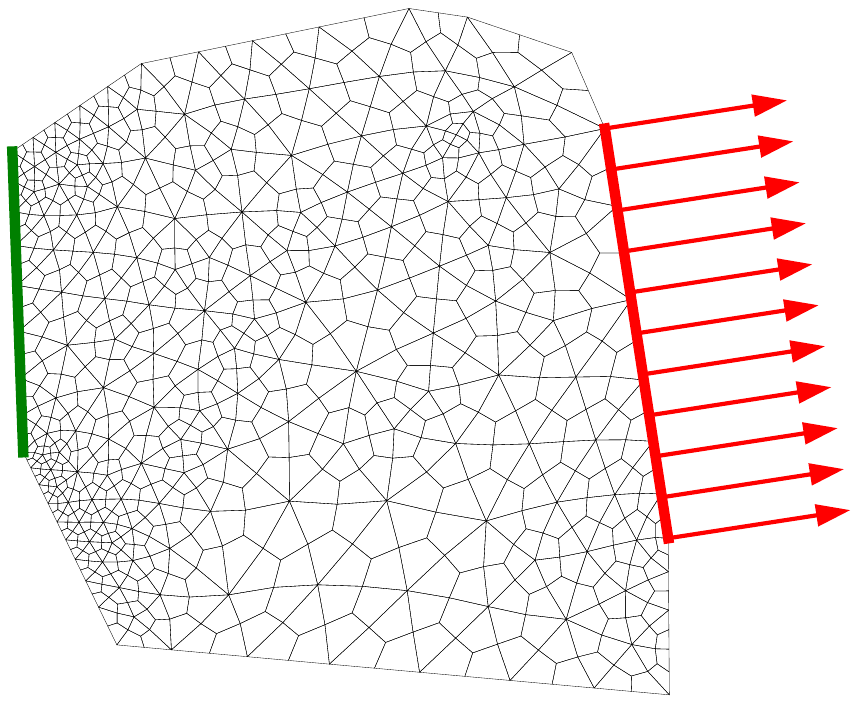}
		\caption{}
		\label{fig:convex_adaptive_2}
	\end{subfigure}
	\begin{subfigure}[b]{0.24\textwidth}
		\includegraphics[width=1\textwidth]{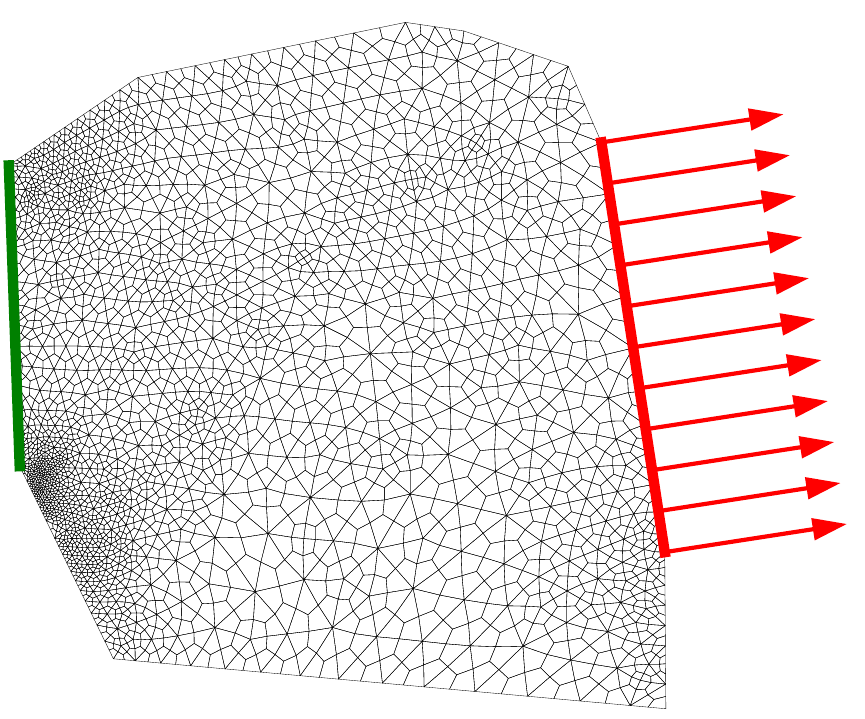}
		\caption{}
		\label{fig:convex_adaptive_3}
	\end{subfigure}
	\caption{Example of adaptive refinement.}
	\label{fig:convex_adaptive}
\end{figure}

To obtain an adaptively refined mesh, we use a uniformly refined fine mesh $\mathcal{M}^{*}$ for a heuristic error estimation. Consider a mesh $\mathcal{M}_{\ell}$ at level $\ell$. The von Mises stress is computed on the mesh  $\mathcal{M}_{\ell}$ and on $\mathcal{M}^{*}$. By using the error between the fine and coarse solution measured at every vertex of $\mathcal{M}_{\ell}$, we generate a mesh size field $\mathcal{B}_{\ell}$ as follows:

\begin{equation}	
	\mathcal{B}_{\ell}=\mathbf{h} \left(1-\frac{\boldsymbol{\varepsilon}_{avg}}{2\varepsilon_{max}} \right),
	\label{eq:mesh_size_field}
\end{equation} 
In \eqref{eq:mesh_size_field}, $\mathbf{h}$ and $\boldsymbol{\varepsilon}_{avg}$ represent the maximum edge length and average von Mises error on each element, while  $\varepsilon_{max}$ denotes the global maximum von Mises error in the mesh.  We refer to the Gmsh tutorial\footnote{\url{https://gitlab.onelab.info/gmsh/gmsh/-/blob/master/tutorials/python/t10.py}} for using the estimated mesh size field $\mathcal{B}_{\ell}$ to create an updated mesh $\mathcal{M}_{\ell+1}$.
An example of adaptive refinement obtained using the method as described above is shown in Figure \ref{fig:convex_adaptive}, where the meshes are adaptively refined from Figures \ref{fig:convex_adaptive_0} - \ref{fig:convex_adaptive_3}. In the figures, the Dirichlet and Neumann boundary are highlighted in green and red respectively, while the red arrows indicate the traction direction. We choose this heuristic refinement approach to reduce the influence of the initial mesh and of the refinement procedure that is usually present for adaptive schemes based on marking elements.

\section{Network architecture}\label{networkarchitecture}

In the following we discuss the structure of the neural network architecture that we use in this paper. Usually, neural networks are built up of layers sequentially. The network receives an input which is transformed through a series of hidden layers followed by an output layer. The hidden layers and the output layer are made up of neurons associated with weights and biases. Each neuron takes as input the outputs of the neurons of the previous layer. The output of the neuron is given as the composition of an affine linear combination of the input values (taking the weights as coefficients and adding/subtracting the bias), composed with an activation function. The composition of several such layers yields a non-linear function if the activation function is non-linear. In case of a fully connected neural network, each neuron takes as input all neurons from the previous layer (with possibly non-zero weights).

The weights and biases of all the neurons are the degrees of freedom which are trained, i.e., for which one optimizes. The entire network is then trained based on a given set of input and target data, and the weights and biases are optimized with respect to a given objective. The objective is described in terms of a scalar-valued loss function. The aim of the training is to find the parameters (weights and biases) which minimize the loss function by using gradient descent type methods. Computing the gradients is accomplished by reverse-mode differentiation (back-propagation).

While fully connected neural networks are useful for solving many optimization problems, when handling higher dimensional data, they tend to result in a huge number of parameters and can lead to overfitting. Therefore, when the data possesses some local structure which is of relevance for the output of the network, one can use neural networks based on convolution.

\subsection{Convolutional neural networks }

\begin{figure}[th]%
	\centering
	\includegraphics[width=0.45\textwidth]{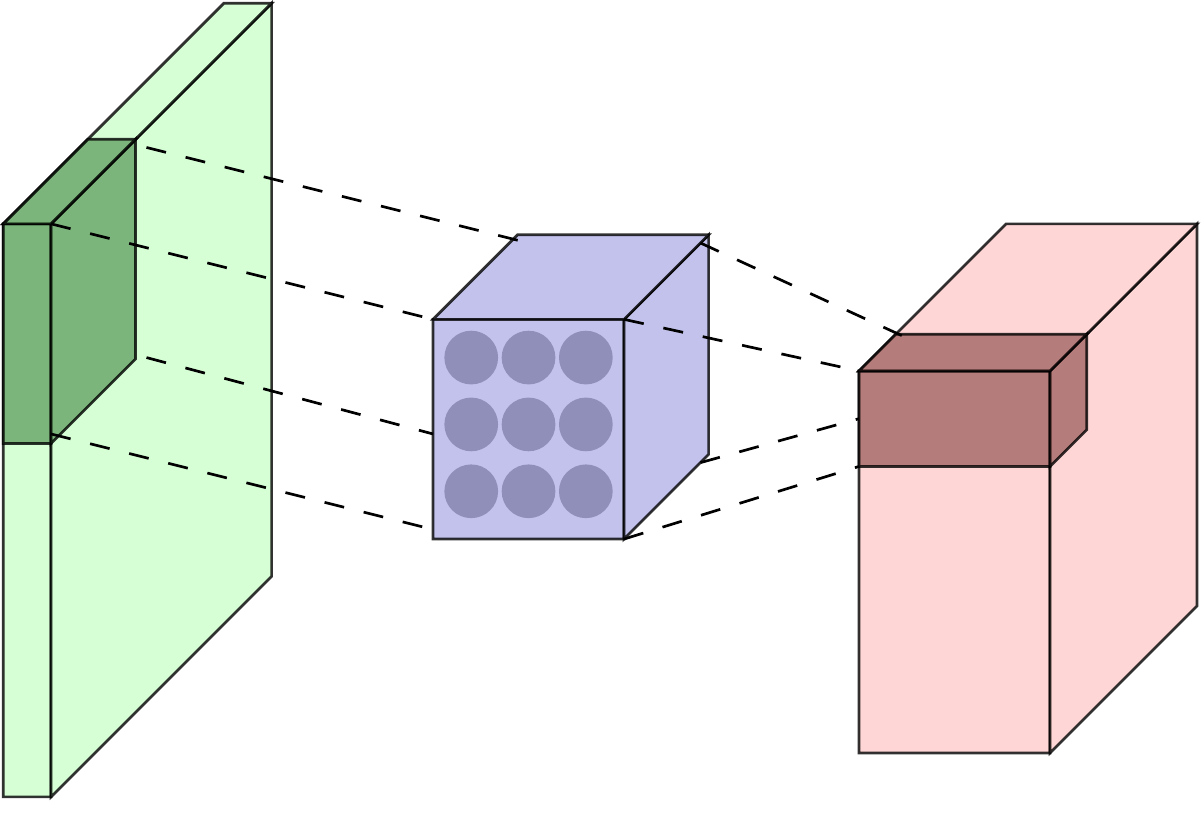}
	\caption{Convolutional neural networks.}
	\label{fig:CNN}
\end{figure}
Convolutional neural networks (CNN) are explicitly used for working with 2D or 3D data such as images. CNNs are commonly used for image classification. Spatially coherent features such as corners and lines in the images are identified and combined to predict an associated label or category. Similarly to fully connected neural networks, they are formed by sequences of layers. An example of a CNN is illustrated in Figure \ref{fig:CNN}. It contains an input (green block), a convolution layer (blue block), and an output (red block). In the convolution layer, the filter is formed by neurons arranged in three dimensions. The filter “slides” over the input, performing element-wise multiplication with the current part of the input (highlighted with dark green in Figure \ref{fig:CNN}), and resulting in a local output (highlighted with dark red in Figure \ref{fig:CNN}). 

The output size is determined by several factors including the filter size, input size, padding and stride. Padding refers to adding additional boundary data surrounding the image. It is commonly used for reducing border effects. Stride is defined as the number of steps the filter moves in the convolution. Suppose the input size is given by $W_I$$\times$$H_I$$\times$$D_I$, convoluted with a  $W_F$$\times$$H_F$$\times$$D_F$ filter, $P$ padding, and $S$ stride will produces output of size $[(W_I-W_F+2P)/S+1]$$\times$$[(H_I-H_F+2P)/S+1]$$\times$$D_F$. Meanwhile, techniques such as pooling and regularization can be added to the network for deep learning enhancement. 

\subsection{The U-net architecture}
\label{sec:Unet}

\begin{figure}[th]%
	\centering
	\includegraphics[width=0.65\textwidth]{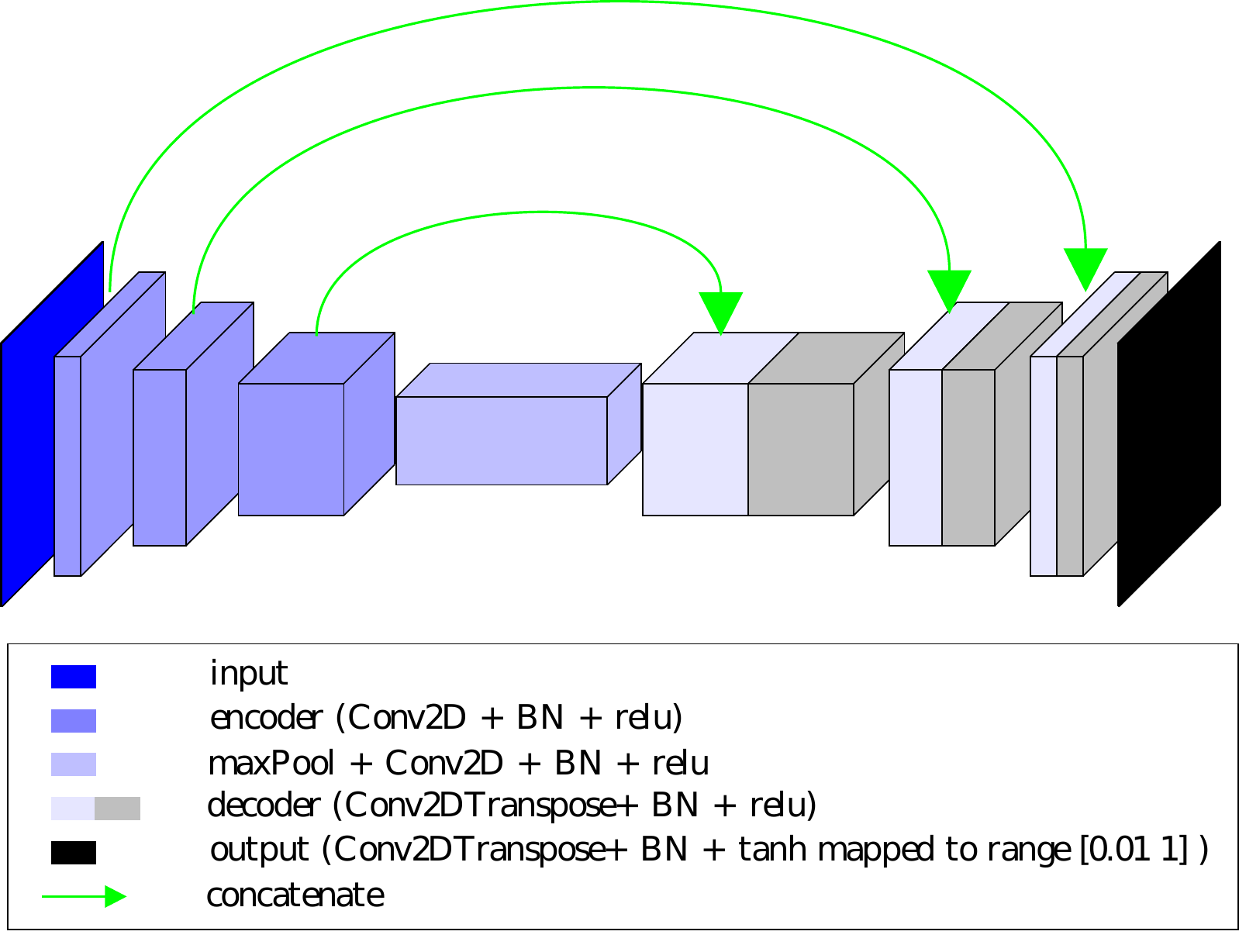}
	\caption{U-net architecture for mesh size prediction. It consists of contracting and expanding parts.}
	\label{fig:unet}
\end{figure}

In this study, we develop a network to learn the mesh density operator by using the U-net architecture~\cite{ronneberger2015}. The network consists of two parts: contracting and expanding. The first part is formed by encoder blocks which downsample the data, while the second part contains decoder blocks which upsample the data. Each block consists of convolution layers (transposed convolutions for the decoder block), and batch normalization. The two parts are connected by a bottleneck block which has a max pooling layer followed by convolution, and batch normalization. The U-net for training the data is illustrated in Figure \ref{fig:unet}. 
We adopt U-net architecture because it is capable of capturing the contents and enables precise localization. The contracting part of the U-net is the same as a general CNN which extract features such as reentrant corners, and the relative location of the Dirichlet and Neumann boundaries from the inputs. The network then up-samples its hidden layers to create a gray scale image. The gray scale image will have intensity values ranging from 0 to 1, which indicate the relative mesh resolution. This architecture allows the network to combine features from different spatial regions of the images and to localize more precisely the mesh resolution at the region of interest.

The main difference between the U-net architecture as in~\cite{ronneberger2015} and the U-net architecture used in this study is that, in our case, there is only one max pooling layer located at the bottleneck block. The low-level feature map from the convolution layers and the precise position of the features in the inputs are important for meshing which is sensitive to geometry. For this reason, we limit the use of pooling having the effect of making the representation become approximately invariant to small translations of the input \cite{Goodfellow2016}. 

\begin{figure}[th]%
	\centering
	\includegraphics[width=0.85\textwidth]{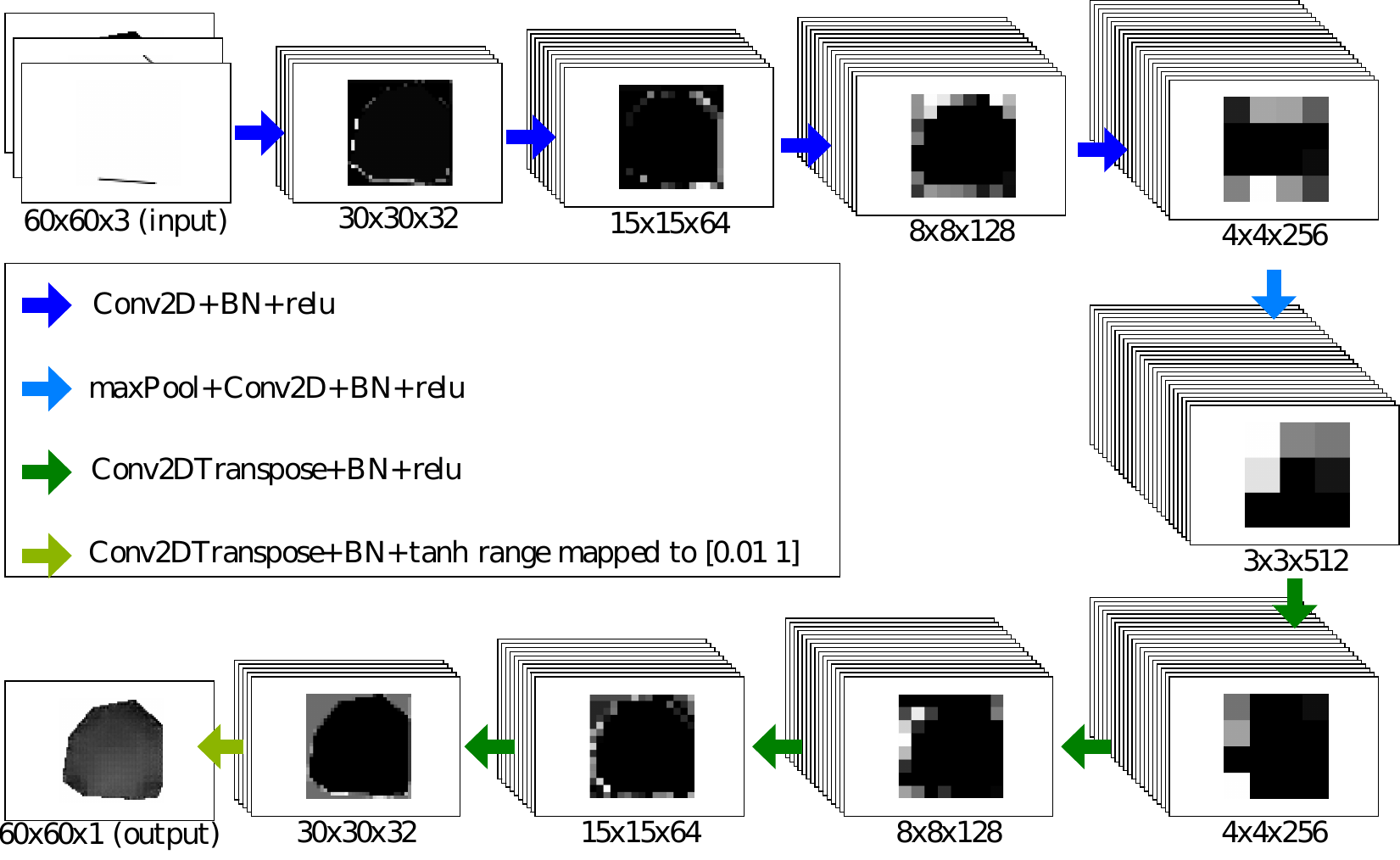}
	\caption{{Representative outputs for the U-net structure shown in Figure~\ref{fig:unet}. The input of the network contains three layers which represent the geometry, Dirichlet boundary and Neumann boundary, respectively, as shown in Figures \ref{fig:geo_convexExp}-\ref{fig:trac_convexExp}.}}
	\label{fig:unet_output}
\end{figure}

We control the down-sampling through padding and striding in the convolution. We apply 2 strides which reduce the size of the convolution output by half, therefore, reducing the image size as the algorithm is moving downward in the contracting path. At the end of the contracting path, the size is reduced to $4\times4$. Max-pooling is then used to further down-sample to $3\times3$ followed by de-convolution, which is the reverse of convolution in terms of  dimensions change at each layer in the expanding path. For the decoder block, the filter size, padding and stride are selected such that the tensor size matches with the contracting path in reverse order. Concatenation is used to add the tensors of matching size between contracting path to expanding path. The final output is a  $60\times60$ image with a single channel. We set the depth of the filter at the first encoder block to $D_F=32$. It is doubled from one block to the next in the contracting part. At the max pool block, the depth is increased to $D_F=512$. The depth is then decreased by half until $D_F=32$ at the last decoder block. {Output samples from the U-net are given in Figure~\ref{fig:unet_output}. Each encoder, decoder and max pool block provides a group of outputs. The image resolution and the number of layers are different after each block. For example, the resolution decreases and the depth is increasing from one encoder block to another. For each group, we show a representative layer (slice through the depth) and label the actual size.}

\subsection{Model evaluation}
In this study, we seek for the operator that maps given geometry and boundary conditions to an optimal (adaptive) mesh density. The inputs and outputs are stored as grayscale images which can be converted to matrices of a fixed size corresponding to the image resolution. 
The mean square error is used as the regression loss function. It is defined by
\begin{equation}	
	MSE=\frac{\sum_{i=1}^n(x_i-x_i^p)^2}{n},
	\label{eq:MSE}
\end{equation} 
where $x$ refers to the intensity value at pixel $i$, and $n$ is the number of pixels. 

\section{Numerical examples}\label{numericalexamples}

We conduct some experiments for testing the proposed mesh density estimation method. By using the U-net as discussed in Section \ref{sec:Unet}, we train six models. The first four models are \emph{Model$(0,0,0)$}, \emph {Model$(0,k,0)$}, \emph {Model$(0,k,0)$}, and \emph {Model$(1,k,0)$}, $k\in\{1,2\}$, where the associated triplets (introduced in Section \ref{complexityGeometry}) indicate the geometry complexity classes used for training. Since each of these models are trained with data restricted to simple geometries or groups of highly related geometries, they are expected to yield more accurate results. However, this choice requires that the underlying complexity class of the example must be known beforehand. Thus, it is more reasonable to have models trained for all possible geometric complexity classes. Here we examine two extended models which are referred to as Model with Reduced Training (\emph{Model RT}) and Model with Complete Training (\emph{Model CT}). They are trained using combinations of classes as listed below:
\begin{itemize}
	\item \emph {Model RT}: $(0,0,0)$, $(1,0,0)$ and  $(0,k,0), k\in\{1,2\}$  
	\item \emph {Model CT}: $(0,0,0)$, $(1,0,0)$, $(0,k,0), k\in\{1,2\}$ and $(1,k,0), k\in\{1,2\}$
\end{itemize} 
While \emph{Model CT} covers all possible cases, model \emph{Model RT} is easier to train, as no training data of high complexity $(1,k,0)$, $k\in\{1,2\}$, needs to be constructed. Nonetheless it may be applied to examples of high complexity.

A comparison of ground truth and output for all the models are given in Section \ref{sec:comparison_rel}. In Sections \ref{sec:exp_000} and \ref{sec:exp_100} show the results of \emph{Model$(0,0,0)$} and  \emph{Model$(1,0,0)$} handing simple geometric domain. The different effect of \emph{Model$(1,k,0)$}, \emph {Model RT} and \emph {Model CT} on a complex geometry domain is given in Section \ref{sec:exp_120}. Meanwhile, an example with smooth boundary is shown in Section \ref{sec:exp_011}, which demonstrates the generalization properties of the corresponding \emph{Model$(0,k,0)$}. Finally, the overall performances are quantified by using histograms in Section \ref{sec:eval_of_model}. We also show the plot of the loss function vs. training epochs for the models at the end of this section. In the following section we introduce the error measures that we consider.

\subsection{Error measures}

We used relative error for evaluation of predicted output. Let $I^G$ and $I^P$ be the $n_p \times n_p$ matrices representing the ground truth and predicted output images respectively, where $n_p$ is the number of pixels in each space direction. Here $I^G_{ij} \in [0,1]$ and $I^{P}_{ij} \in (0,1]$ denote the value of the mesh size for the subdomain covered by the pixel with index $ij$. We separate the image into two disjoint subdomains $\Omega^{I}$ and $\Omega^{O}$, where  $\Omega^{I}$ denotes the pixels which are (at least partially) in the computational domain and $\Omega^{O}$ denotes the pixels which are completely outside the domain. Furthermore, we set $I^G_{ij}=1$ for the pixels which belong to $\Omega^{O}$, whereas $I^G_{ij}<1$ for pixels $\Omega^{I}$. For a given domain and set of boundary conditions, the accuracy of predictive mesh size is evaluated by using relative error ($\varepsilon_{rel}$) between the ground truth and the preprocessed output ($\tilde{I}^{P}$) defined by: 
\begin{equation}	
	\varepsilon_{rel}=\frac{\left\|I^G- \tilde{I}^{P}\right\|}{\left\|I^G\right\|}
	\label{eq:rel_error}
\end{equation} 
where $\tilde{I}^{P}$ is obtained by setting $I^{P}_{ij}$ within $\Omega^O$ to be 1 and $||\cdot||$ denotes the $L^2$ norm. $\tilde{I}^{P}$ is adopted for calculating $\varepsilon_{rel}$ such that only the computational domain of interest is considered when computing the error. 

For validation of the linear elasticity solution, we calculate the error in the relative energy norm of the finite element solution \cite{sukumar2004finite}. It is defined in terms of exact stress $\boldsymbol{\sigma}$, computed stress $\boldsymbol{\sigma}^h$, and coefficient matrix $\boldsymbol{D}$ written in Voigt notation as:
\begin{equation}
	e_{rel} =\frac{\sqrt{\frac{1}{2}\int_\Omega (\boldsymbol{\sigma}-\boldsymbol{\sigma}^h)^T \boldsymbol{D}^{-1} (\boldsymbol{\sigma}-\boldsymbol{\sigma}^h) d\Omega}}{\sqrt{\frac{1}{2}\int_\Omega \boldsymbol{\sigma}^T \boldsymbol{D}^{-1} \boldsymbol{\sigma} d\Omega}}
	\label{eq:energy_norm}
\end{equation} 
and 
\[ \boldsymbol{D} = \frac{E}{1-\nu^2}\begin{bmatrix} 1 & \nu & 0 \\ \nu & 1 & 0 \\ 0 & 0 & 0.5(1-\nu) \end{bmatrix}.
\]

\subsection{Comparison of relative error}\label{sec:comparison_rel}

\begin{figure}[th]%
	\centering
	\includegraphics[width=0.75\textwidth]{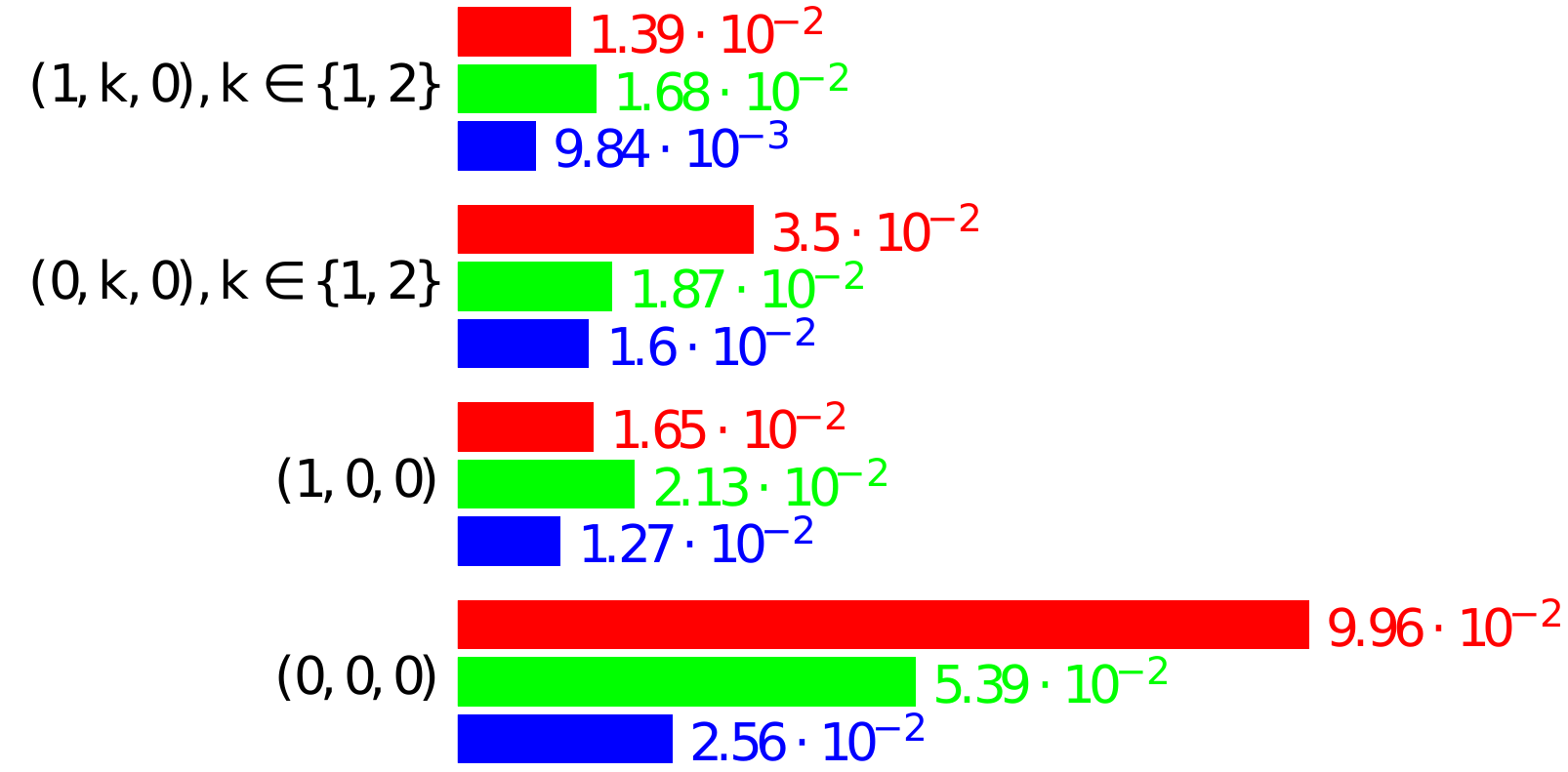}
	\caption{Comparison of average $\varepsilon_{rel}$ of test data sets for each geometry group on different models. The length of blue bars correspond to average $\varepsilon_{rel}$ from \emph{Model $(0,0,0)$}, \emph{Model $(1,0,0)$}, \emph{Model $(0,k,0)$}, and \emph{Model $(1,k,0)$} respectively. Green bars and red bars shows the average $\varepsilon_{rel}$ associate to \emph{Model RT} and \emph{Model CT}.}
	\label{fig:chart}
\end{figure}

As mentioned previously, the idea of the proposed method is to provide the geometry and boundary conditions as input (in the form of images) to the network and receive a gray scale image (the predicted mesh size density) as output. Here, we compare of relative error between the predicted output and ground truth from the different models. Figure \ref{fig:chart} shows the average relative error computed on the test data set for geometric complexity classes as labeled on left side of figure. 

The blue bars in the chart are the results from models trained with the respective data sets, i.e. \emph{Model $(0,0,0)$}, \emph{Model $(1,0,0)$}, \emph{Model $(0,k,0)$}, and \emph{Model $(1,k,0)$}. They always have the smallest average relative error. The green bars are the results from \emph{Model RT}. \emph{Model RT} results in higher relative error for test data sets $(1,0,0)$ and $(1,k,0)$. This is because there are fewer data with non-convex boundaries used for training this model. The red bars are the results from \emph{Model CT}. Although this model is trained with data from all possible geometry complexity classes, it results in higher relative error. The reason is that the network learns better when the training data is consistent. Nonetheless, \emph{Model CT} could be useful for prediction of mesh size density when there is no information given on the geometry complexity class.

\subsection{Example for geometric complexity class (0,0,0)}\label{sec:exp_000}
\begin{figure}[hbt!]
	\centering
	\begin{subfigure}[b]{0.24\textwidth}
		\includegraphics[width=1\textwidth]{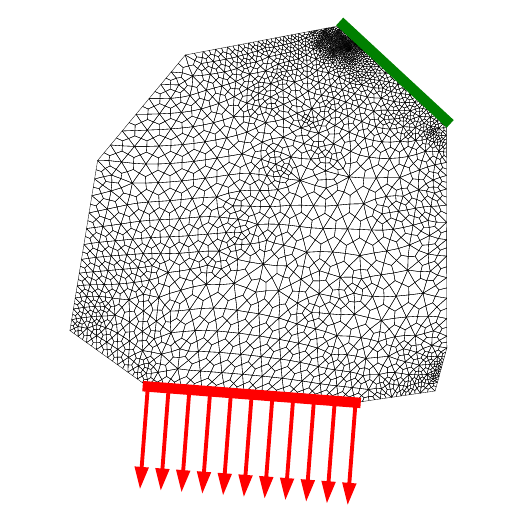}
		\caption{}
		\label{fig:convex_adaptiveMesh}
	\end{subfigure} 
	\begin{subfigure}[b]{0.24\textwidth}
		\includegraphics[width=1\textwidth]{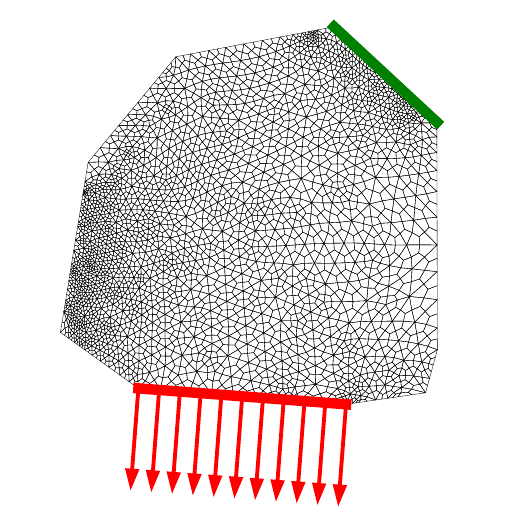}
		\caption{}
		\label{fig:convex_outputMesh}
	\end{subfigure} 
	\begin{subfigure}[b]{0.24\textwidth}
		\includegraphics[width=1\textwidth]{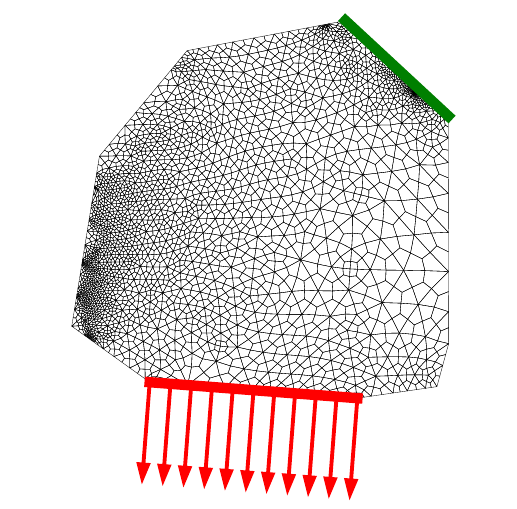}
		\caption{}
		\label{fig:convex_outputMesh_scaled}
	\end{subfigure} 
	\begin{subfigure}[b]{0.24\textwidth}
		\includegraphics[width=1\textwidth]{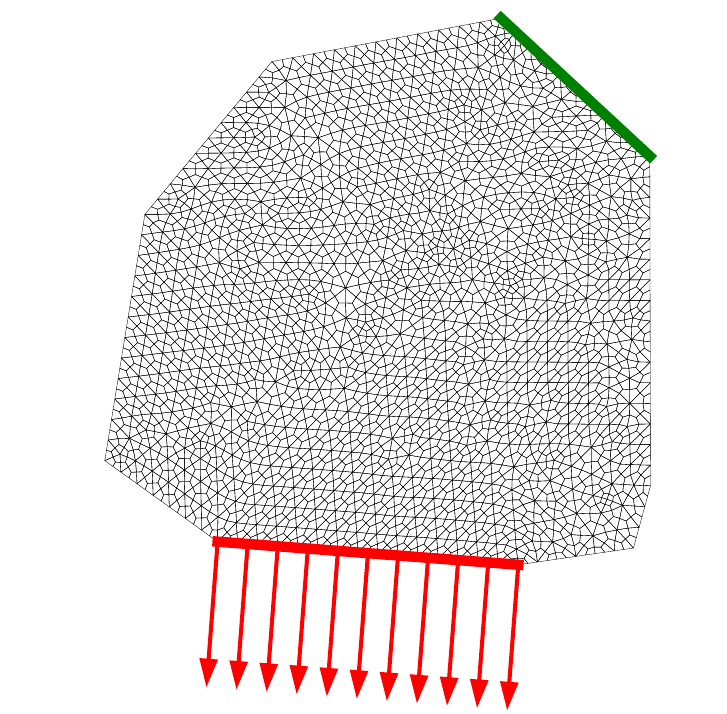}
		\caption{}
		\label{fig:convex_uniformMesh}
	\end{subfigure} 
	\begin{subfigure}[b]{0.2\textwidth}
		\includegraphics[width=1\textwidth]{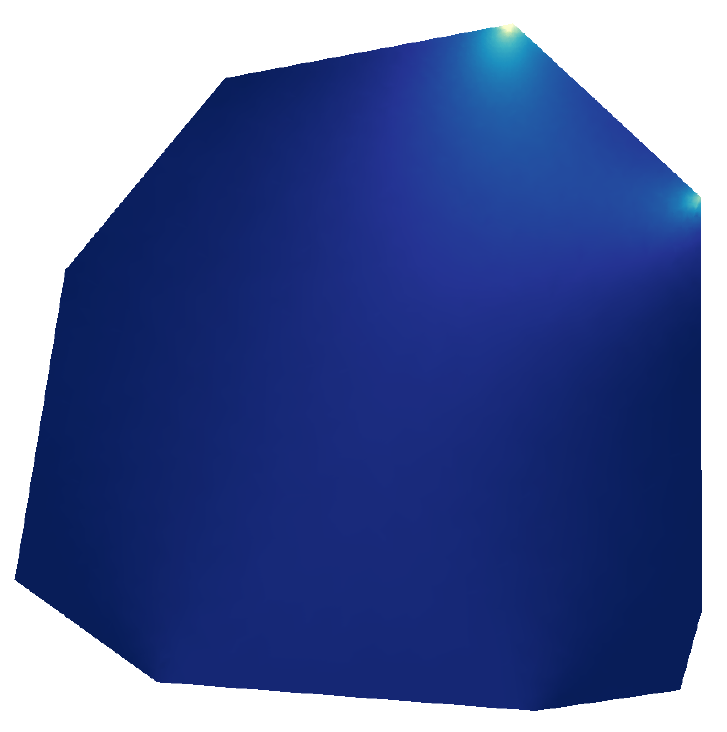}
		\caption{}
		\label{fig:convex_adaptiveVM}
	\end{subfigure} 
	\begin{subfigure}[b]{0.21\textwidth}
		\includegraphics[width=1\textwidth]{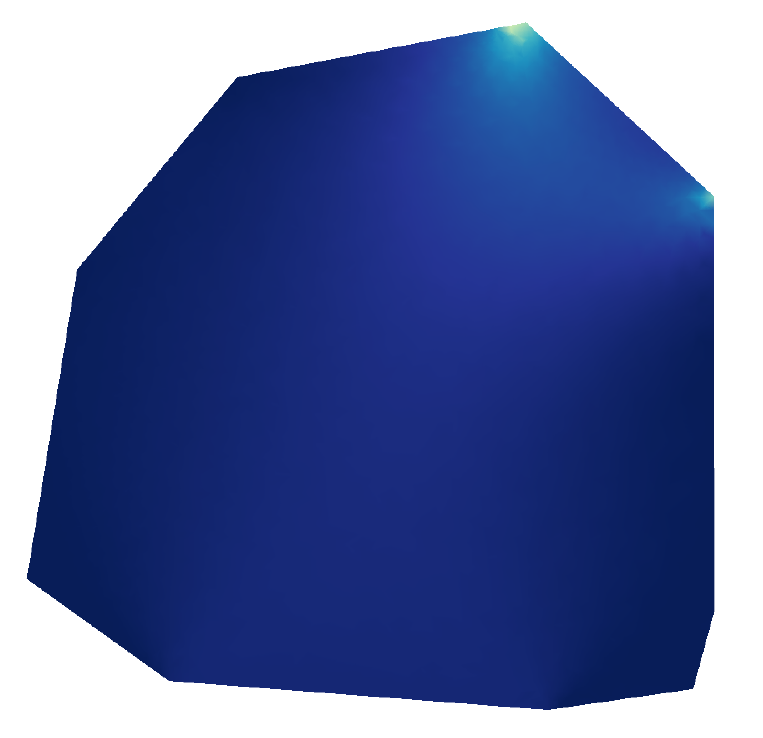}
		\caption{}
		\label{fig:convex_outputVM}
	\end{subfigure}
	\begin{subfigure}[b]{0.21\textwidth}
		\includegraphics[width=1\textwidth]{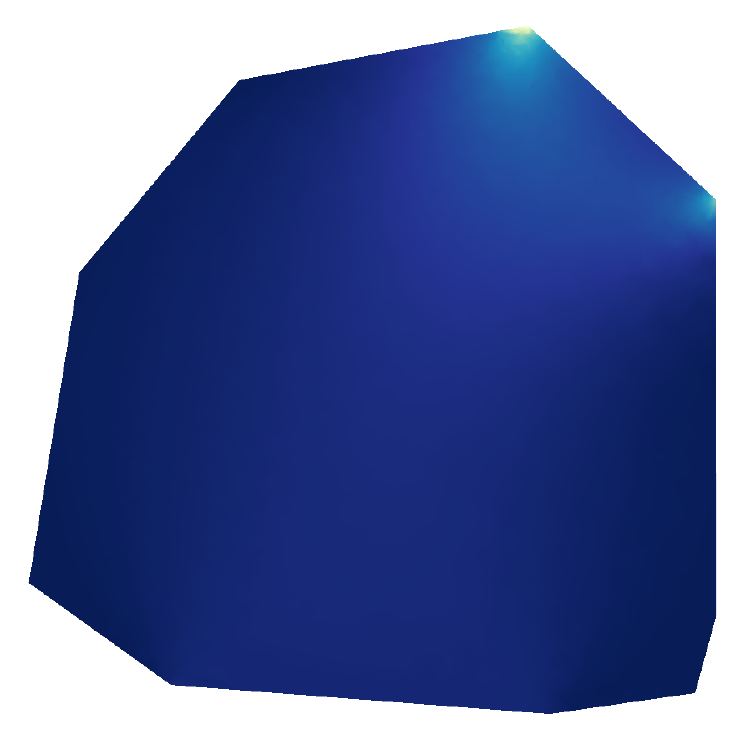}
		\caption{}
		\label{fig:convex_outputVM_scaled}
	\end{subfigure} 
	\begin{subfigure}[b]{0.27\textwidth}
		\includegraphics[width=1\textwidth]{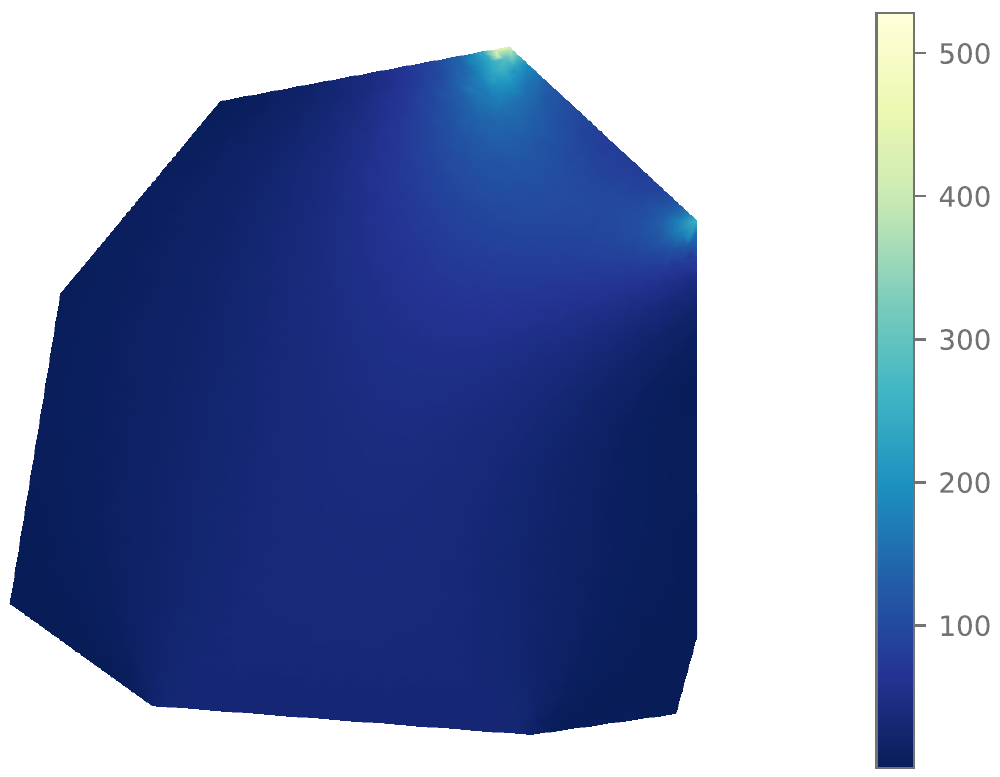}
		\caption{}
		\label{fig:convex_uniformVM}
	\end{subfigure}
	\caption{Example \emph{Model$(0,0,0)$}. First row: mesh from (a) adaptive refinement, (b) predicted output, (c)   post-processed output, and (d) uniform refinement. Second row: von Mises stress for (e) adaptive mesh, (f) predicted output mesh, (g) post-processed output mesh, and (h) uniform refinement mesh.}
	\label{fig:convex_exp}
\end{figure}

In the following, we examine the mesh discretization constructed based on the predicted output. We will start with showing the result from \emph {Model$(0,0,0)$} which is trained with a simple geometry complexity class. The adaptive mesh, the mesh constructed from the predicted output, and the uniform refinement mesh are shown in the first row of Figure \ref{fig:convex_exp}. The meshes are constructed such that they have approximately the same number of elements. The boundary marked with green color represents the Dirichlet boundary. The Neumann boundary is highlighted with red and the arrows denote the traction direction. The figures below each mesh show the corresponding von Mises stress. The relative error between the ground truth and predicted output is $\varepsilon_{rel}=0.019728$. 

It is shown in Figures \ref{fig:convex_adaptiveMesh} and \ref{fig:convex_outputMesh} that the predicted output has very similar mesh density distribution to the adaptive mesh (finer mesh size around the Dirichlet boundary and at the end points of the Neumann boundary). 
{However, it does not capture the high stress around the end points of the Dirichlet boundary well (See Figure \ref{fig:convex_outputVM}). On the other hand, the uniform mesh in Figure \ref{fig:convex_uniformMesh} shows a higher stress at the corresponding areas (See Figure \ref{fig:convex_outputVM}). The reason is that for certain geometries, the uniform mesh contains smaller quads around the corner for better shape approximation. Therefore it is able to capture the high stress if it occurs at those corners. Nonetheless, we can improve the predicted mesh by post-processing the network output. The output is scaled using a piecewise linear function such that for intensity values lower than 0.02, we scale the value by a factor of 0.25. The post-processed result is shown in Figure \ref{fig:convex_outputMesh_scaled}. It is shown in the corresponding von Mises stress distribution (Figure \ref{fig:convex_outputVM_scaled}) that the mesh better resolves the high stress region near the Dirichlet boundary. Meanwhile, this locally increased refinement also results in fewer and coarser elements in other areas, as a result of imposing a constraint on the number of elements.} The comparison of number of elements and relative error in energy norm ($e_{rel}$) are given in Table \ref{tab:convex_exp}. Note that the error is mostly concentrated in the region with high stresses near the fixed boundary.

\begin{table}[ht]
	\centering
	\begin{tabular}[t]{lccc}
		\hline
		Refinement & Number of elements& Relative error in the energy norm & Maximum von Mises stress\\
		\hline
		adaptive  & 3639 & 0.069448 & 520.839698\\
		predicted & 3633 & 0.111299 & 400.140106\\
		predicted (postprocessed)& 3651 & 0.100510 & 527.925733\\
		uniform & 3606 & 0.112577 & 440.274180\\
		\hline		
	\end{tabular}
	\caption{{Example \emph{Model$(0,0,0)$} : comparison of number of elements, relative error in the energy norm and maximum value of von Mises stress.}}
	\label{tab:convex_exp}
\end{table}%

\subsection{Example for geometric complexity class (1,0,0)}\label{sec:exp_100}
\begin{figure}[hbt!]
	\centering
	\begin{subfigure}[b]{0.3\textwidth}
		\includegraphics[width=1\textwidth]{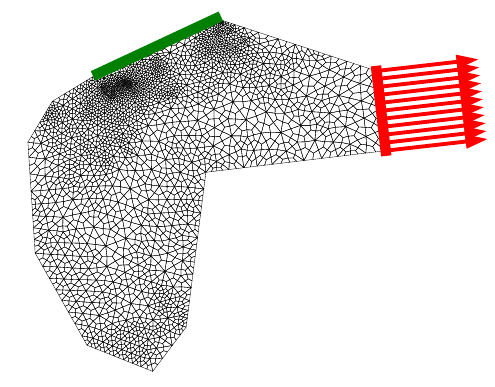}
		\caption{}
		\label{fig:concave_adaptiveMesh}
	\end{subfigure} 
	\begin{subfigure}[b]{0.3\textwidth}
		\includegraphics[width=1\textwidth]{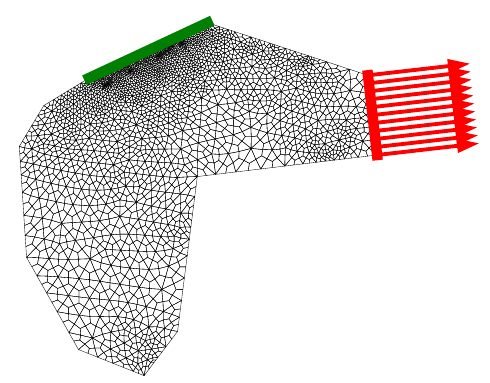}
		\caption{}
		\label{fig:concave_outputMesh}
	\end{subfigure}
	\begin{subfigure}[b]{0.3\textwidth}
		\includegraphics[width=1\textwidth]{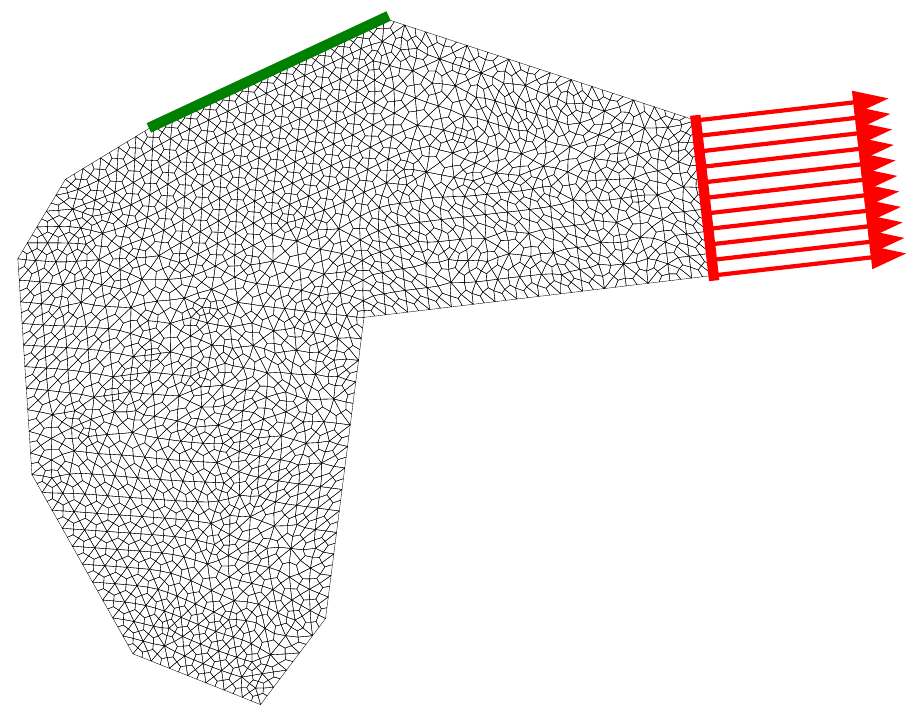}
		\caption{}
		\label{fig:concave_uniformMesh}
	\end{subfigure} 
	\begin{subfigure}[b]{0.25\textwidth}
		\includegraphics[width=1\textwidth]{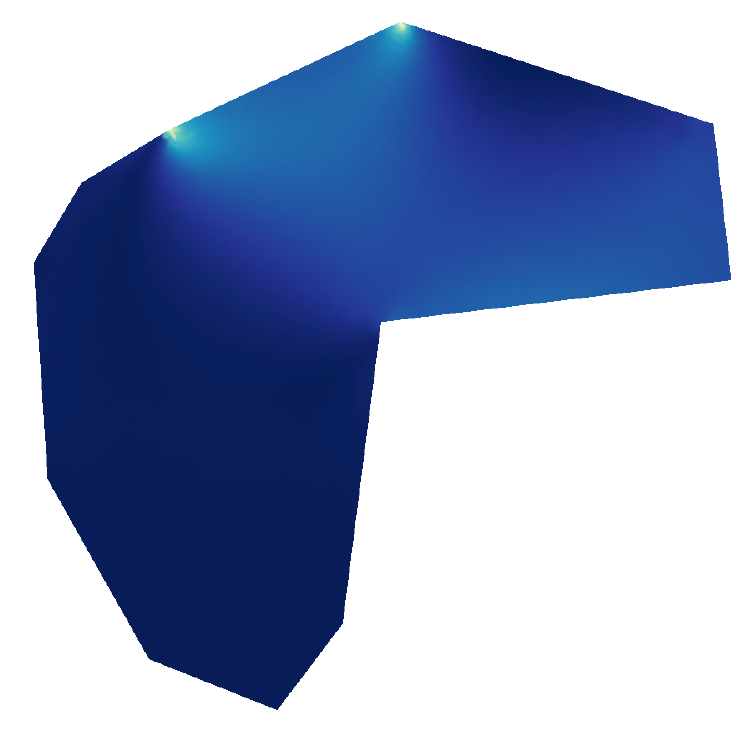}
		\caption{}
		\label{fig:concave_adaptiveVM}
	\end{subfigure} \quad\quad
	\begin{subfigure}[b]{0.25\textwidth}
		\includegraphics[width=1\textwidth]{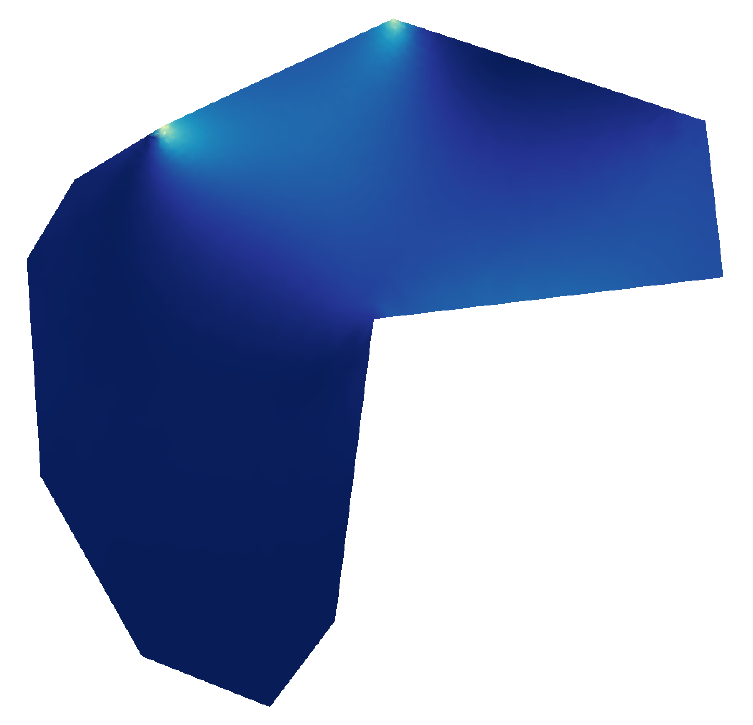}
		\caption{}
		\label{fig:concave_outputVM}
	\end{subfigure} \quad\quad
	\begin{subfigure}[b]{0.325\textwidth}
		\includegraphics[width=1\textwidth]{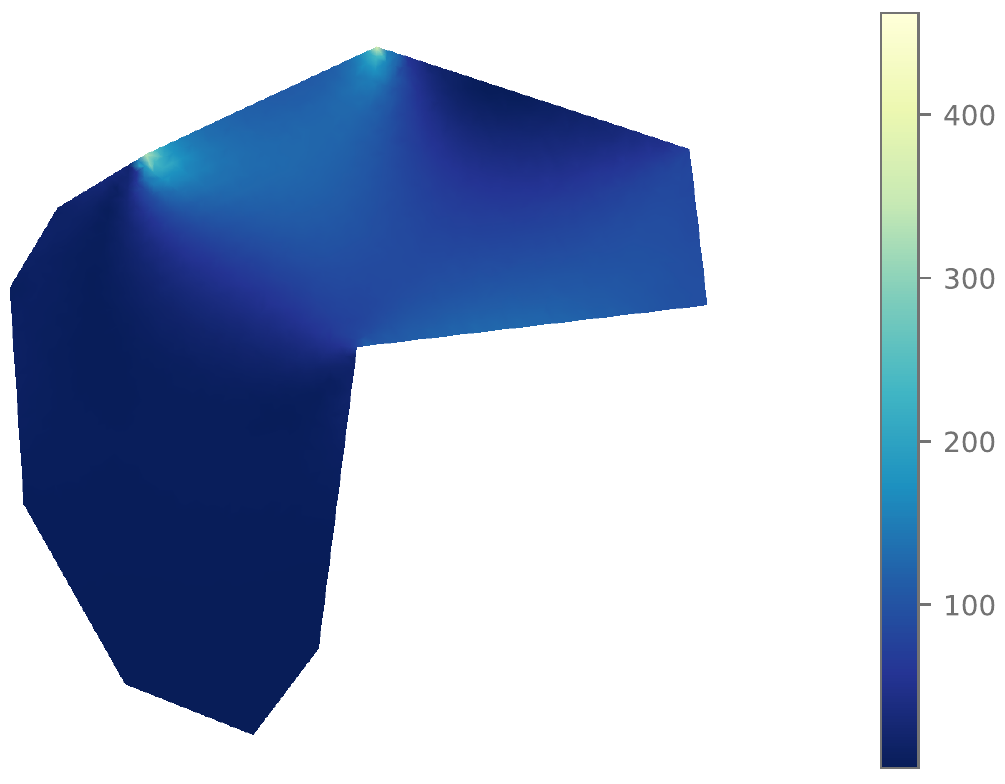}
		\caption{}
		\label{fig:concave_uniformVM}
	\end{subfigure}
	\caption{\emph{Model$(1,0,0)$} example. First row: mesh from (a) adaptive refinement, (b) predicted output, and (c) uniform refinement. Second row: von Mises stress for (d) adaptive mesh, (e) predicted output mesh, and (f) uniform refinement mesh.}
	\label{fig:concave_exp}
\end{figure}

In this example shows the result from \emph{Model$(1,0,0)$}. The comparison of adaptive, predicted and fine uniform meshes and their von Mises stress are given in Figure \ref{fig:concave_exp}. For this example, the relative error between predicted output and ground truth is given by $\varepsilon_{rel}= 0.009530$. The predicted and adaptive meshes match fairly well in the domain especially around the Dirichlet boundary. The relative errors in energy norm for the adaptive, predicted, and uniform meshes as listed in Table \ref{tab:concave_exp} also demonstrate the capability of proposed method for mesh density prediction. 

\begin{table}[hbt!]
	\centering
	\begin{tabular}[t]{lccc}
		\hline
		Refinement & Number of elements& Relative error in the energy norm & Maximum value of von Mises stress\\
		\hline
		adaptive  & 3447 & 0.078086 & 462.201620\\
		predicted &3432 & 0.088255 & 428.388343\\
		uniform & 3438 & 0.106043 &  373.382649\\
		\hline		
	\end{tabular}
	\caption{{\emph{Model$(1,0,0)$} example: comparison of number of elements, relative error in the energy norm, and  maximum value of von Mises stress.}}
	\label{tab:concave_exp}
\end{table}%

\subsection{Example for geometric complexity class (1,2,0)}\label{sec:exp_120}
\begin{figure}[hbt!]
	\centering	
	\vspace{-0.5cm}
	\begin{subfigure}[b]{0.25\textwidth}
		\includegraphics[width=1\textwidth]{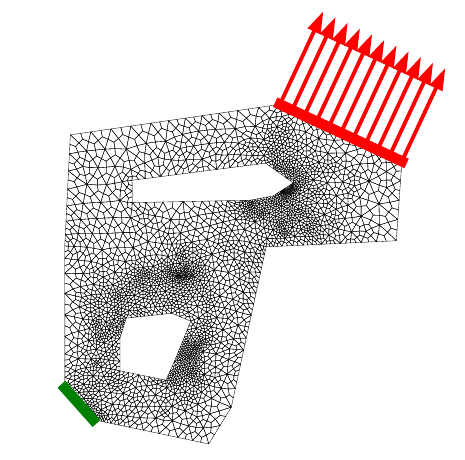}
		\caption{}
		\label{fig:concave_void_adaptiveMesh} 
	\end{subfigure} \quad\quad
	\begin{subfigure}[b]{0.25\textwidth}
		\includegraphics[width=1\textwidth]{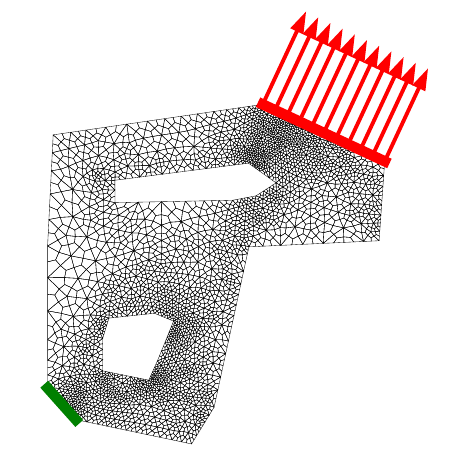}
		\caption{}
		\label{fig:concave_void_outputMesh}
	\end{subfigure} \quad\quad
	\begin{subfigure}[b]{0.25\textwidth}
		\includegraphics[width=1\textwidth]{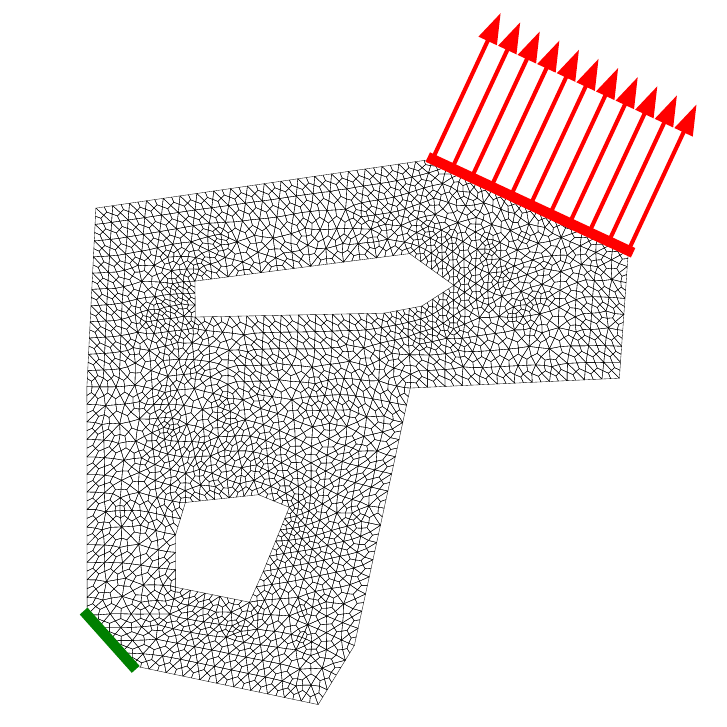}
		\caption{}
		\label{fig:concave_void_uniformMesh}
	\end{subfigure} 
	\begin{subfigure}[b]{0.2\textwidth}
		\includegraphics[width=1\textwidth]{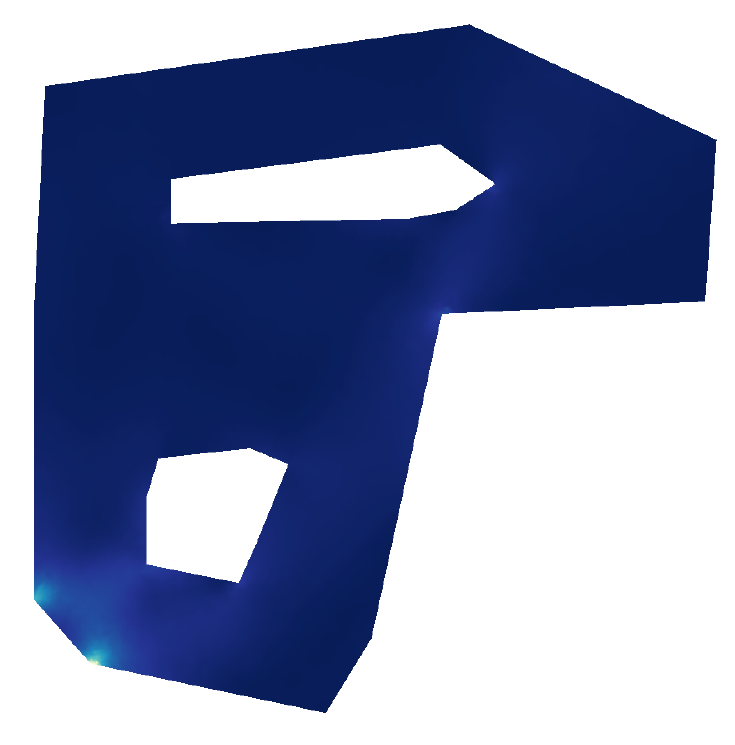}
		\caption{}
		\label{fig:concave_void_adaptiveVM}
	\end{subfigure}\quad\quad
	\begin{subfigure}[b]{0.2\textwidth}
		\includegraphics[width=1\textwidth]{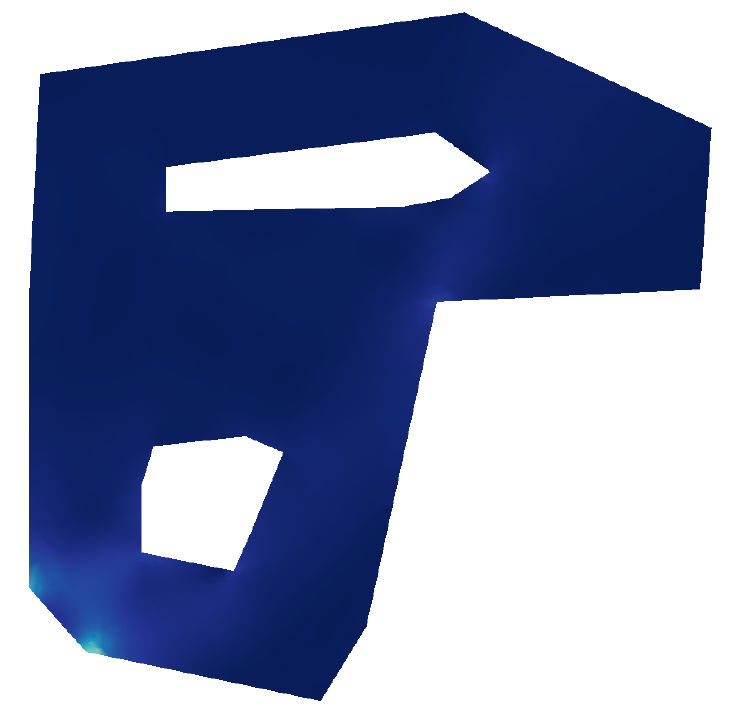}
		\caption{}
		\label{fig:concave_void_outputVM}
	\end{subfigure}\quad\quad
	\begin{subfigure}[b]{0.3\textwidth}
		\includegraphics[width=1\textwidth]{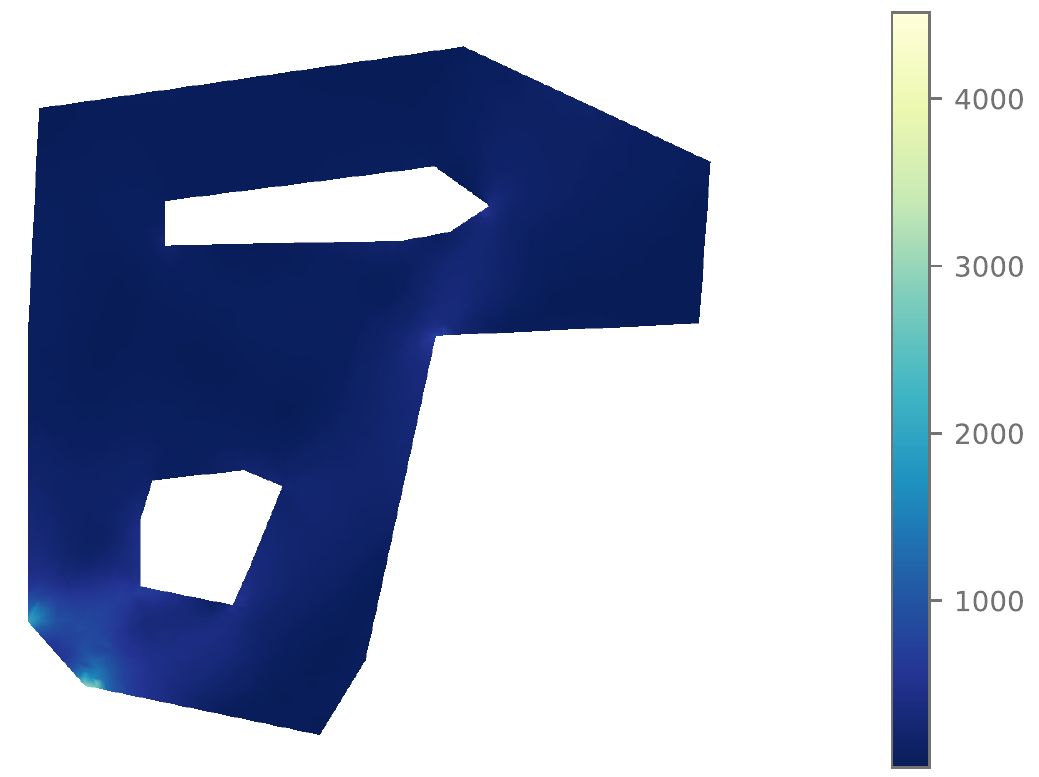}
		\caption{}
		\label{fig:concave_void_uniformVM}
	\end{subfigure}
	\begin{subfigure}[b]{0.25\textwidth}	
		\includegraphics[width=1\textwidth]{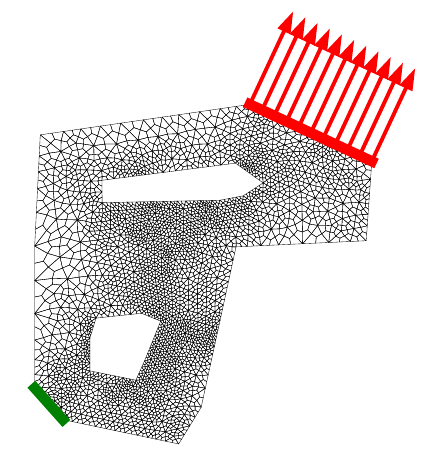}
		\caption{}
		\label{fig:modelI_outputMesh}
	\end{subfigure} \quad
	\begin{subfigure}[b]{0.25\textwidth}
		\includegraphics[width=1\textwidth]{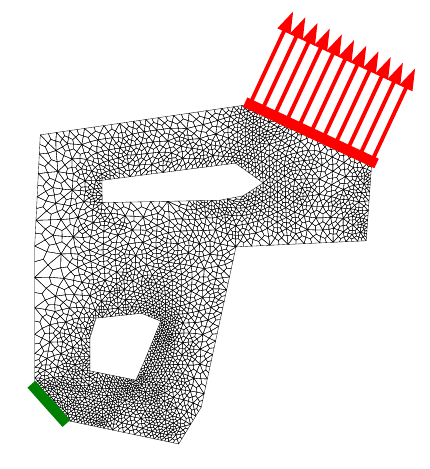}
		\caption{}
		\label{fig:modelII_outputMesh}
	\end{subfigure} \\
	\begin{subfigure}[b]{0.2\textwidth}
		\includegraphics[width=1\textwidth]{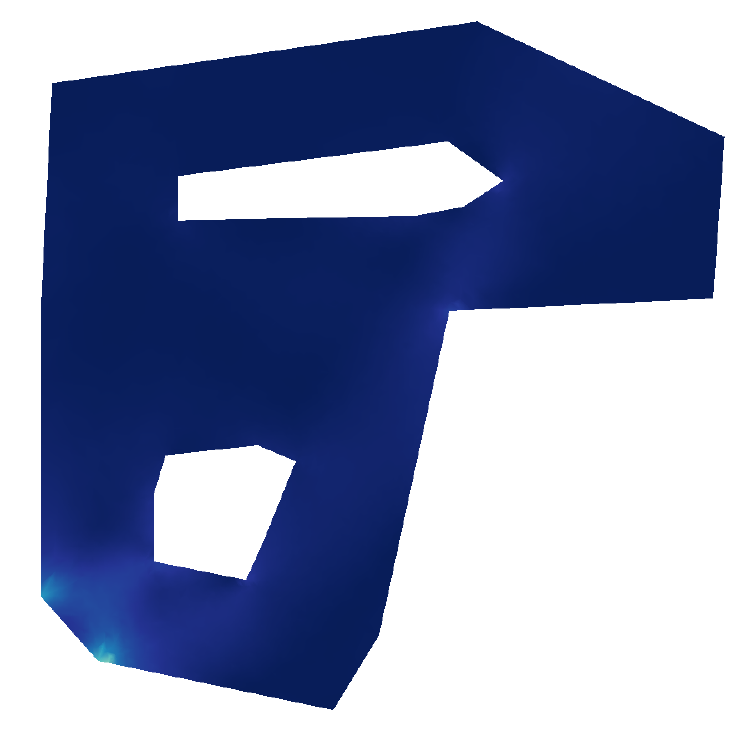}
		\caption{}
		\label{fig:modelI_outputVM}
	\end{subfigure} \quad \quad
	\begin{subfigure}[b]{0.2\textwidth}
		\includegraphics[width=1\textwidth]{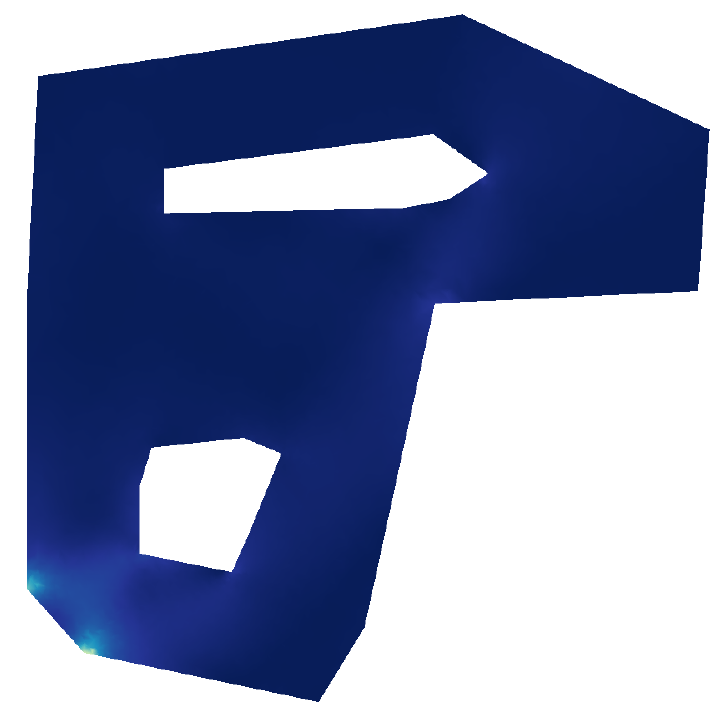}
		\caption{}
		\label{fig:modelII_outputVM}
	\end{subfigure}
	\caption{Comparison of results from trained networks. First row: (a) adaptive mesh, (b) output mesh, and (c) uniform refinement mesh. Second row: von Mises stress computed on (d) adaptive mesh, (e) output mesh and (f) uniform refinement mesh.  Third and fourth rows: mesh constructed using predicted output from (g) \emph{Model RT} and (h) \emph{Model CT} and the corresponding von Mises stress.}
	\label{fig:concave_void_exp}
\end{figure}

Here we show an example with complex geometry. It is a non-convex domain with two voids. The comparison of adaptive mesh, uniform mesh and predicted output from \emph{Model$(1,k,0)$} are shown in the first row of Figure \ref{fig:concave_void_exp}. The figures below the meshes are the corresponding von Mises stresses.
In this example, we also compare predicted outputs from \emph{Model RT} and \emph{Model CT}. The resulting meshes and von Mises stress are shown in Figures \ref{fig:modelI_outputMesh} - \ref{fig:modelII_outputVM}. The predicted output from \emph{Model $(1,k,0)$} has the relative error $\varepsilon_{rel}= 0.014139$. \emph{Model CT} results in a slightly lower relative error $\varepsilon_{rel}= 0.012991$. An interesting observation is that the results from \emph{Model RT} ($\varepsilon_{rel}=0.014128$) are very close to \emph{Model$(1,k,0)$}. Geometric complexity class $(1,2,0)$ is not included in the training data set for \emph{Model RT}, but the Unet is able to approximately resemble the higher mesh density around the holes. This demonstrates the generalization capability of the network, which is an advantage of the proposed data-driven mesh density prediction method.  

The number of elements and relative errors in energy norm for all the meshes are listed in Table \ref{tab:genus_exp}. We can observe that when the input has complex geometry, the complete training model (\emph{Model CT}) can indeed provide a better discretization results. 

\begin{table}[!htbp]
	\centering
	\begin{tabular}[t]{lccc}
		\hline
		Refinement & Number of elements& Relative error in the energy norm & Maximum value of von Mises stress\\
		\hline
		adaptive  & 4038 & 0.171800 & 4513.818979\\
		predicted (\emph{Model$(1,k,0)$}) & 4008 &0.205166 & 2850.473993\\	
		predicted  (\emph{Model RT}) & 4020 &0.208462 & 2766.330910\\	
		predicted  (\emph{Model CT}) & 4008 &0.176951 & 3442.309417\\
		uniform  & 4044 &0.221010 & 2869.988917\\	
		\hline		
	\end{tabular}
	\caption{{Example for geometric complexity class $(1,2,0)$: number of elements, relative error in the energy norm, and maximum value of von Mises stress.}}
	\label{tab:genus_exp}
\end{table}%

\subsection{Evaluation of models} \label{sec:eval_of_model}

\begin{figure}[hbt!]
	\centering
	\begin{subfigure}[b]{0.45\textwidth}
		\includegraphics[width=1\textwidth]{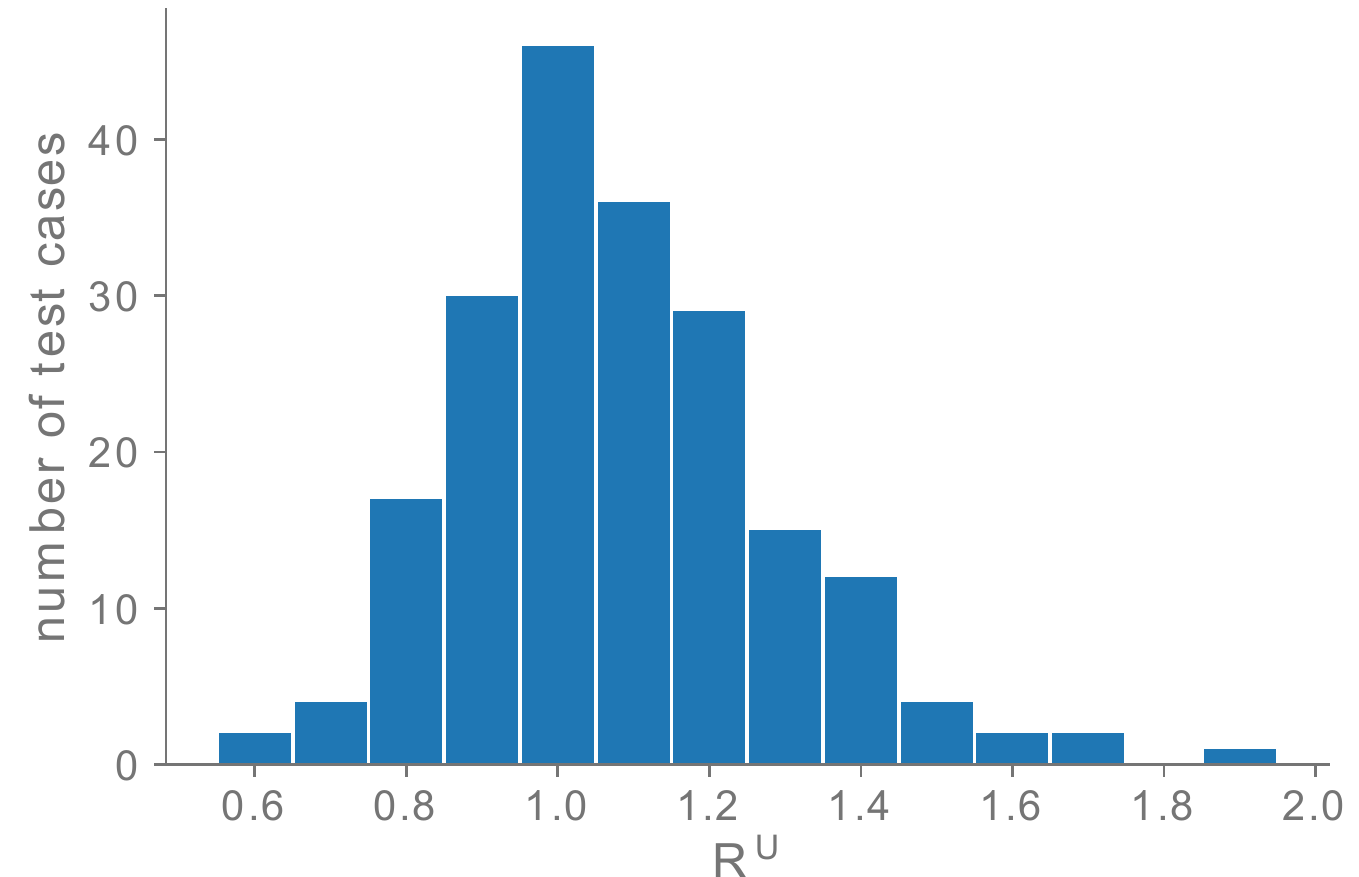}
		\caption{}
		\label{fig:uniform_output_error_convex}
	\end{subfigure} 
	\begin{subfigure}[b]{0.45\textwidth}
		\includegraphics[width=1\textwidth]{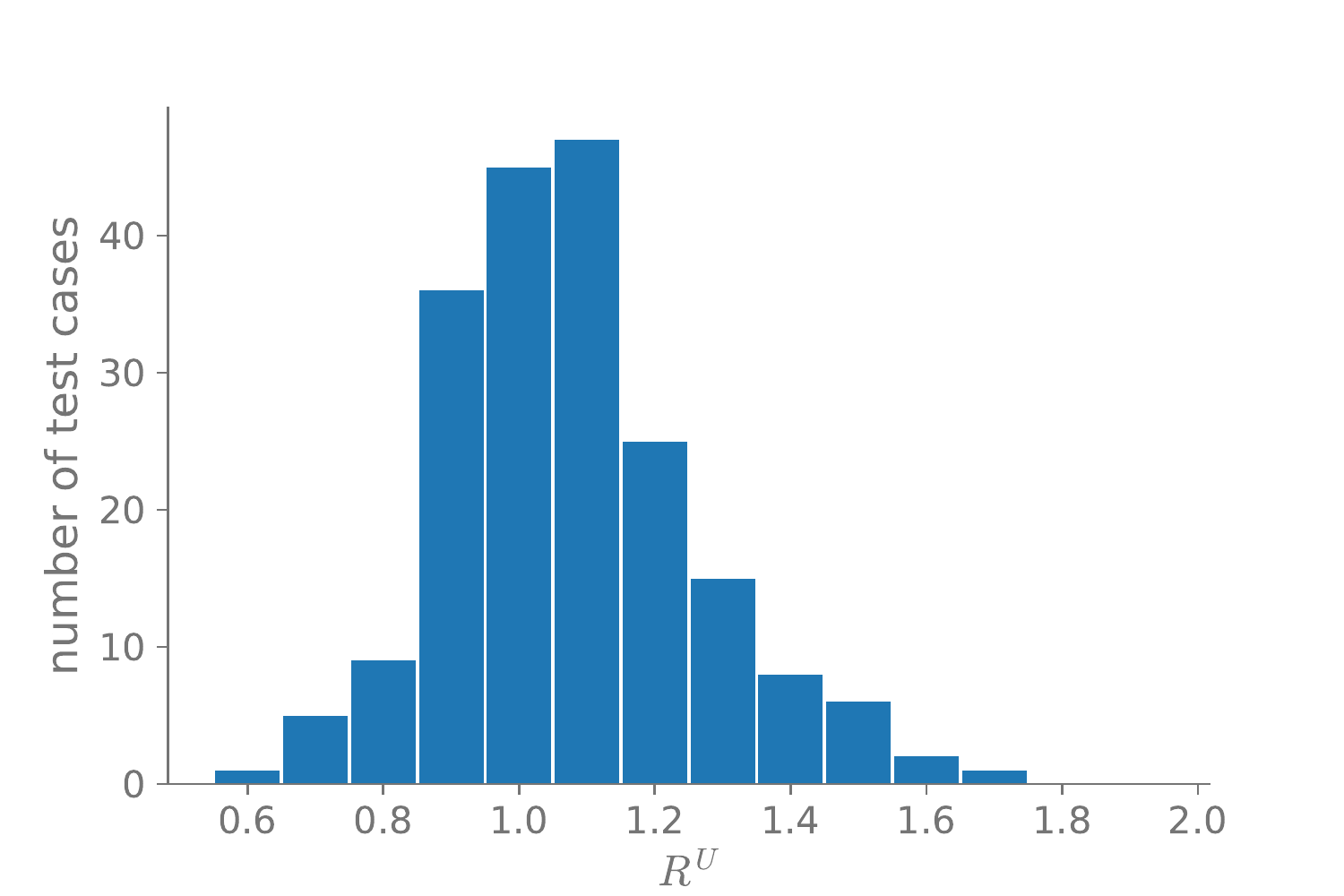}
		\caption{}
		\label{fig:uniform_output_error_concave}
	\end{subfigure} 
	\begin{subfigure}[b]{0.45\textwidth}
		\includegraphics[width=1\textwidth]{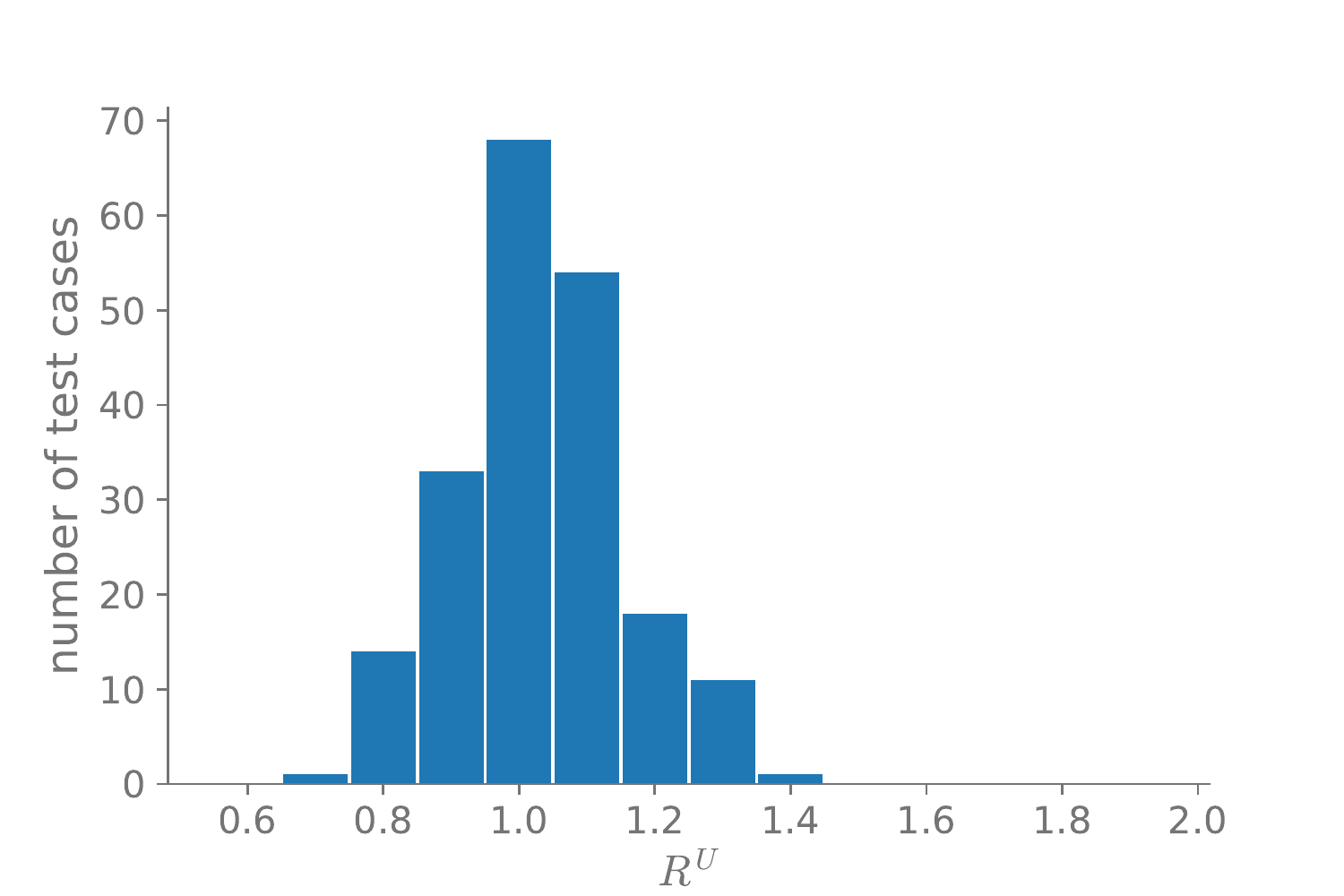}
		\caption{}
		\label{fig:uniform_output_error_convex_void}
	\end{subfigure} 
	\begin{subfigure}[b]{0.45\textwidth}
		\includegraphics[width=1\textwidth]{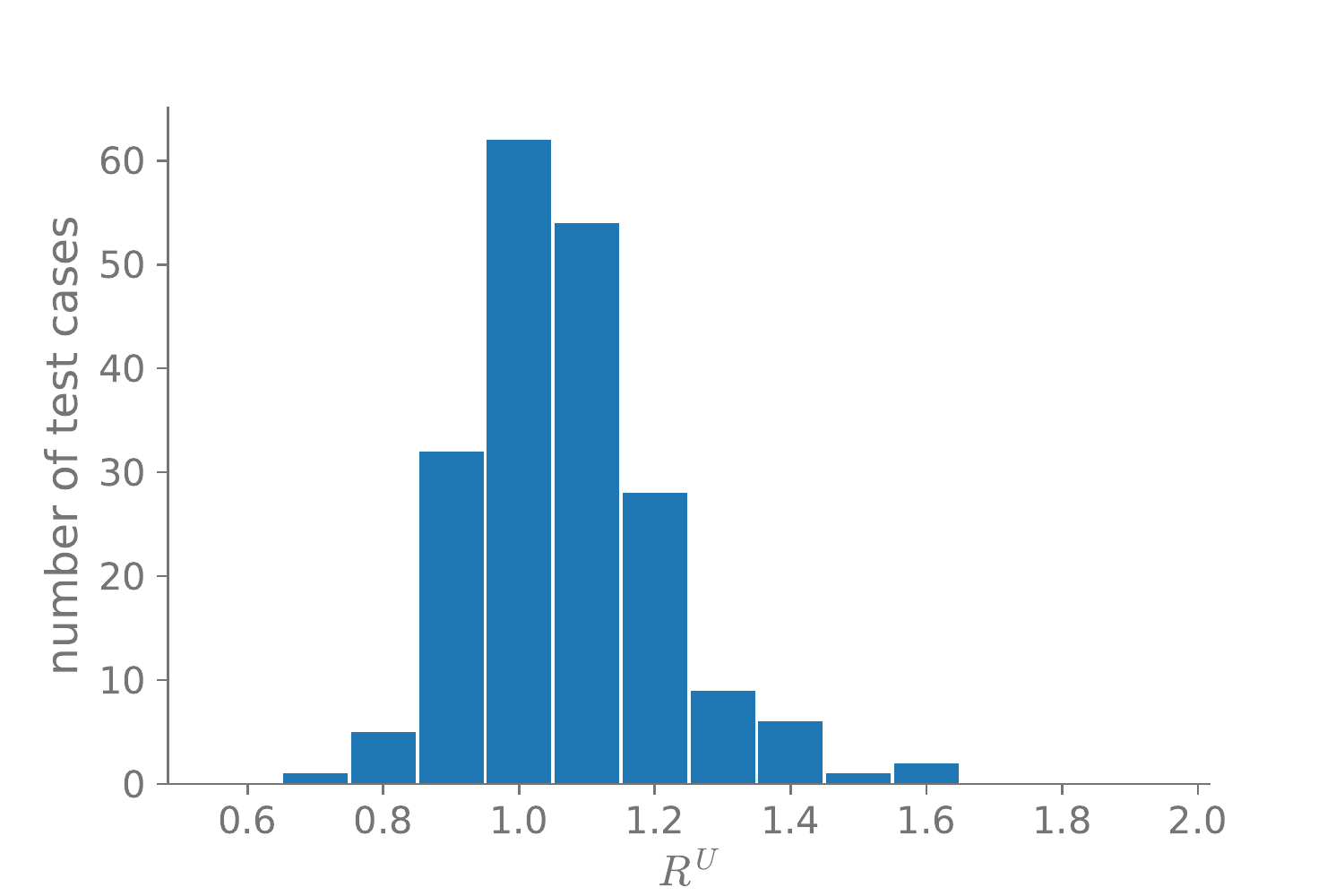}
		\caption{}
		\label{fig:uniform_output_error_concave_void}
	\end{subfigure} 
	\caption{Distribution of $R^{U}$ for (a) \emph{Model$(0,0,0)$} , (b) \emph{Model$(1,0,0)$}, (c) \emph{Model$(0,k,0)$}, and (d) \emph{Model$(1,k,0)$}, $k\in\{1,2\}$.}
	\label{fig:uniform_output_error}
\end{figure}

For a given predicted output, we construct a uniform mesh with approximately the same number of degrees of freedom, such that: 
$\frac{\|DOFs^{U}-DOFs^{P}\|}{DOFs^{P}}$<0.05 where $DOFs^{U}$ and $DOFs^{P}$ represent the degrees of freedom of the uniform and predicted meshes. We then quantify the solution computed on the predicted mesh over the uniform mesh by using the ratio of errors in the energy norm as follows:
\begin{equation}	
	R^{U}=\frac{ e_{rel}^{U}}{ e_{rel}^{P}}
	\label{eq:adapt_vs_pred}
\end{equation} 
where $e_{rel}^{U}$ and $e_{rel}^{P}$ refer to the relative error in the energy norm \eqref{eq:energy_norm} of the uniform and predicted mesh, respectively. 
The rational functions $R^{U}$ is interpreted as:  
\begin{itemize} 
	\item $R^{U}$ $\approx$ 1 : both meshes have similar analysis accuracy.
	\item $R^{U}$ $<$ 1 : the uniform mesh possesses better analysis accuracy.
	\item $R^{U}$ $>$ 1 : the predicted mesh has better analysis accuracy.
\end{itemize}
The distributions of $R^{U}$ for \emph{Model$(0,0,0)$}, \emph{Model$(1,0,0)$}, \emph{Model$(0,k,0)$}, and \emph{Model$(1,k,0)$} are shown in Figure \ref{fig:uniform_output_error}. Each histogram shows the results of $R^{U}$ computed on $100$ sets randomly selected test data from the corresponding geometric complexity classes. There are some cases where both meshes have similar accuracy, or the uniform mesh has better accuracy. However, the right-skewed distribution indicate that the predicted mesh performed better in most of the cases.

\begin{figure}[hbt!]
	\centering
	\begin{subfigure}[b]{0.45\textwidth}
		\includegraphics[width=1\textwidth]{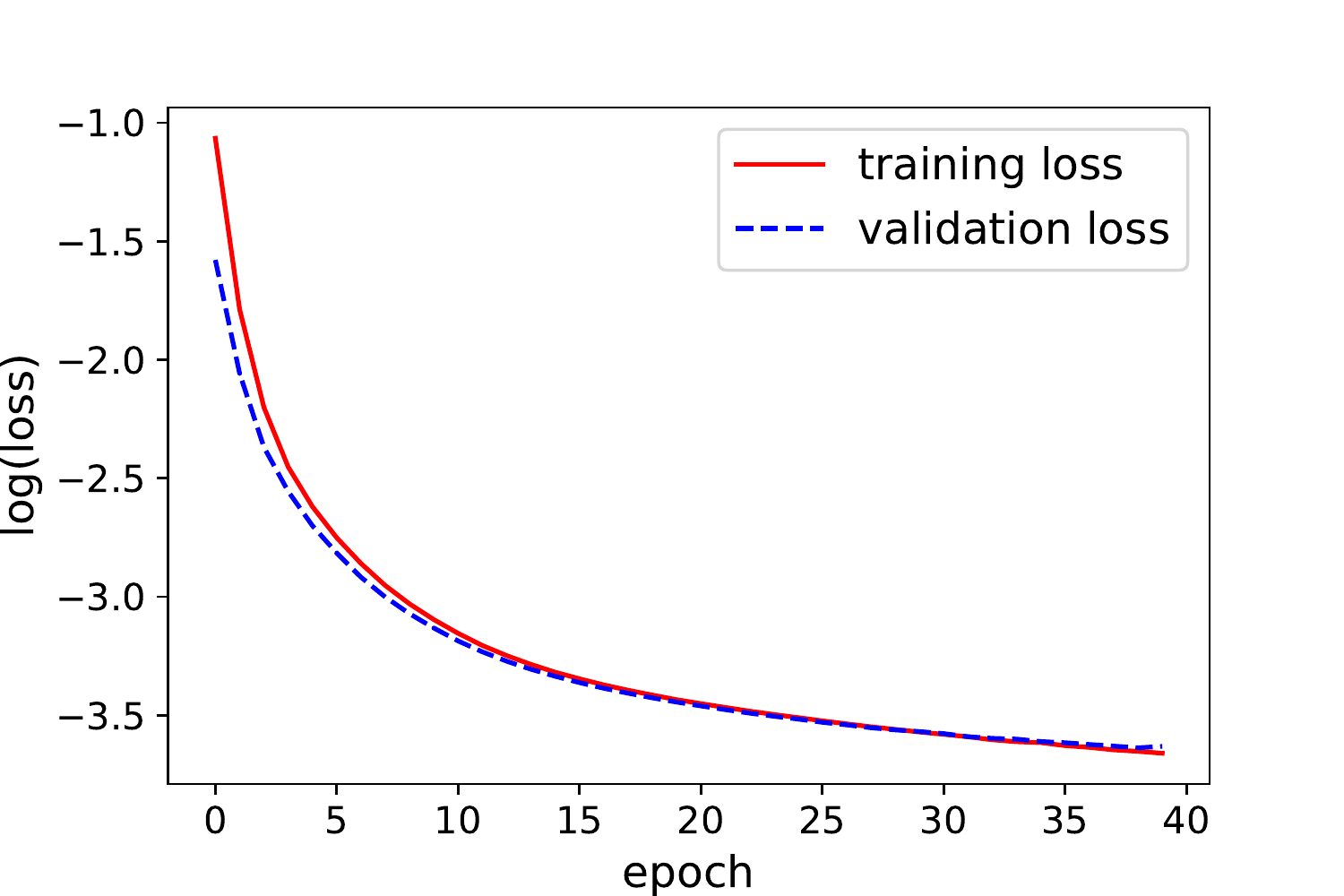}
		\caption{Training log for \emph{Model$(0,0,0)$}.}
		\label{fig:training_log_convex}
	\end{subfigure}	\quad\quad
	\begin{subfigure}[b]{0.45\textwidth}
		\includegraphics[width=1\textwidth]{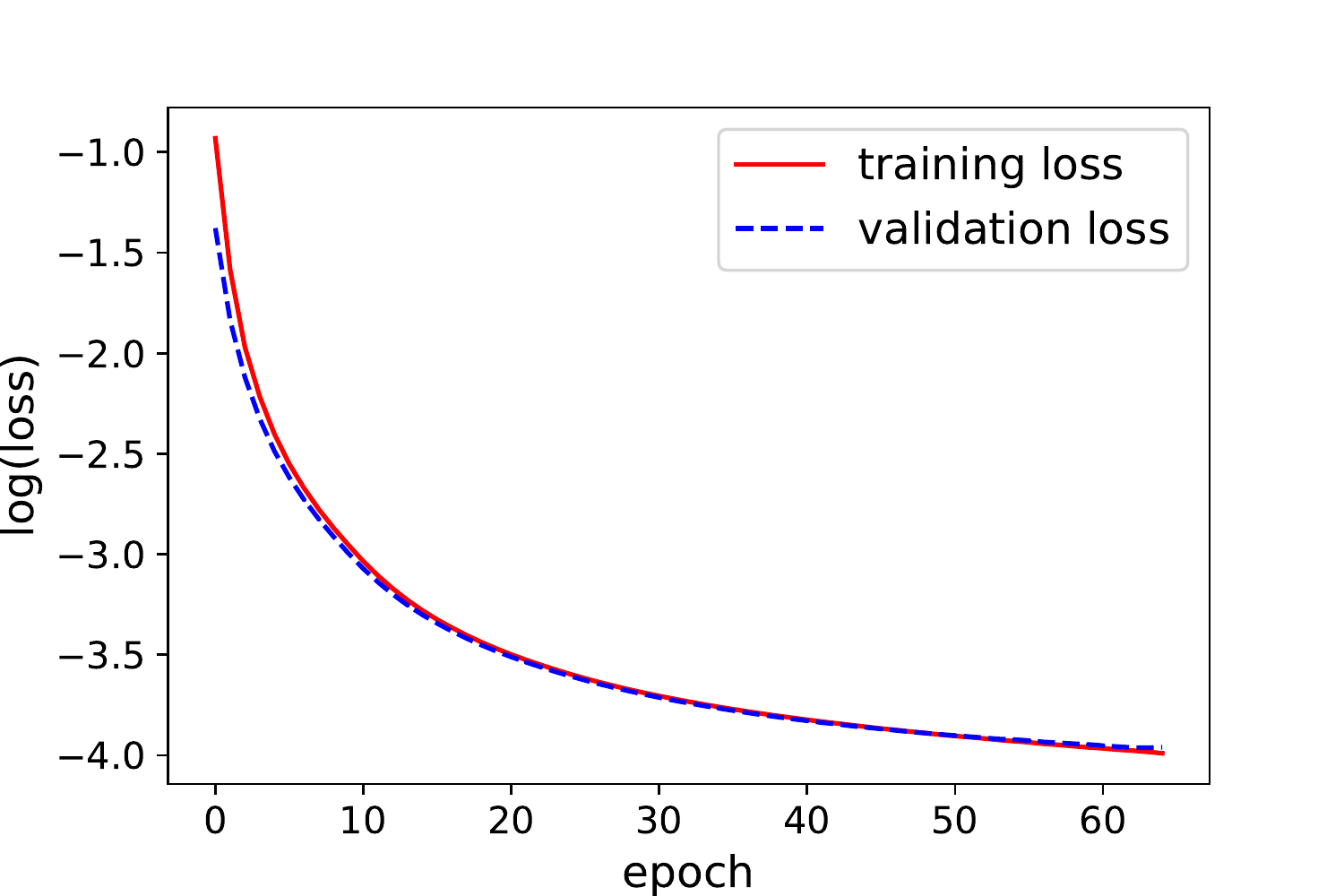}
		\caption{Training log for \emph{Model$(1,0,0)$}.}
		\label{fig:training_log_concave}
	\end{subfigure} 
	\begin{subfigure}[b]{0.45\textwidth}
		\includegraphics[width=1\textwidth]{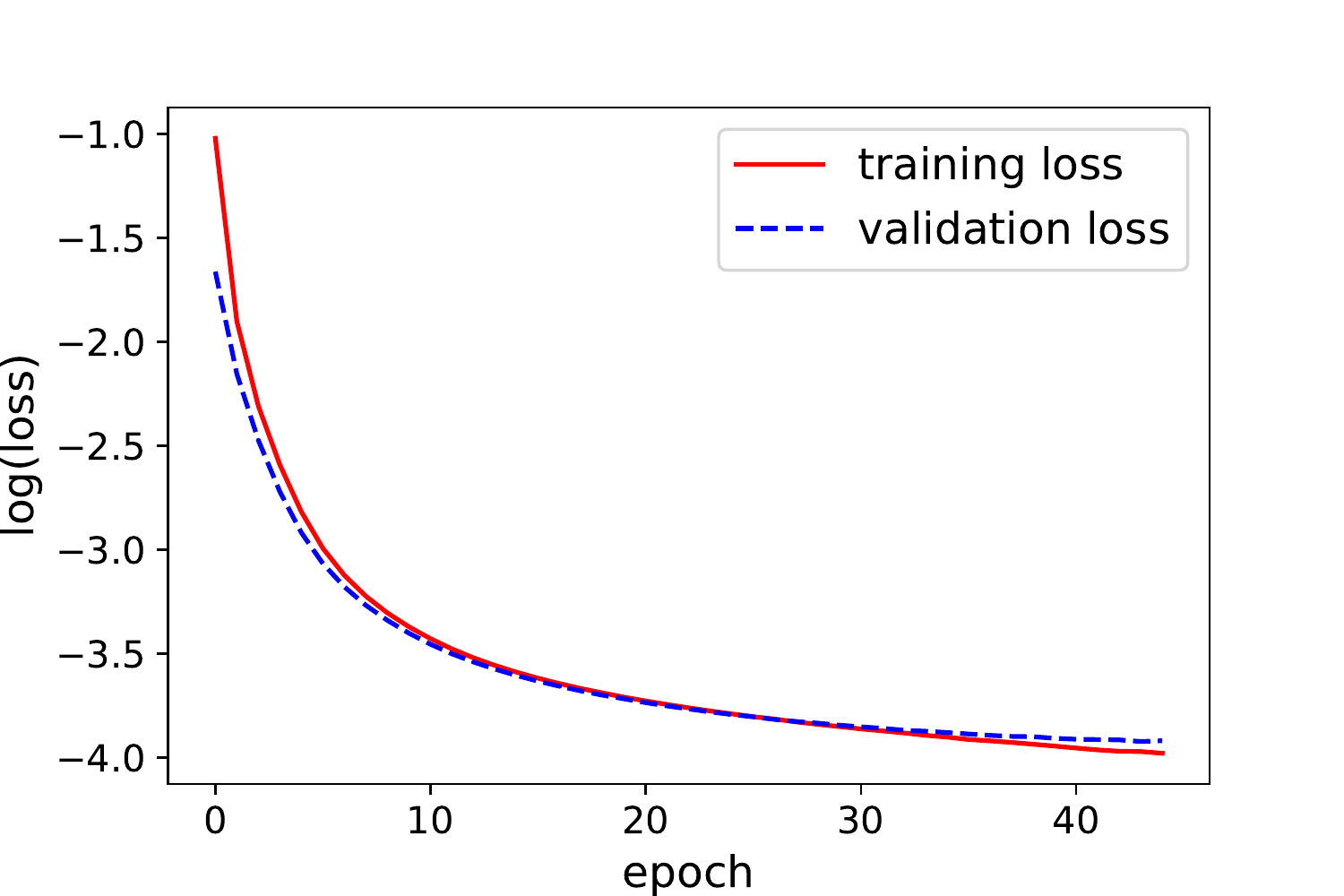}
		\caption{Training log for \emph{Model$(0,k,0)$}}
		\label{fig:training_log_convex_Nvoids}
	\end{subfigure} \quad\quad	
	\begin{subfigure}[b]{0.45\textwidth}
		\includegraphics[width=1\textwidth]{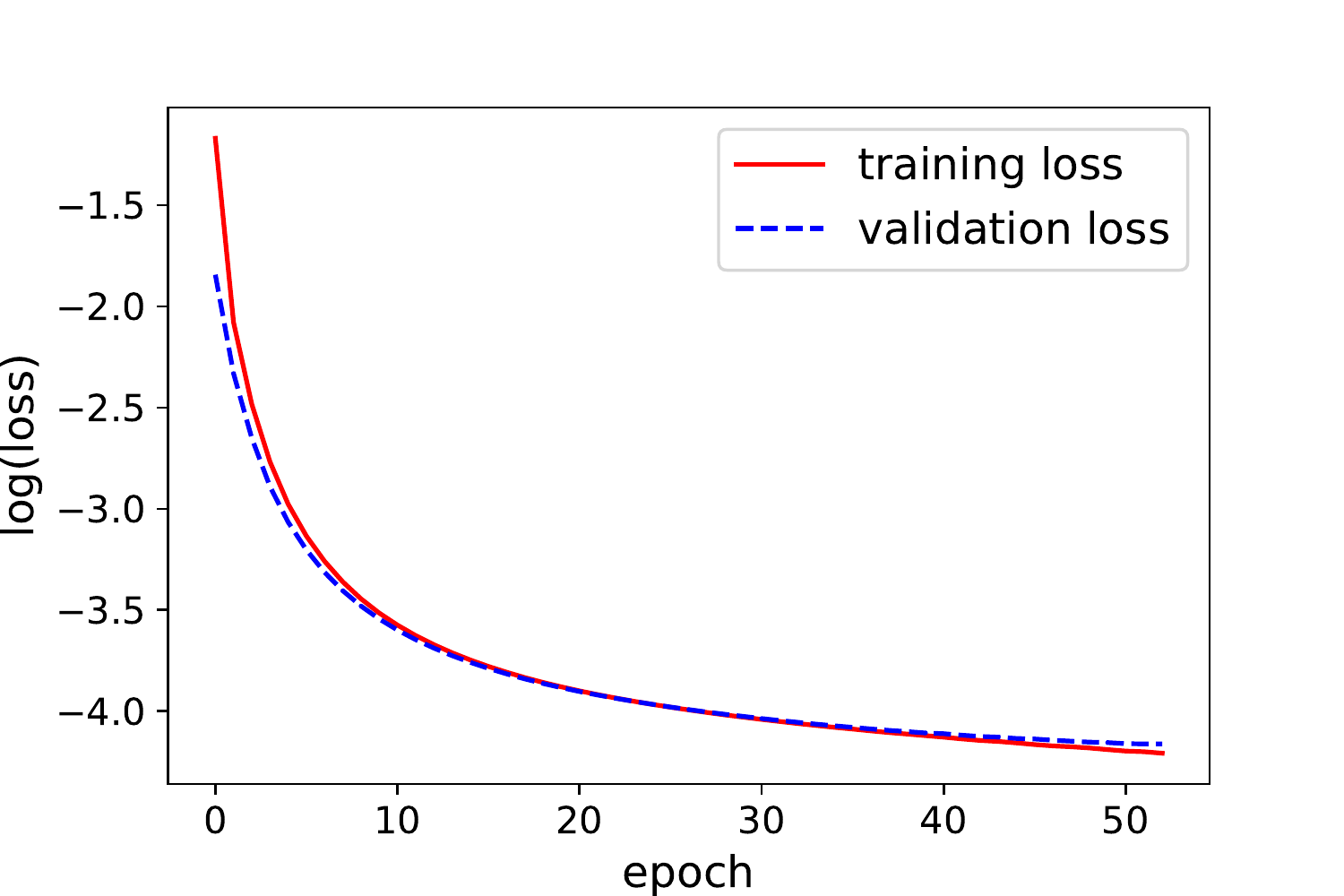}
		\caption{Training log for \emph{Model$(1,k,0)$}}
		\label{fig:training_log_concave_Nvoids}
	\end{subfigure} 
	\begin{subfigure}[b]{0.45\textwidth}
		\includegraphics[width=1\textwidth]{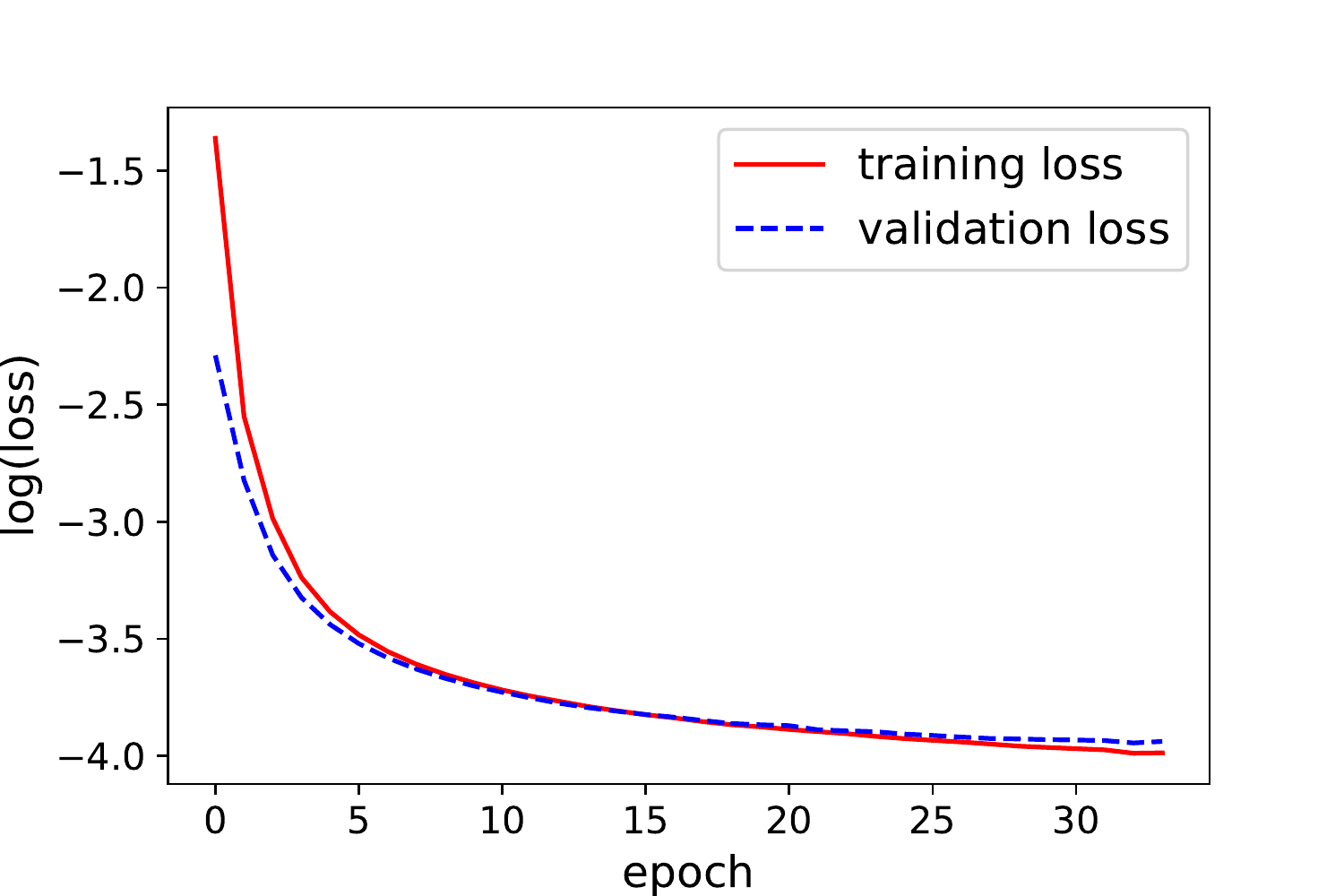}
		\caption{Training log for \emph{Model RT}}
		\label{fig:training_log_modelI}
	\end{subfigure} \quad\quad
	\begin{subfigure}[b]{0.45\textwidth}
		\includegraphics[width=1\textwidth]{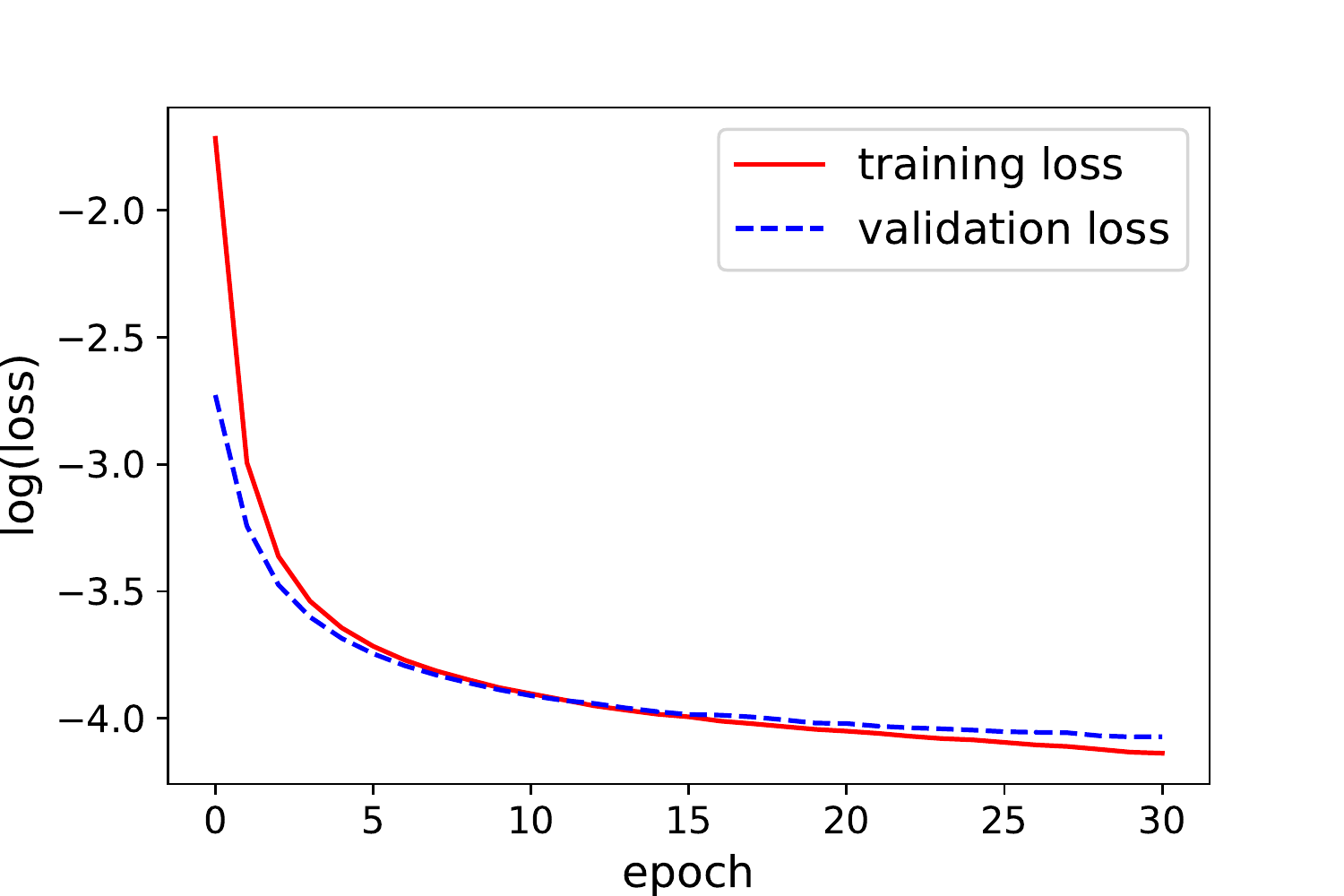}
		\caption{Training log for \emph{Model CT}}
		\label{fig:training_log_modelII}
	\end{subfigure}
	\caption{Training log for the six models trained.}
	\label{fig:training_log}
\end{figure}

Finally, the training log for all the models that we discussed above is given in Figure \ref{fig:training_log}. In the figures, the red line and blue dashes represent the training log and validation log respectively.
Note that we used the U-net architecture as discussed in Section \ref{sec:Unet} to train different group of data (i.e. different geometry complexity classes). 
In order to avoid over-fitting, the training process is terminated when validation loss started to increase. It is shown in the figures that all the training process converge at a very similar rate.

\begin{figure}[hbt!]
	\centering
	\begin{subfigure}[b]{0.25\textwidth}
		\includegraphics[width=1\textwidth]{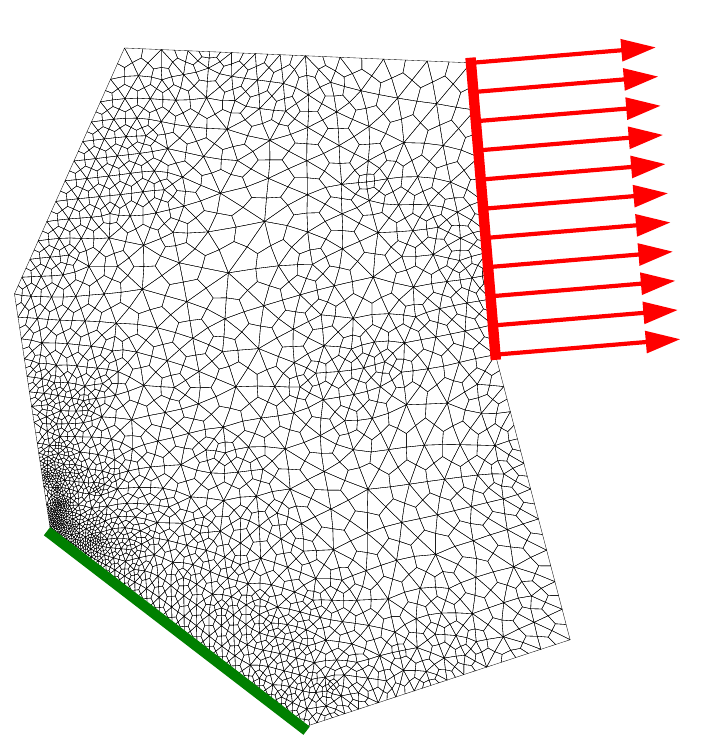}
		\caption{}
		\label{fig:meshNet_exp_adapt}
	\end{subfigure} \quad\quad\quad
	\begin{subfigure}[b]{0.25\textwidth}
		\includegraphics[width=1\textwidth]{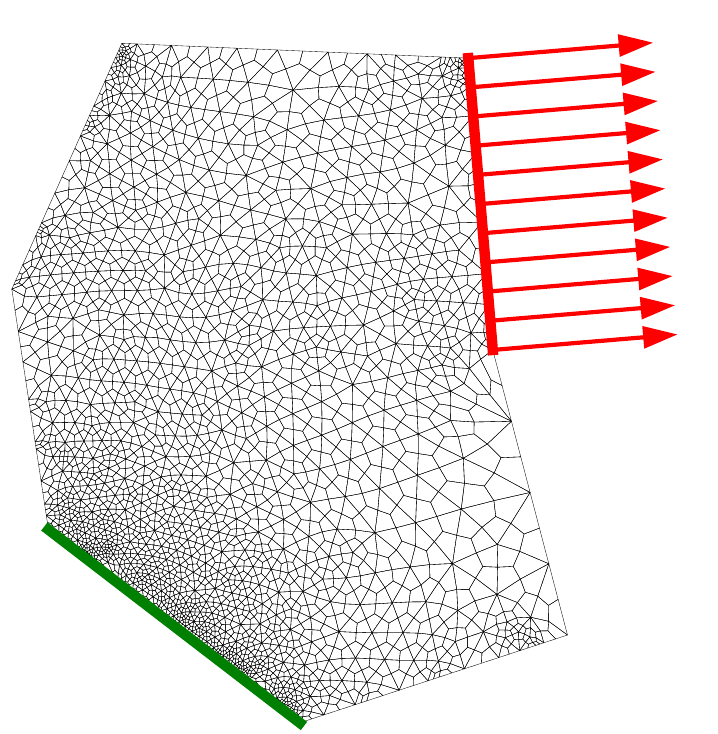}
		\caption{}
		\label{fig:meshNet_exp_output}
	\end{subfigure} \quad\quad\quad
	\begin{subfigure}[b]{0.25\textwidth}
		\includegraphics[width=1\textwidth]{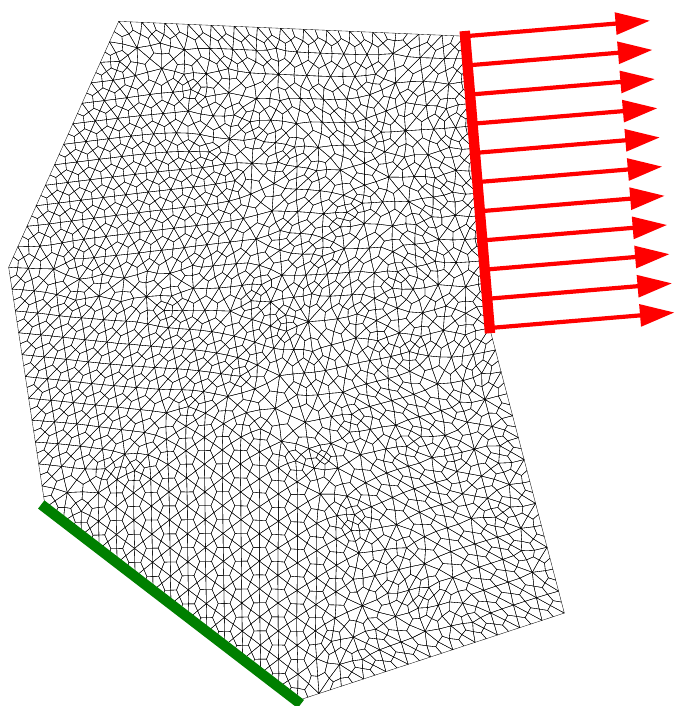}
		\caption{}
		\label{fig:meshNet_exp_uniform}
	\end{subfigure} 
	\begin{subfigure}[b]{0.22\textwidth}
		\includegraphics[width=1\textwidth]{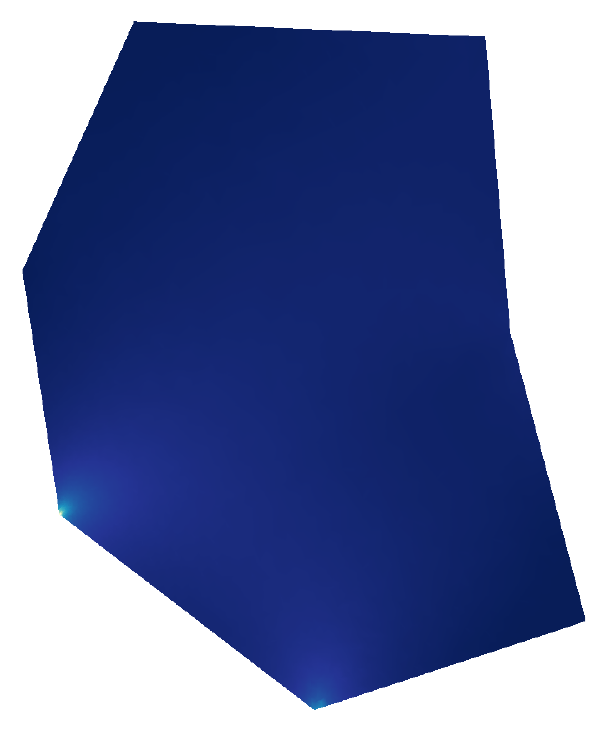}
		\caption{}
		\label{fig:meshNet_exp_adapt_VM}
	\end{subfigure} \quad\quad
	\begin{subfigure}[b]{0.225\textwidth}
		\includegraphics[width=1\textwidth]{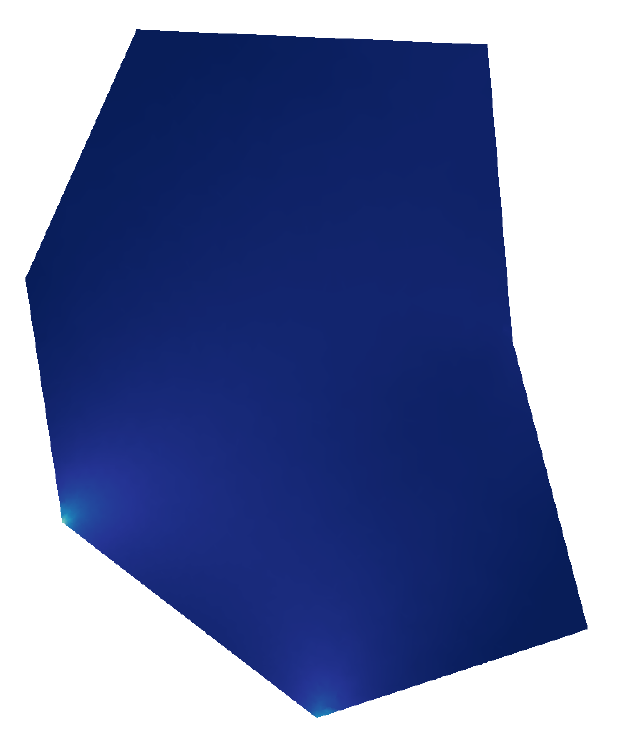}
		\caption{}
		\label{fig:meshNet_exp_output_VM}
	\end{subfigure} \quad\quad
	\begin{subfigure}[b]{0.325\textwidth}
		\includegraphics[width=1\textwidth]{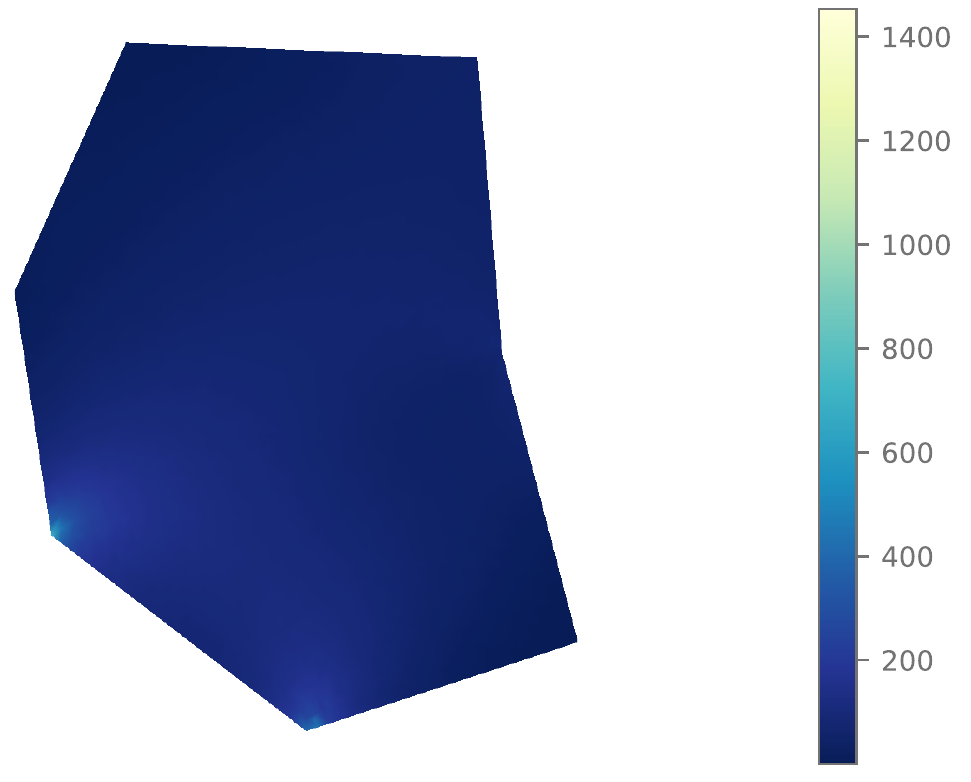}
		\caption{}
		\label{fig:meshNet_exp_uniform_VM}
	\end{subfigure}
	\caption{{Geometry from \cite{zhang2020meshingnet}, with the lower left edge of the domain set as the Dirichlet boundary. First row: mesh from (a) adaptive refinement, (b) predicted output, and (c) uniform refinement. Second row: von Mises stress for (d) adaptive mesh, (e) predicted output mesh, and (f) uniform refinement mesh.}}
	\label{fig:meshNet_exp}
\end{figure}

\begin{figure}[hbt!]
	\centering
	\begin{subfigure}[b]{0.25\textwidth}
		\includegraphics[width=1\textwidth]{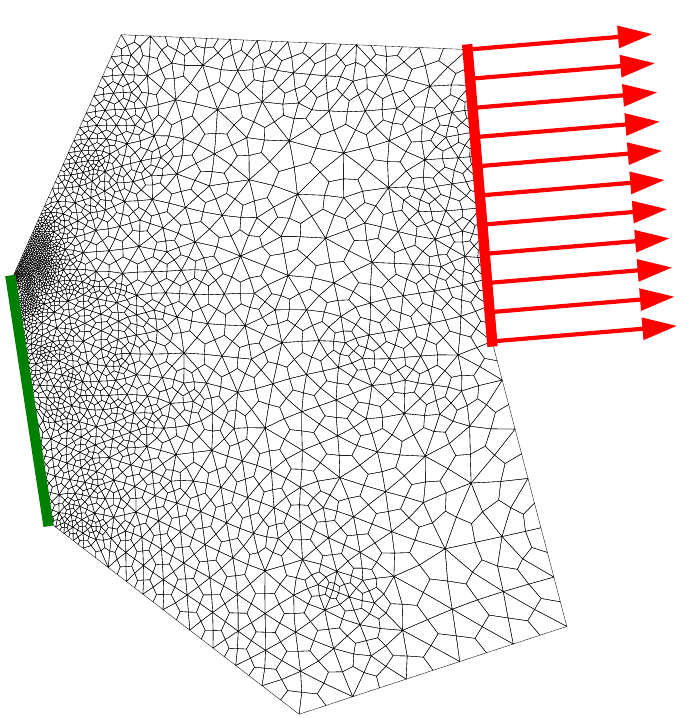}
		\caption{}
		\label{fig:meshNet_exp_adapt_b}
	\end{subfigure}  \quad\quad\quad
	\begin{subfigure}[b]{0.25\textwidth}
		\includegraphics[width=1\textwidth]{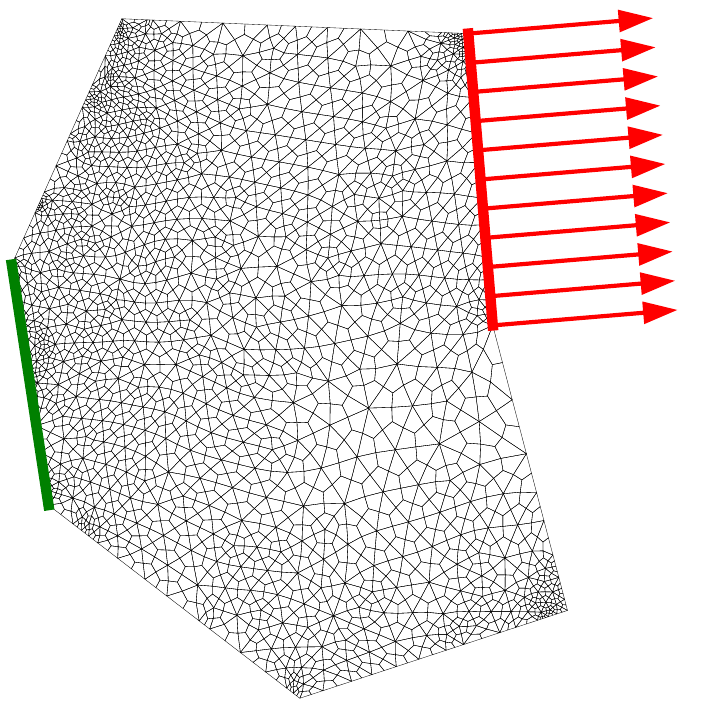}
		\caption{}
		\label{fig:meshNet_exp_output_b}
	\end{subfigure}  \quad\quad\quad
	\begin{subfigure}[b]{0.25\textwidth}
		\includegraphics[width=1\textwidth]{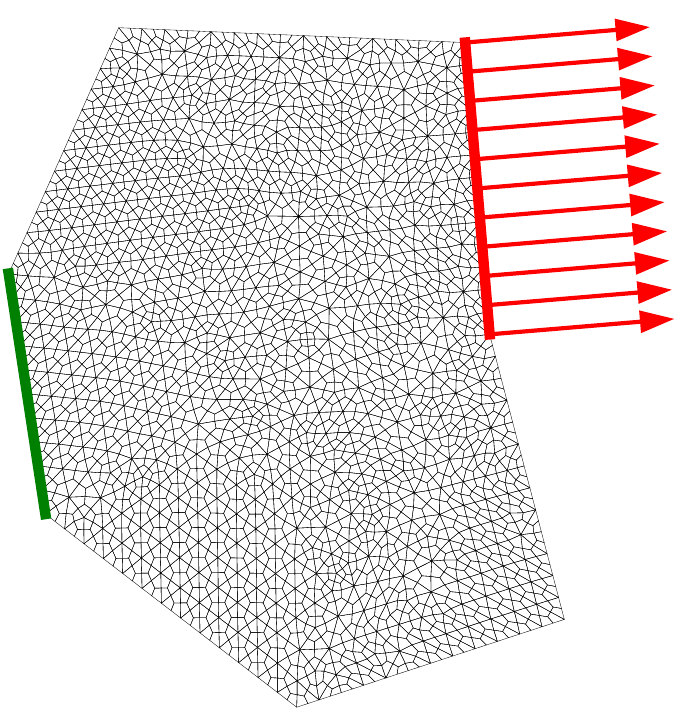}
		\caption{}
		\label{fig:meshNet_exp_uniform_b}
	\end{subfigure} 
	\begin{subfigure}[b]{0.235\textwidth}
		\includegraphics[width=1\textwidth]{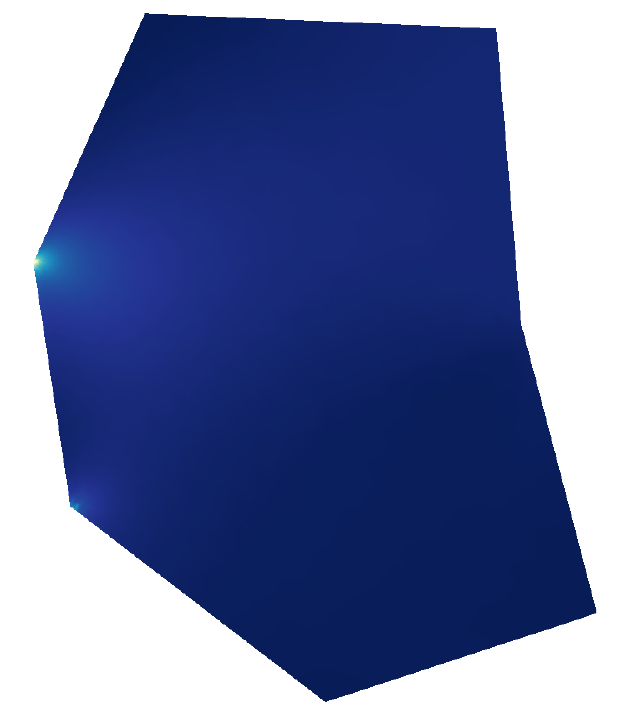}
		\caption{}
		\label{fig:meshNet_exp_adapt_VM_b}
	\end{subfigure} \quad\quad
	\begin{subfigure}[b]{0.225\textwidth}
		\includegraphics[width=1\textwidth]{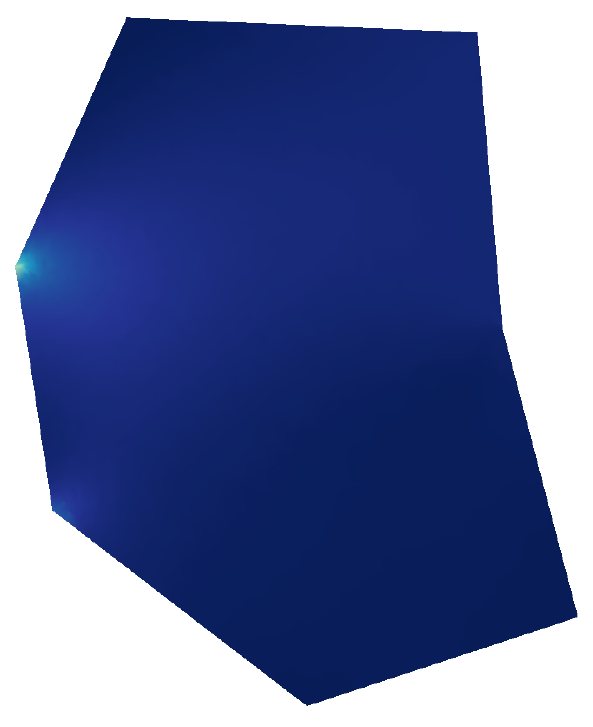}
		\caption{}
		\label{fig:meshNet_exp_output_VM_b}
	\end{subfigure} \quad\quad
	\begin{subfigure}[b]{0.355\textwidth}
		\includegraphics[width=1\textwidth]{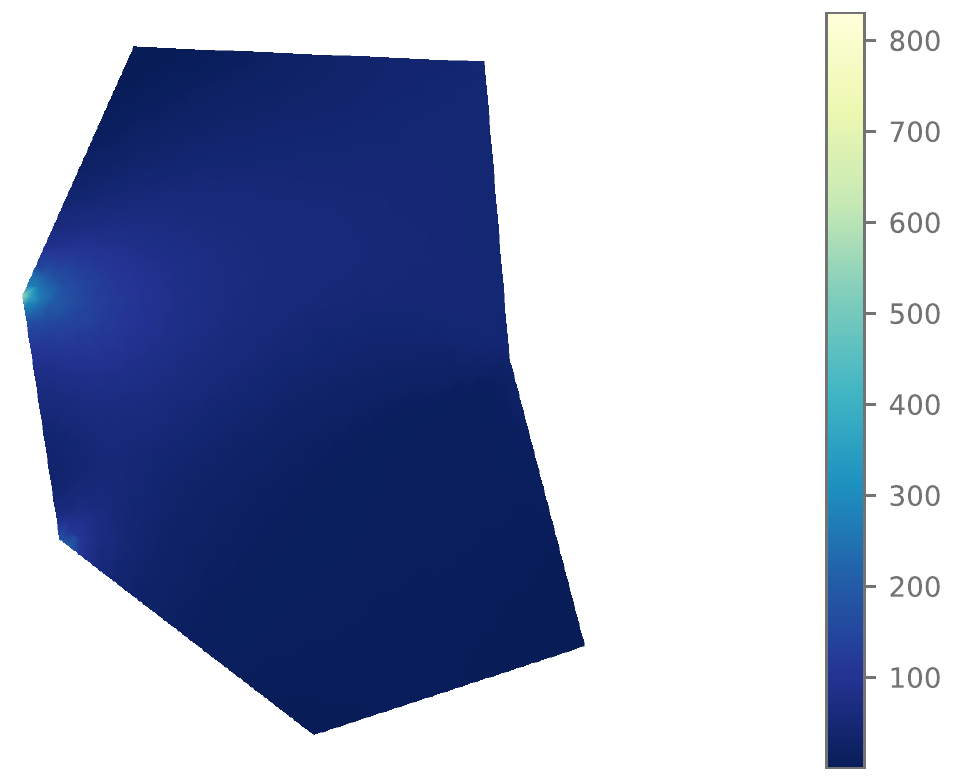}
		\caption{}
		\label{fig:meshNet_exp_uniform_VM_b}
	\end{subfigure}
	\caption{{Geometry from \cite{zhang2020meshingnet}, with the edge at the left side of the domain set as the Dirichlet boundary. First row: mesh from (a) adaptive refinement, (b) predicted output, and (c) uniform refinement. Second row: von Mises stress for (d) adaptive mesh, (e) predicted output mesh, and (f) uniform refinement mesh.}}
	\label{fig:meshNet_exp_b}
\end{figure}

\newpage
{
\subsection{Comparison with MeshingNet on a test geometry}\label{MeshNet}
 }
{In the following, we demonstrate the proposed method's capabilities for handling geometries obtained from image data. We compare our approach with the approach presented in~\cite{zhang2020meshingnet}. We present the results of our method in Figures~\ref{fig:meshNet_exp} and~\ref{fig:meshNet_exp_b}, which are both obtained from the image of the domain given in Figure~7 of~\cite{zhang2020meshingnet}. Note that the linear elasticity problems in~\cite{zhang2020meshingnet} are set in such a way that the Dirichlet boundary is formed by two edges and the Neumann boundary is formed by one edge. Since our models are trained with one edge for each of the Dirichlet and Neumann boundaries, we set the boundary conditions in the following examples to be compatible with the trained  models.}

{Figures \ref{fig:meshNet_exp} and \ref{fig:meshNet_exp_b} show two examples where each is assigned one of the Dirichlet edges given in \cite{zhang2020meshingnet}; whereas the Neumann boundary is set at the same position. In the two figures, the Dirichlet and Neumann boundaries are highlighted with green and red, respectively. We can observe that both predicted meshes (Figures \ref{fig:meshNet_exp_output} and \ref{fig:meshNet_exp_output_b}) are finer around the Dirichlet boundary. The locations where stress values are expected to be high can be predicted accordingly and the analysis results from the output meshes have lower relative energy norms compared to uniform meshes. Tables \ref{tab:meshNet_exp} and \ref{tab:meshNet_exp_b} show the comparison of the number of elements for each mesh and the relative energy norm computed by using the adaptive mesh as the reference solution. We point out that the predicted mesh might not match perfectly with the adaptive refinement mesh. This is because the latter is computed using exact geometry, while the predicted local refinement is based on coarse resolution images.}
 
{This study and \cite{zhang2020meshingnet} both introduce machine learning approaches for meshing problems. However, we focus on different aspects. In \cite{zhang2020meshingnet}, the domain is restricted to be a polygon with 6-8 edges. In addition, the Dirichlet boundaries are fixed at the 4th and 5th edge, while the 1st edge is always set as the Neumann boundary. But, the model has the option of setting different (homogeneous) material properties and traction with random amplitude up to 1000. Since the underlying problem is linear, we assume that these inputs do not significantly affect the relative mesh density. However, the optimal mesh grading may be affected. Thus, to improve the results we also discuss extensions using varying PDE parameters in Section~\ref{futurework}. Our method emphasizes handling a larger variety of geometries ranging from convex to non-convex and domains with voids. We restrict our model to impose boundary conditions only on non-adjacent edges to avoid the singularity formed at the points where the Dirichlet and Neumann boundaries join.  Nonetheless, the proposed method can also be extended to work with other boundary conditions by simply adding new data in the training process.}
\begin{table}[hbt!]
	\centering
	\begin{tabular}[t]{lccc}
		\hline
		Refinement & Number of elements& Relative error in the energy norm & Maximum value of von Mises stress\\
		\hline
		adaptive & 3324 & - &1453.0795\\
		predicted & 3294 &  0.047580 & 882.712778\\
		uniform & 3330 &  0.069678 & 670.116392\\
		\hline		
	\end{tabular}
	\caption{{Comparison between the number of elements and relative errors in the energy norm, for the example from Figure~\ref{fig:meshNet_exp}, on a geometry from~\cite{zhang2020meshingnet}.}}
	\label{tab:meshNet_exp}
\end{table}
\begin{table}[hbt!]
	\centering
	\begin{tabular}[t]{lccc}
		\hline
		Refinement & Number of elements& Relative error in the energy norm & Maximum value of von Mises stress\\
		\hline
		adaptive & 3153 & - & 829.819385\\
		predicted & 3147 &  0.067127 & 570.675438\\
		uniform & 3138 &  0.090880 & 539.449472\\
		\hline		
	\end{tabular}
	\caption{{Comparison between the number of elements and relative errors in the energy norm, for the example from Figure~\ref{fig:meshNet_exp_b}, on a geometry from~\cite{zhang2020meshingnet}.}}
	\label{tab:meshNet_exp_b}
\end{table}

\section{Discussion and future work}\label{futurework}
A key strength of our method is that it can be extended in a vast number of possible directions.
These include the application to more complicated model problems as well as the use of different discretization techniques for solving partial differential equations.

{There are two possible ways of extending our approach: Either, we generate training data for more general cases or we modify the network architecture itself, summarized in Sections~\ref{sec:extension-training} and~\ref{sec:extension-network}, respectively.}

{
\subsection{Extension by adding training data}\label{sec:extension-training}
{The measures of geometric complexity that we consider in Section~\ref{complexityGeometry} can be extended in many ways by generating more training data sets.
For example, in Section~\ref{sec:mesh_generation} we restricted our data set to contain only problems where the fixed and the traction boundary are not adjacent. However, this restriction may be dropped by properly extending the training data set.
Similarly, one may train with several edges marked as fixed and/or several edges marked as traction boundary.
}
}
\subsubsection{Extension to curved boundaries}\label{sec:exp_011}

\begin{figure}[hbt!]
	\centering	
	\begin{subfigure}[b]{0.3\textwidth}
		\includegraphics[width=1\textwidth]{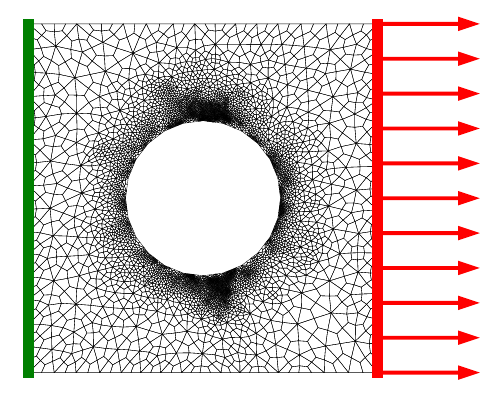}
		\caption{}
		\label{fig:square_circular_hole_adaptive}
	\end{subfigure} 
	\begin{subfigure}[b]{0.3\textwidth}
		\includegraphics[width=1\textwidth]{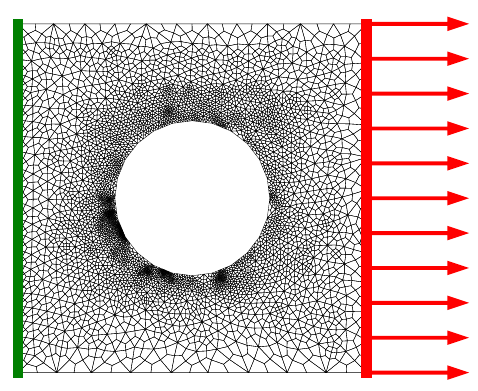}
		\caption{}
		\label{fig:square_circular_hole_output}
	\end{subfigure} 
	\begin{subfigure}[b]{0.3\textwidth}
		\includegraphics[width=1\textwidth]{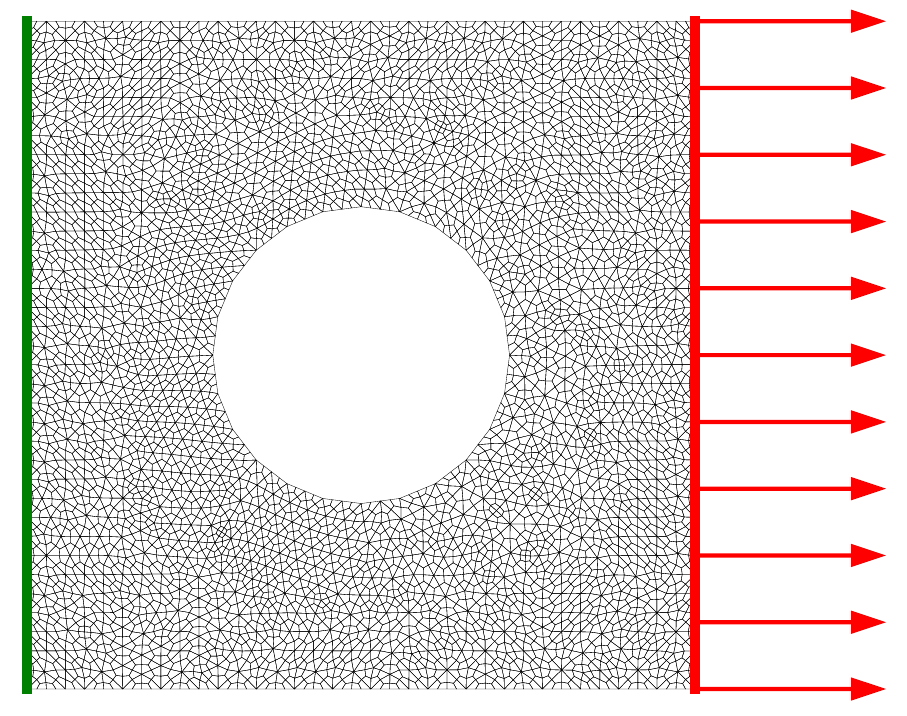}
		\caption{}
		\label{fig:square_circular_hole_uniform}
	\end{subfigure} 
	\begin{subfigure}[b]{0.225\textwidth}
		\includegraphics[width=1\textwidth]{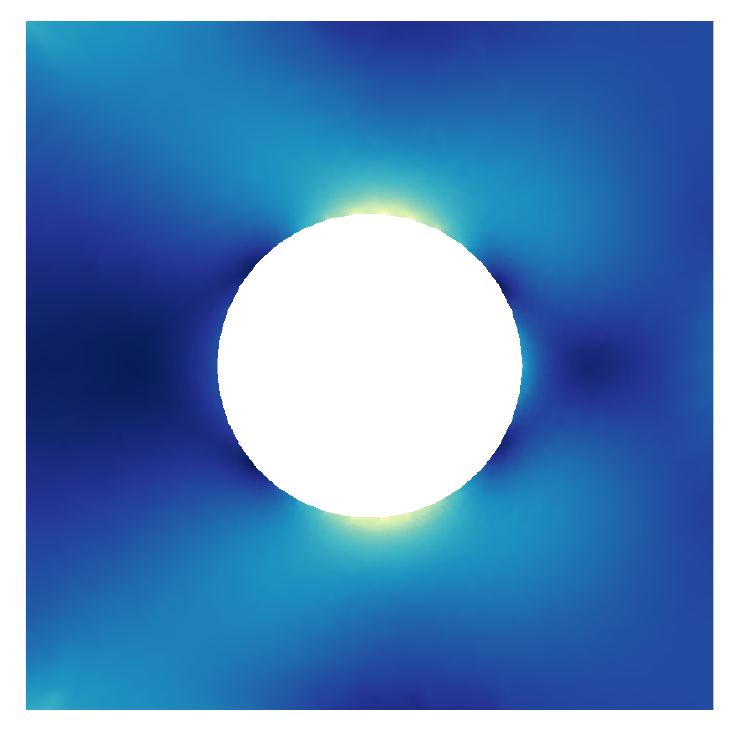}
		\caption{}
		\label{fig:square_circular_hole_adaptiveVM}
	\end{subfigure} \quad\quad
	\begin{subfigure}[b]{0.225\textwidth}
		\includegraphics[width=1\textwidth]{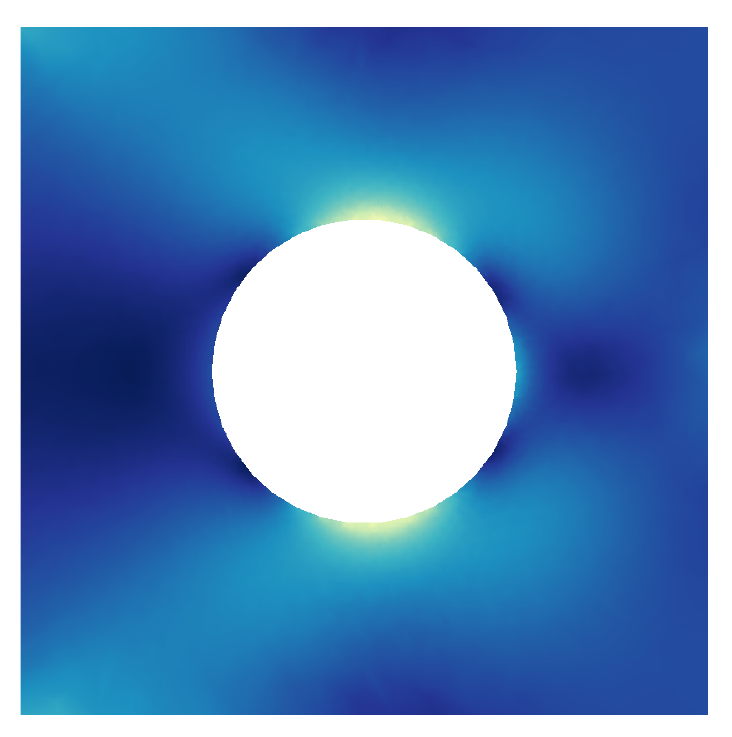}
		\caption{}
		\label{fig:square_circular_hole_outputVM}
	\end{subfigure} \quad\quad
	\begin{subfigure}[b]{0.31\textwidth}
		\includegraphics[width=1\textwidth]{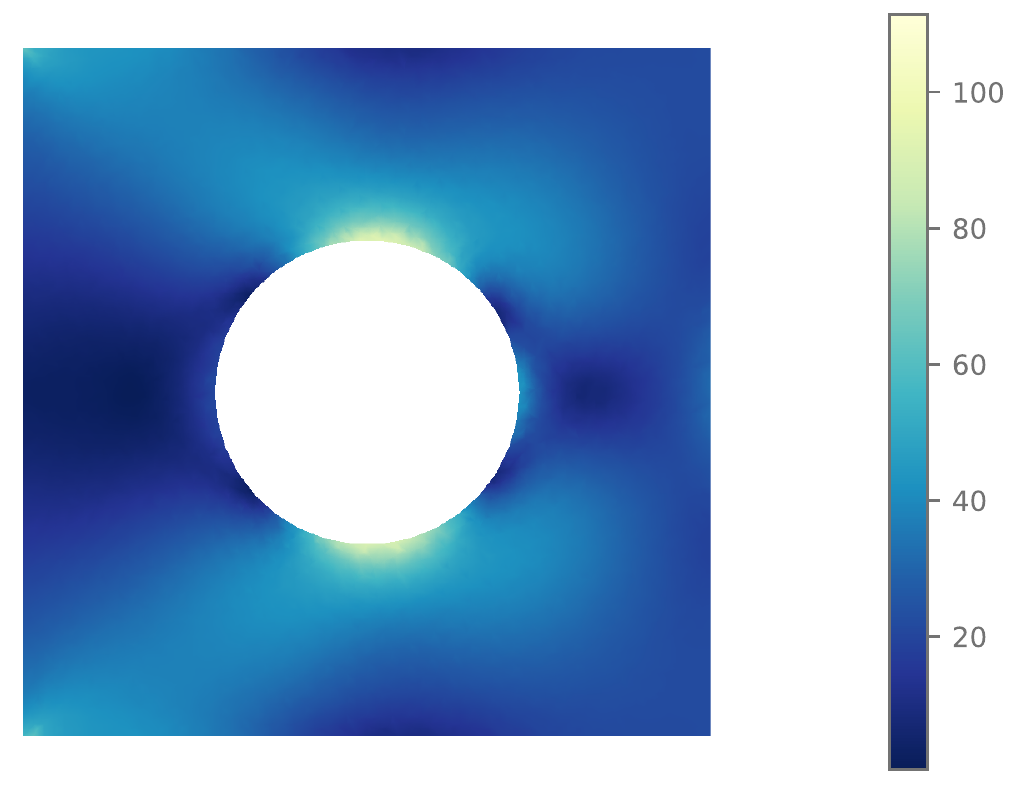}
		\caption{}
		\label{fig:square_circular_hole_uniformVM}
	\end{subfigure}
	\begin{subfigure}[b]{0.235\textwidth}
		\includegraphics[width=1\textwidth]{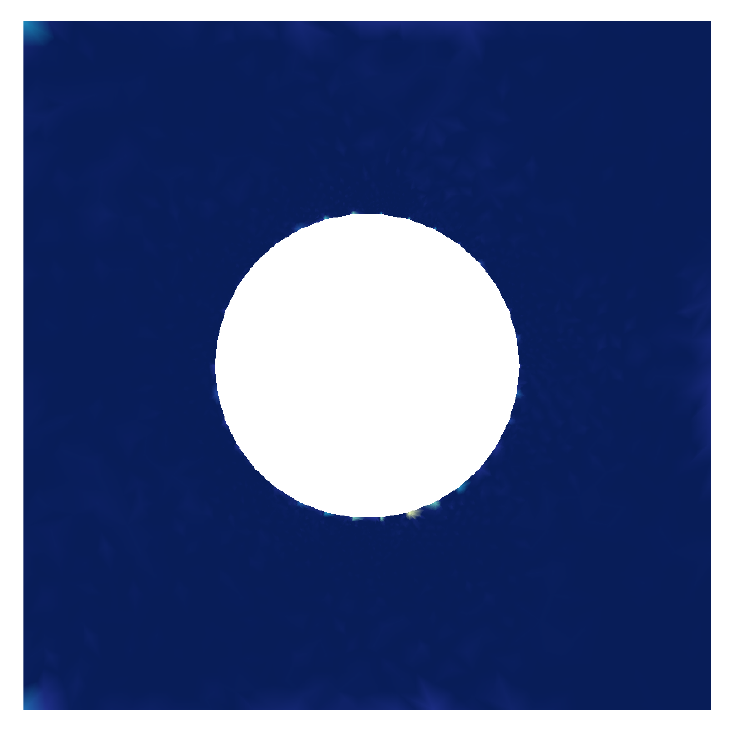}
		\caption{}
		\label{fig:square_circular_hole_VM_error_adaptive}
	\end{subfigure} \quad\quad
	\begin{subfigure}[b]{0.235\textwidth}
		\includegraphics[width=1\textwidth]{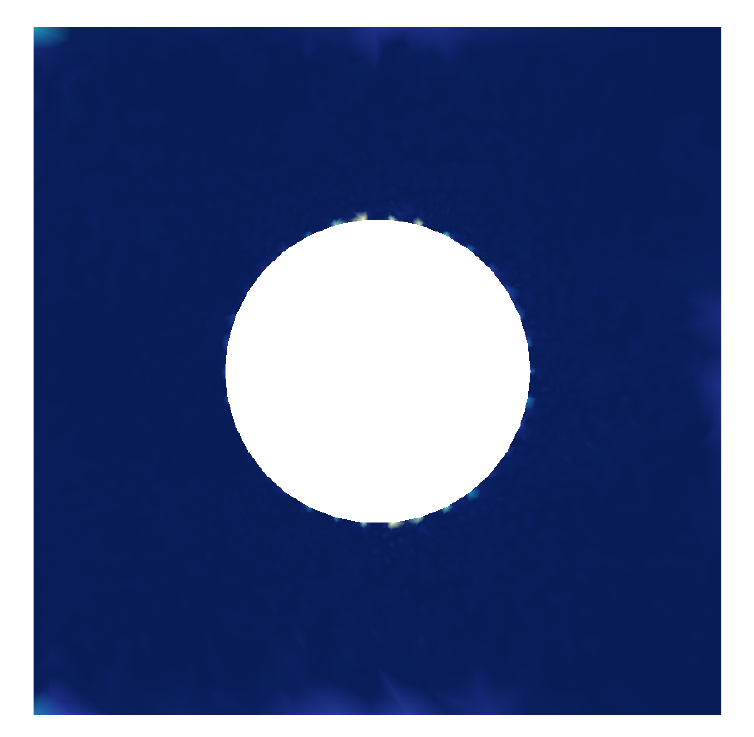}
		\caption{}
		\label{fig:square_circular_hole_VM_error_output}
	\end{subfigure} \quad\quad
	\begin{subfigure}[b]{0.3\textwidth}
		\includegraphics[width=1\textwidth]{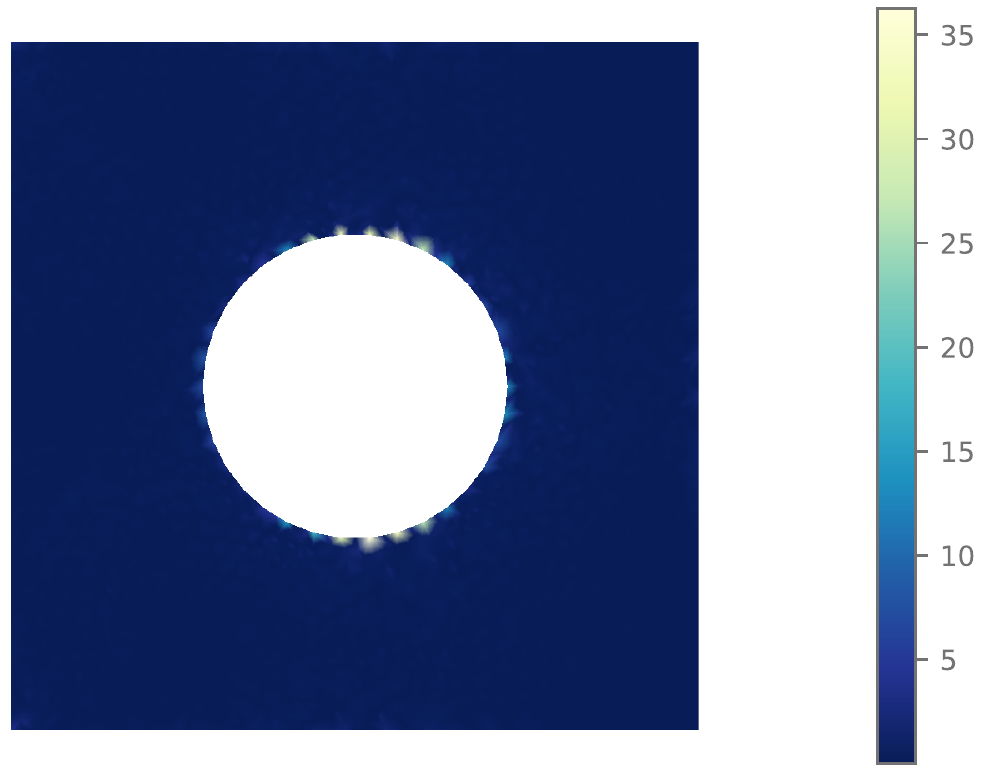}
		\caption{}
		\label{fig:square_circular_hole_VM_error_uniform_lc}
	\end{subfigure}
	\caption{Square with circular hole example. First row: (a) adaptive mesh, (b) output mesh, and (c) uniform refinement mesh. Second row: von Mises stress computed on (d) adaptive mesh, (e) output mesh and (f) uniform refinement mesh. Third row: von Mises stress difference between reference solution with (g) adaptive mesh, (h) output mesh and (i) uniform refinement mesh.}
	\label{fig:square_circular_hole}
\end{figure}

Another possible extension is the application to domains with curved boundaries. To visualize, we present an example for geometric complexity class (0,1,1). We use a well-known benchmark problem (the square plate with circular hole) to demonstrate this generalization of the proposed data driven method. As mentioned in Section \ref{sec:mesh_generation}, the data used for training are defined by polygons. Smooth boundaries such as circular voids are 'new' to the trained model. We observe that the model predicts a larger mesh size for such a type of geometry that it has not encountered before; this leads to a higher relative error $\varepsilon_{rel}=0.08481937$. The mesh constructed using the predicted output has very similar mesh size distribution pattern as the adaptive mesh (i.e. a finer mesh around the hole). While the relative error in the energy norm from the predicted mesh and uniform mesh are very close (see Table \ref{tab:square_energyNorm}); it is shown in Figure \ref{fig:square_circular_hole_VM_error_output} that the predicted output mesh has better von Mises stress approximation around the hole compared to a uniform mesh (see Figure \ref{fig:square_circular_hole_VM_error_uniform_lc}).  

\begin{table}[hbt!]
	\centering
	\begin{tabular}[t]{lccc}
		\hline
		Refinement & Number of elements& Relative error in energy norm & Maximum value of von Mises stress\\
		\hline
		adaptive &  6634 & 0.023551 & 111.425341\\
		predicted & 6424 & 0.026283 & 98.606861\\
		uniform & 6536 & 0.026432 & 92.712228\\
		\hline		
	\end{tabular}
	\caption{{Example for geometric complexity class $(0,1,1)$: number of elements, relative error in the energy norm, and maximum value of von Mises stress.}}
	\label{tab:square_energyNorm}
\end{table}%

On this example one can see that the network can be applied to unseen data, such as curved boundaries, and that the output mesh density is suitable, but not necessarily as good as an adaptive fit. Visually, the output density distribution is similar to the adaptive mesh used as training data, but more uniform. This effect may be reduced by post-processing the output image. Moreover, including domains with curved boundaries also in the training data for the network will most likely improve the result.

\subsubsection{Beyond finite element methods}
Most importantly, we intend to extend our method to discretization techniques beyond finite element methods, such as isogeometric analysis, which was introduced in~\cite{hughes2005isogeometric}. While the network can be trained on data generated from any kind of adaptive numerical method, it can then be applied to any space of locally refinable discretizations, such as linear finite elements, multi-patch splines, THB-splines~\cite{giannelli2012thb} or (unstructured) T-splines~\cite{sederberg2003t,casquero2020seamless}.
Since B-splines (and derived concepts) rely on a patch structure, where each patch is tensor-product, refinement becomes more involved. Each patch may be refined in a tensor-product fashion or using e.g. hierarchical B-splines or T-splines. While the latter two are more flexible, a fully tensor-product refinement for each patch is significantly easier to generate. Moreover, the subsequent assembly and solution of the resulting linear system is also faster.
Thus, the multi-patch segmentation should be selected a-priori such that a tensor-product refinement is as efficient (with respect to the number of degrees of freedom) as a local refinement scheme using e.g. THB-splines. Hence, it may be of interest to encode not only the local mesh size, but to also obtain information on the mesh anisotropy. If in a region of the domain the local refinement direction is clear, one may fit a patch such that its parameter lines correspond to that direction.
Such a segmentation may be performed using classical methods~\cite{pauley2015isogeometric, li2021deep} or using machine learning as well.

{
\subsubsection{Improved training strategies}\label{improvedtraining}
The framework of our approach is given in such a way that, once the family of PDE problems (including material parameters, possible boundary conditions etc.) is set up and the image dimensions are chosen, the neural network can be trained based on the requirements of the users. We envision a continuous training strategy, always improving the network by including additional model geometries (or families of geometries) based on user experience and user demand. This way, the network architecture does not need to change when new types of input data are considered.

Another way to continuously improve the output of the network is by generating training data obtained from improved adaptive refinement schemes. Currently, we use a heuristic scheme, presented in Section~\ref{generating-target-density}, which may be replaced by a (quasi-)optimal local refinement scheme.
}
{
\subsection{Generalizations that require a modified network structure}\label{sec:extension-network}
}
{
In the following we discuss possible generalizations that require a new network structure, increasing the number of input/output layers or the changing their dimensions. In some cases, the previously trained networks may be modified, in other cases a new network may need to be set up.}

{

\subsubsection{More general model problems}\label{generalmodel}
In this paper we restrict ourselves to linear elasticity on a planar domain without body forces and constant boundary conditions along the traction boundary and fixed boundary.
While this is already a challenging problem, our method is by no means restricted to this case.

Since we focus on a representation of the domain as an image, it is possible to encode varying parameters of the PDE or of the boundary conditions as color (or grayscale) images.
Similarly, body forces can be included. Besides the variations in geometric complexity, it is then also of interest to vary the complexity of the PDE. 
As with the geometric models, one can define suitable complexity classes by hand, such as vanishing or non-vanishing body forces, constant or varying parameters etc.
Alternatively, one may use a library of problems with known (or desired) solutions as a guide.

\subsubsection{Pre- and post-processing of input and output}\label{prepost}

In the following we also want to highlight how the algorithm may be improved by properly processing the input and output. Since the network is set up such that both input and output are images, pre- and post-processing steps from image processing may be used to improve the results. Note that the network allows any image of the correct size as an input. Thus, the geometry does not need to be from one of the considered classes. This property can be used to train the network also with segments/cutouts of larger images, which may be split into appropriately sized parts. Moreover, image segments of different resolution may be considered and appropriately merged.

The output of the network may be smoothened to obtain meshes with improved mesh grading. Similarly, one may use edge/highlight detection algorithms to find features in the output image that need to be refined more. This approach is also proposed in Section~\ref{sec:exp_000} and tested, see e.g. Table~\ref{tab:convex_exp}. Since we consider images, it may also be feasible to train the network with noisy data to make the output more robust. However, such modifications must be performed carefully and studied more deeply, as small details in the geometry have a significant effect on the expected local mesh size.

}

{
\subsubsection{Higher dimensions}\label{higherdim}
The proposed method can be generalized to problems of higher dimensions  such as solving PDE on volumetric domains or solving time dependent problems such as parabolic PDE using space-time methods.
While the method can be set up for problems on volumetric domains in a straight-forward way, an extension to surface domains or space-time problems can be more difficult.
The extension to space-time domains leads to 3D or multi-channel images, instead of just grayscale. 
Note that for time-dependent PDEs, other types of networks like recurrent neural networks or LSTM might be more suitable.

{Multi-channel images may also be used to generalize the approach to surface data when the surface can be parameterized over a planar quadrilateral domain.
In this case, the differential equation on the parameterized surface can be pulled back to the parametric domain, resulting in a PDE with varying coefficients that can be encoded as separate channels of the image data.
A typical example for this case are single-patch B-spline surfaces that are employed in isogeometric analysis.}

In order to process surface data with more complicated methods from geometric deep learning~\cite{bronstein2017geometric} may be used for surface data.
These methods generalize the notion of image CNN to data that is represented by graphs or discrete manifolds.
The challenge is here to define operators that act on these data types and that perform similarly to the convolution operator of matrices.
CNN for discrete for discrete surface data are an active field of research \cite{fey2018splinecnn, 10.1145/3326362} and it would be an interesting direction for future work to employ these techniques to generate locally refined surface parameterizations.
}
\section{Conclusions}\label{conclusions}

In this paper we developed a local refinement method for PDEs using machine learning. The geometry of the domain as well as the boundary data of the PDE are encoded as images. These images are then used as input data for a neural network based on the U-net architecture. The output of the network is another image which encodes the local mesh size of a finite element mesh. From this information we then construct a mesh over the domain using Gmsh.
The quality of the locally refined mesh can then be measured by computing the discretization error of the resulting finite element solution.

We compare the quality of the output for different training strategies. To do so, we categorize possible input data sets using different geometric complexity classes. One can then pursue different training strategies. It is expected that the resulting output yields best results, if the network is trained exactly with those training data sets from the same complexity class as the input. This may however be quite expensive, since it means that one has to generate training data for a wide range of geometric complexities. We could show that it suffices to train the network with data sets which are complex only with respect to one dimension and simple with respect to other dimensions of complexity. For instance, if the network is trained for L-shaped domains without holes as well as domains with a convex boundary and a hole, it can produce a suitable output for an L-shaped domain with a hole, even though it has never seen such a domain. This means that one does not necessarily need to train for all possible geometric complexity classes but can restrict to a subset for which training data is easier to produce.

In the future we want to extend the approach to other types of PDEs and to a larger class of possible domains and data of the PDE. Most importantly, the approach can, in principle, also be extended to 3D. However, the classification of possible domains and the generation of training data becomes more involved. In addition, we want to extend the method also to generating multi-patch isogeometric discretizations with curved boundaries, which take into account the estimated mesh density.

\backmatter

\bmhead{Acknowledgments}

The authors were supported by the Linz Institute of Technology (LIT) and the government of Upper Austria through the project LIT-2019-8-SEE-116 entitled ``PARTITION -- PDE-aware isogeometric discretization based on neural networks''. This support is gratefully acknowledged.
F.~Scholz additionally acknowledges the support by the JST CREST Grant no. JPMJCR1911.

\end{document}